\documentclass[12pt,british]{article}

\usepackage{babel}
\usepackage{amsmath,amsfonts,amssymb,amsthm,amscd}
\usepackage{eucal,enumerate,textcomp,url}
\usepackage{tikz}
\usetikzlibrary{calc}
\usepackage[matrix,arrow]{xy}

\usepackage{ucs} 
\usepackage[utf8x]{inputenc} 
\usepackage[T1]{fontenc} 

\usepackage[noBBpl]{mathpazo}
\usepackage{mathabx}
\let\newbigast=\bigast
\renewcommand{\bigast}{\mathbin{\newbigast}}

\usepackage{braket}
\changenotsign

\usepackage[pdftex,                %
breaklinks=true,
bookmarksnumbered = true,
colorlinks= true,
urlcolor= green,
anchorcolor = yellow,
citecolor=blue,
pdfborder= {0 0 0}
]{hyperref}                   

\makeatletter

\renewcommand{\p@enumii}{}
\renewcommand{\p@enumiii}{}

\def\@enum@{\list{\csname label\@enumctr\endcsname}%
           {\usecounter{\@enumctr}\def\makelabel##1{
\normalfont\ignorespaces\emph{{##1}~}}
\setlength{\labelsep}{3pt}
\setlength{\parsep}{0pt}
\setlength{\itemsep}{0pt}
\setlength{\leftmargin}{0pt}
\setlength{\labelwidth}{0pt}
\setlength{\listparindent}{\parindent}
\setlength{\itemsep}{0pt}
\setlength{\itemindent}{0pt}
\topsep=3pt plus 1pt minus 1 pt}}

\def\@map#1#2[#3]{\mbox{$#1 \colon\thinspace #2 \to #3$}}
\def\map#1#2{\@ifnextchar [{\@map{#1}{#2}}{\@map{#1}{#2}[#2]}}
\makeatother

\newcommand{\altarrow}[4]{\mbox{$#1 \colon\thinspace #2 #3 #4$}}

\renewcommand{\epsilon}{\ensuremath{\varepsilon}}
\renewcommand{\phi}{\ensuremath{\varphi}}
\newcommand{\vide}{\ensuremath{\varnothing}}
\renewcommand{\to}{\ensuremath{\longrightarrow}}
\renewcommand{\mapsto}{\ensuremath{\longmapsto}}

\DeclareMathOperator{\id}{\text{Id}}

\renewcommand{\ker}[1]{\ensuremath{\operatorname{\text{Ker}}\left({#1}\right)}}

\newcommand{\rp}{\ensuremath{\mathbb{R}P^2}}
\newcommand{\dt}{\ensuremath{\mathbb D}^{2}}
\newcommand{\Z}{\ensuremath{\mathbb Z}}
\newcommand{\dbK}{\ensuremath{\mathbb K}}
\newcommand{\N}{\ensuremath{\mathbb N}}

\newcommand{\C}{\ensuremath{\mathbb C}}
\newcommand{\R}{\ensuremath{\mathbb R}}

\newcommand{\St}[1][2]{\ensuremath{\mathbb S}^{#1}}

\newcommand{\FF}{\ensuremath{\mathbb F}}
\newcommand{\F}[1][n]{\ensuremath{\FF_{{#1}}}}

\renewcommand{\to}{\ensuremath{\longrightarrow}}
\renewcommand{\ker}[1]{\ensuremath{\operatorname{\text{Ker}}\left({#1}\right)}}

\DeclareRobustCommand*{\up}[1]{\textsuperscript{#1}}
\newcommand{\ft}[1][n]{\ensuremath{{\Delta_{#1}^2}}}

\newcommand{\out}[1]{\ensuremath{\operatorname{\text{Out}}\left({#1}\right)}}
\newcommand{\aut}[1]{\ensuremath{\operatorname{\text{Aut}}\left({#1}\right)}}
\newcommand{\sn}[1][n]{\ensuremath{S_{#1}}}
\newcommand{\an}[1][n]{\ensuremath{A_{#1}}}

\newcommand{\hooklongrightarrow}{\lhook\joinrel\longrightarrow}

\newcommand{\lhra}{\mathrel{\lhook\joinrel\to}}

\newcommand{\ang}[1]{\ensuremath{\left\langle #1\right\rangle}}
\newcommand{\normcl}[1]{\ensuremath{\ang{\!\ang{#1}\!}}}
\newcommand{\setangr}[2]{\ensuremath{\ang{#1 \,\left\lvert \, #2 \right.}}}

\newcommand{\setr}[2]{\ensuremath{\brak{#1 \,\left\lvert \, #2 \right.}}}
\newcommand{\setl}[2]{\ensuremath{\brak{\left. #1 \,\right\rvert \, #2}}}
\newcommand{\im}[1]{\ensuremath{\operatorname{Im}(#1)}}
\newcommand{\Int}[1]{\ensuremath{\operatorname{Int}(#1)}}

\renewcommand{\epsilon}{\varepsilon}
\renewcommand{\th}{\ensuremath{\up{th}}}

\newcommand{\tstar}{\ensuremath{\operatorname{T}^{\ast}}}

\newcommand{\istar}{\ensuremath{\operatorname{I}^{\ast}}}
\newcommand{\dic}[1]{\ensuremath{\operatorname{\text{Dic}}_{#1}}}
\newcommand{\dih}[1]{\ensuremath{\operatorname{\text{Dih}}_{#1}}}
\newcommand{\quat}[1][8]{\ensuremath{\mathcal{Q}_{#1}}}
\newcommand{\mcggen}[2][n]{\ensuremath{\operatorname{\mathcal{MCG}}(#2,#1)}}
\newcommand{\mcgzero}[1]{\ensuremath{\operatorname{\mathcal{MCG}}(#1)}}

\newcommand{\mcg}[1][n]{\ensuremath{\operatorname{\mathcal{MCG}}(\St,#1)}}

\newcommand{\tonestar}{\ensuremath{T^{\ast}}}
\newcommand{\oonestar}{\ensuremath{O^{\ast}}}

\newcommand{\garside}[1][n]{\ensuremath{\Delta_{#1}}}

\DeclareRobustCommand*{\up}[1]{\textsuperscript{#1}}
\renewcommand{\th}{\ensuremath{\up{th}}}

\newcommand{\brak}[1]{\ensuremath{\left\{ #1 \right\}}}
\renewcommand{\set}[2]{\ensuremath{\Set{ #1 \, | \, #2}}}

\newcommand{\evc}[1]{\ensuremath{\underline{\underline E} #1}}
\newcommand{\bvc}[1]{\ensuremath{\underline{\underline B} #1}}

\theoremstyle{plain}
\newtheorem{thm}{Theorem}
\newtheorem{lem}[thm]{Lemma}
\newtheorem{prop}[thm]{Proposition}
\newtheorem{cor}[thm]{Corollary}

\newtheorem*{fjconjecture}{Isomorphism Conjecture (IC)}
\newtheorem*{ficconj}{Fibred Isomorphism Conjecture (FIC)}

\newtheoremstyle{newremark}
  {}
  {}
  {}
  {}
  {\bfseries}
  {.}
  {.5em}
  {}

\theoremstyle{newremark}
\newtheorem{rem}[thm]{Remark}

\newtheorem*{defn}{Definition}
\newtheorem{rems}[thm]{Remarks}

\newcommand{\reth}[1]{Theorem~\protect\ref{th:#1}}
\newcommand{\relem}[1]{Lemma~\protect\ref{lem:#1}}
\newcommand{\repr}[1]{Proposition~\protect\ref{prop:#1}}
\newcommand{\reco}[1]{Corollary~\protect\ref{cor:#1}}
\newcommand{\resec}[1]{Section~\protect\ref{sec:#1}}
\newcommand{\rerem}[1]{Remark~\protect\ref{rem:#1}}

\newcommand{\rerems}[1]{Remarks~\protect\ref{rem:#1}}
\newcommand{\req}[1]{equation~(\protect\ref{eq:#1})}
\newcommand{\reqref}[1]{(\protect\ref{eq:#1})}

\begin{document}

\title{A survey of surface braid groups and the lower algebraic $K$-theory of their group rings}

\author{John~Guaschi\vspace*{1mm}\\ 
Normandie Universit\'e, UNICAEN,\\
Laboratoire de Math\'ematiques Nicolas Oresme UMR CNRS~\textup{6139},\\
14032 Caen Cedex, France\\
e-mail:~\texttt{john.guaschi@unicaen.fr}\vspace*{4mm}\\
Daniel Juan-Pineda\vspace*{1mm}\\
Centro de Ciencias Matem\'aticas,\\
Universidad Nacional Aut\'onoma de M\'exico,
Campus Morelia,\\
Morelia, Michoac\'an, M\'exico 58089\\
e-mail:~\texttt{daniel@matmor.unam.mx}\vspace*{4mm}}

\date{5th February 2013}

\maketitle

\begingroup
\renewcommand{\thefootnote}{}
\endgroup 

\begin{abstract}
\noindent
\emph{In this article, we give a survey of the theory of surface braid groups and the lower algebraic $K$-theory of their group rings. We recall several definitions and describe various properties of surface braid groups, such as the existence of torsion, orderability, linearity, and their relation both with mapping class groups and with the homotopy groups of the $2$-sphere. The braid groups of the $2$-sphere and the real projective plane are of particular interest because they possess elements of finite order, and we discuss in detail their torsion and the classification of their finite and virtually cyclic subgroups. Finally, we outline the methods used to study the lower algebraic $K$-theory of the group rings of surface braid groups, highlighting recent results concerning the braid groups of the $2$-sphere and the real projective plane.}
\end{abstract}

\pagebreak

\tableofcontents



\section{Introduction}

The braid groups $B_n$ were introduced by E.~Artin in~1925~\cite{A1},
in a geometric and intuitive manner, and further studied in 1947 from
a more rigourous and algebraic standpoint~\cite{A2,A3}. These groups
may be considered as a geometric representation of the standard
everyday notion of braiding strings or strands of hair. As well as
being fascinating in their own right, braid groups play an
important rôle in many branches of mathematics, for example in
topology, geometry, algebra, dynamical systems and theoretical
physics, and notably in the study of knots and links~\cite{BZ}, in the
definition of topological invariants (Jones polynomial, Vassiliev
invariants)~\cite{Jon1,Jon2}, of the mapping class groups~\cite{Bi2,Bi3,FM}, and of
configuration spaces~\cite{CG,FH1}. They also have potential applications to
biology, robotics and cryptography, for example~\cite{BCHWW}. 

The Artin braid groups have been generalised in many different directions, such as Artin-Tits groups~\cite{Bri,BrS,Del}, surface braid groups, singular braid monoids and groups, and virtual and welded braid groups. One recent exciting topological development is the discovery of a connection between braid groups and the homotopy groups of the $2$-sphere via the notion of Brunnian braids~\cite{BCHWW,BCWW}. Although there are many surveys on braid groups~\cite{BB,GM4,Mag,MK,P,R,Ve} as well as some books and monographs~\cite{Bi2,Ha,KT,MK}, for the most part, the theory of surface braid groups is discussed in little detail in these works. The aim of this article is two-fold, the first being to survey various aspects of this theory and some recent results, highlighting the cases of the $2$-sphere and the real projective plane, and the second being to discuss current developments in the study of the lower algebraic $K$-theory of the group rings of surface braid groups. In \resec{basic}, we give various definitions of surface braid groups, and recall their relationship with mapping class groups. In \resec{properties}, we describe a number of properties of these groups, including the existence of Fadell-Neuwirth short exact sequences of their pure and mixed braid groups, which play a fundamental rôle in the theory. In \resec{present}, we recall some presentations of surface braid groups, and in Sections~\ref{sec:centre} and~\ref{sec:embeddings}, we survey known results about their centre and their embeddings in other braid groups. Within the theory of surface braid groups, those of the sphere $\St$ and the real projective plane $\rp$ are interesting and important, one reason being that their configuration spaces are not Eilenberg-Mac~Lane spaces. In \resec{homotype}, we study the homotopy type of these configuration spaces and the cohomological periodicity of the braid groups of $\St$ and $\rp$, and we describe some of the results mentioned above concerning Brunnian braids and the homotopy groups of $\St$. In Sections~\ref{sec:orderable} and~\ref{sec:linear}, we discuss orderability and linearity of surface braid groups.

\resec{virtually} is devoted to the study of the structure of the braid groups of $\St$ and $\rp$, notably their torsion, their finite subgroups and their virtually cyclic subgroups. Finally, in \resec{ktheory}, we discuss recent work concerning the $K$-theory of the group rings of surface braid groups. The existence of torsion in the braid groups of $\St$ and $\rp$ leads to new and interesting behaviour in the lower algebraic $K$-theory of their group rings. Recent techniques provided by the Fibred Isomorphism Conjecture (FIC) of Farrell and Jones have brought to light examples of of intricate group rings whose lower algebraic $K$-groups are trivial, see \reth{kthasp} for example, as well as highly-complicated algebraic $K$-theory groups. A fairly complete example of the latter is that of the $4$-string braid group $B_4(\St)$ of the sphere, for which we show that $K_i(\Z[B_4(\St)])$ is infinitely generated for $i=0,1$ (see \reth{b4s2kth}). We conjecture that a similar result is probably true for all $i>1$. On the other hand, it is known that $\operatorname{rank}(K_i(\Z[B_4(\St)]))<\infty$ for all $i\in \Z$~\cite{JS}. It is interesting to observe that the geometrical aspects of a group largely determine the structure of the algebraic $K$-groups of its group ring. We include up-to-date results on the algebraic $K$-groups of surface braid groups, and mention possible extensions of these computations. The main obstructions to extending our results from $B_{4}(\St)$ to the general case are the lack of appropriate models for their classifying spaces, as well the complicated structure of the Nil groups.

\enlargethispage{4mm}

\subsection*{Acknowledgements}

Both authors are grateful to the French-Mexican International Laboratory ``LAISLA'' for its financial support. The first author was partially supported by the international Cooperation Capes-Cofecub project numbers Ma~733-12 (France) and Cofecub~1716/2012 (Brazil). The second author would like to acknowledge funding from CONACyT and PAPIIT-UNAM. The first author wishes to thank Daciberg Lima Gon\c{c}alves for interesting and helpful conversations during the preparation of this paper.




\section{Basic definitions of surface braid groups}\label{sec:basic}

One of the interesting aspects about surface braid groups is that they may be defined from various viewpoints, each giving a different insight into their nature~\cite{R}. The notion of surface braid group was first introduced by Zariski, and generalises naturally Artin's geometric definition~\cite{Z1,Z2}. Surface braid groups were rediscovered during the 1960's by Fox who proposed a powerful (and equivalent) topological definition in terms of the fundamental group of configuration spaces. We recall these and other definitions below. Unless stated otherwise, in the whole of this manuscript, we shall use the word \emph{surface} to denote a connected surface, orientable or non orientable, with or without boundary, and of the form $M=N\setminus Y$, where $N$ is a compact, connected surface, and $Y$ is a finite (possibly empty) subset lying in the interior $\Int{N}$ of $N$.

\subsection{Surface braids as a collection of strings}\label{sec:defstring}

Let $M$ be a surface, and let $n\in \N$. We fix once and for all a finite $n$-point subset $X=\brak{x_{1},\ldots,x_{n}}$ of $\Int{M}$ whose elements shall be the base points of our braids.

\begin{figure}[h]
\hfill
\begin{tikzpicture}[scale=0.55,thick]
{\draw (-1,0) to[bend left] (1,0);
\draw (-1.2,.1) to[bend right] (1.2,.1);
\draw[rotate=0] (0,0) ellipse (100pt and 50pt);}

\draw (-1,6) to[bend left] (1,6);
\draw (-1.2,6.1) to[bend right] (1.2,6.1);


\node at (-6, 6) {$M\times \brak{0}$};
\node at (-6, 0) {$M\times \brak{1}$};

\node at (-1.9, 5.6) {$x_{1}$};
\node at (-0.65, 5.3) {$x_{2}$};
\node at (1.6, 5.6) {$x_{3}$};

\draw[draw=white,line width=5pt] (1.5,5) .. controls (1.5,0) and (0,2) .. (0,-1);
\draw[very thick] (1.5,5) .. controls (1.5,0) and (0,2) .. (0,-1);
\draw[draw=white,line width=5pt] (-1.5,5) .. controls (-1.5,0) and (1.5,4) .. (1.5,-1);
\draw[very thick] (-1.5,5) .. controls (-1.5,0) and (1.5,4) .. (1.5,-1);

\draw[draw=white,line width=5pt] (0,5) .. controls (0,0) and (-1.5,2) .. (-1.5,-1);
\draw[very thick] (0,5) .. controls (0,0) and (-1.5,2) .. (-1.5,-1);
\node at (-2.3, 3.5) {$\beta_{1}$};
\node at (-2.3, 0.2) {$\beta_{2}$};
\node at (2.3, 3) {$\beta_{3}$};

\draw[rotate=0,draw=white,line width=5pt] (0,6) ellipse (100pt and 50pt);
\draw[rotate=0] (0,6) ellipse (100pt and 50pt);

\foreach \k in {-1.5,0,1.5}
{\draw[fill] (\k,5) circle [radius=0.05];};

\foreach \k in {-1.5,0,1.5}
{\draw[fill] (\k,-1) circle [radius=0.05];};
\end{tikzpicture}
\hspace*{\fill}
\caption{A geometric $3$-braid with $M$ equal to the $2$-torus.}\label{fig:braid1}
\end{figure}
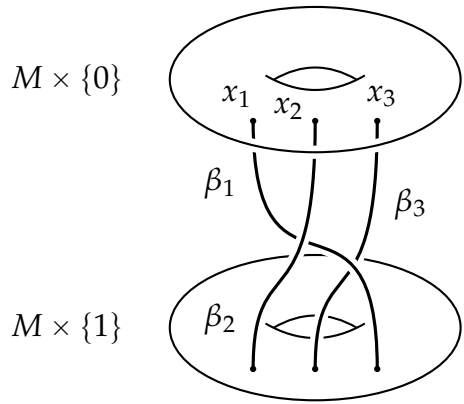

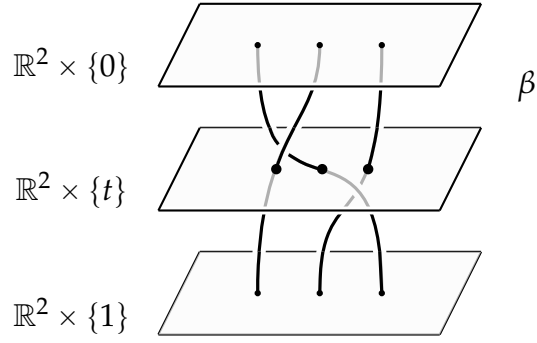
\begin{figure}[h]
\hfill
\begin{tikzpicture}[scale=0.55,thick]
\draw (-3.9,-2) -- (2.9,-2);
\draw (-2.9,0) -- (3.9,0);
\draw (-2.9,0) -- (-3.9,-2);
\draw (3.9,0) -- (2.9,-2);

\draw (-2.9,3) -- (3.9,3);

\draw[draw=white,line width=5pt] (1.5,5) .. controls (1.5,0) and (0,2) .. (0,-1);
\draw[very thick] (1.5,5) .. controls (1.5,0) and (0,2) .. (0,-1);
\draw[draw=white,line width=5pt] (-1.5,5) .. controls (-1.5,0) and (1.5,4) .. (1.5,-1);
\draw[very thick] (-1.5,5) .. controls (-1.5,0) and (1.5,4) .. (1.5,-1);

\draw[draw=white,line width=5pt] (0,5) .. controls (0,3) and (-1.5,3) .. (-1.5,-1);
\draw[very thick] (0,5) .. controls (0,3) and (-1.5,3) .. (-1.5,-1);

\path[fill=gray!2, opacity=0.7] (-3.9,4) -- (2.9,4) -- (3.9,6) -- (-2.9,6)-- cycle;
\path[fill=gray!2, opacity=0.7] (-3.9,1) -- (2.9,1) -- (3.9,3) -- (-2.9,3)-- cycle;
\path[fill=gray!2, opacity=0.7] (-3.9,-2) -- (2.9,-2) -- (3.9,0) -- (-2.9,0) -- cycle;

\draw[draw=white,line width=3pt] (-3.9,1) -- (2.9,1);
\draw (-3.9,1) -- (2.9,1);
\draw (-2.9,3) -- (-3.9,1);
\draw (3.9,3) -- (2.9,1);

\draw[line width=3pt,draw=white] (-3.9,4) -- (2.9,4);
\draw (-3.9,4) -- (2.9,4);
\draw (-2.9,6) -- (3.9,6);
\draw (-2.9,6) -- (-3.9,4);
\draw (3.9,6) -- (2.9,4);

\node at (-6,4.5) {$\R^2\times \brak{0}$};
\node at (-6,1.5) {$\R^2\times \brak{t}$};
\node at (-6,-1.55) {$\R^2\times \brak{1}$};

\draw[very thick,draw=white] (-1.5,-1) .. controls (-1.5,-0.5) and (-1.5,-0.5) .. (-1.46,0);
\draw[very thick] (-1.5,-1) .. controls (-1.5,-0.5) and (-1.5,-0.5) .. (-1.46,0);

\draw[very thick,draw=white] (0,-1) .. controls (0.02,-0.45) and (0.0,-0.53) .. (0.075,0);
\draw[very thick] (0,-1) .. controls (0.02,-0.45) and (0.0,-0.53) .. (0.075,0);

\draw[very thick] (1.5,-1) -- (1.45,0);

\draw[very thick,draw=white] (-1.06,2) .. controls (-0.88,2.5) and (-0.87,2.5) .. (-0.62,3);
\draw[very thick] (-1.06,2) .. controls (-0.88,2.5) and (-0.87,2.5) .. (-0.62,3);

\draw[very thick,draw=white] (-1,2.6) .. controls (-1.11,2.8) and (-1.125,2.77) .. (-1.21,3);
\draw[very thick] (-1,2.6) .. controls (-1.11,2.75) and (-1.125,2.75) .. (-1.21,3);

\draw[very thick] (0.06,1.98) .. controls (-0.42,2.13) and (-0.42,2.15) .. (-0.75,2.35);

\draw[very thick,draw=white] (1.17,2) .. controls (1.3,2.5) and (1.3,2.5) .. (1.39,3);
\draw[very thick] (1.17,2) .. controls (1.3,2.5) and (1.3,2.5) .. (1.39,3);

\draw[fill] (1.17,2) circle [radius=0.1];
\draw[fill] (0.06,2) circle [radius=0.1];
\draw[fill] (-1.05,2) circle [radius=0.1]; 

\foreach \k in {-1.5,0,1.5}
{\draw[fill] (\k,-1) circle [radius=0.05];};

\foreach \k in {-1.5,0,1.5}
{\draw[fill] (\k,5) circle [radius=0.05];};

\node at (5,4) {$\beta$};

\end{tikzpicture}
\hspace*{\fill}
\caption{A $3$-braid in $\R^2$ illustrating condition~(\ref{it:geombraid3}) of the definition of geometric braid.}\label{fig:braid2}
\end{figure}

\begin{defn}\label{def:geombraid}
A \emph{geometric $n$-braid} in $M$ is a collection $\beta= \brak{\beta_1,\ldots,\beta_n}$ consisting of $n$ arcs $\map{\beta_i}{[0,1]}[M\times [0,1]]$, $i=1,\ldots,n$, called \emph{strings} (or \emph{strands}) such that:
\begin{enumerate}[(a)]
\item for $i=1,\ldots,n$, $\beta_i(0)=(x_{i},0)$ and $\beta_i(1)\in X\times \brak{1}$ (the strings join the elements of $X$ belonging to the copies of $M$ corresponding to $t\in \brak{0,1}$).
\item for all $t\in [0,1]$ and for all $i,j\in \brak{1,\ldots,n}$, $i\neq j$, $\beta_i(t)\neq \beta_j(t)$ (the strings are pairwise disjoint).
\item\label{it:geombraid3} for all $t\in [0,1]$, each string meets the subset $M\times \brak{t}$ in exactly one point (the strings are strictly monotone with respect to the $t$-coordinate).
\end{enumerate}
See Figure~\ref{fig:braid1} for an example of a geometric $3$-braid in the $2$-torus, and Figure~\ref{fig:braid2} for an example of a geometric $3$-braid that illustrates condition~(\ref{it:geombraid3}).
\end{defn}

In the case where $M$ is the plane, a braid is often depicted by a projection (taken to be in general position) onto the plane $xz$ such as that depicted in Figure~\ref{fig:braid2}, so that there are only a finite number of points where the strings cross, and such that the crossings occur at distinct values of $t$. We distinguish between under- and over-crossings.
Our convention is that such a braid is to be read from top to bottom, the top of the braid corresponding to $t=0$, and the bottom to $t=1$. Similar pictures may be drawn for other surfaces of small genus (see~\cite{M,MK} for example).

Two geometric $n$-braids of $M$ are said to be \emph{equivalent} if there exists an isotopy (keeping the endpoints of the strings fixed) from one to the other \emph{through $n$-braids}. In particular, under the isotopy, the strings remain pairwise disjoint. This defines an equivalence relation, and the equivalence classes are termed \emph{$n$-braids}. The set of $n$-braids of $M$ is denoted by $B_n(M)$. By a slight abuse of terminology, we shall not distinguish between a braid and its geometric representatives.

The \emph{product} of two $n$-braids $\alpha$ and $\beta$, denoted $\alpha\beta$, is their concatenation, defined by glueing the endpoints of $\alpha$ to the respective initial points of $\beta$ (formally, $\alpha$ should be `squashed' into the slab $0\leq t\leq \frac{1}{2}$, and $\beta$ into the slab $\frac{1}{2}\leq t\leq 1$). One may check that this operation does not depend on the choice of geometric representatives of $\alpha$ and $\beta$, and that it is \emph{associative}. The \emph{identity element} $\id$ of $B_n(M)$ is the braid all of whose strings are vertical. The \emph{inverse} of an $n$-braid $\beta= \brak{(\beta_1(t),\ldots,\beta_n(t))}_{t\in [0,1]}$ is given by $\beta^{-1}= \brak{(\beta_1(1-t),\ldots,\beta_n(1-t))}_{t\in [0,1]}$ (its mirror image with respect to $M \times \brak{\frac{1}{2}}$). Equipped with this operation, $B_n(M)$ is thus a group, which we call the \emph{$n$-string braid group of $M$}. 

To each $n$-braid $\beta=(\beta_1,\ldots, \beta_n)$, one may associate a permutation $\tau_n(\beta)\in \sn$ defined by $\beta_{i}(1)=(x_{\tau_n(\beta)(i)},1)$, and the following correspondence:
\begin{equation}\label{eq:bnsn}
\begin{aligned}
\tau_n\colon\thinspace B_n(M) & \to \sn\\
\beta& \mapsto \tau_n(\beta)
\end{aligned}
\end{equation}
is seen to be a surjective group homomorphism. The kernel $P_n(M)$ of $\tau_n$ is known as the \emph{$n$-string pure braid group} of $M$, and so $\beta\in P_n(M)$ if and only if $\beta_i(1)=i$ for all $i=1,\ldots,n$. Clearly $P_n(M)$ is a normal subgroup of $B_n(M)$ of index $n!$, and we have the following short exact sequence:
\begin{equation}\label{eq:permseq}
1\to P_n(M)\to B_n(M) \stackrel{\tau_n}{\to} \sn\to 1.
\end{equation}
It is well known that if $M$ is equal to $\R^2$ or to the $2$-disc $\dt$ then $B_{n}(M)$ and $P_{n}(M)$ are isomorphic to the usual \emph{Artin braid groups} $B_{n}$ and $P_{n}$~\cite[Theorem~1.5]{Ha}.

\begin{rem}\label{rem:pnmbnm}
The exact sequence~\reqref{permseq} is frequently used to reduce the study of certain problems in $B_{n}(M)$ to that in $P_{n}(M)$ (see for example Theorems~\ref{th:pnconfig},~\ref{th:centre},~\ref{th:Braidfic},~\ref{th:fics2} and~\ref{th:ficrp}, as well as \repr{nocentre}). The group $B_n(M)$ is also sometimes known as the \emph{permuted} or \emph{full} braid group of $M$, and $P_n(M)$ as the \emph{unpermuted} or \emph{coloured} braid group.
\end{rem}

\subsection{Surface braids as trajectories of non-colliding particles}\label{sec:particles}

\begin{defn}\label{def:noncollide}
Consider $n$ particles which move on the surface $M$, whose initial points are $\gamma_{i}(0)=x_{i}$ for $i=1,\ldots,n$, and whose trajectories are $\gamma_{i}(t)$ for $t\in [0,1]$. A \emph{braid} is thus the collection $\gamma=(\gamma_{1}(t),\ldots,\gamma_{n}(t))_{t\in [0,1]}$ of trajectories satisfying the following two conditions:
\begin{enumerate}[(a)]
\item\label{it:traja} the particles do not collide, \emph{i.e.}\ for all $t\in [0,1]$ and for all $i,j\in \brak{1,\ldots,n}$, $i\neq j$, $\gamma_{i}(t)\neq \gamma_{j}(t)$.
\item\label{it:trajb} they return to their initial points, but possibly undergoing a permutation: $\gamma_{i}(1)\in X$ for all $i\in \brak{1,\ldots,n}$.
\end{enumerate}
\end{defn}

There is a clear bijective correspondence between this definition of braid and the definition of geometric $n$-braid in \resec{defstring}. Indeed, if $\gamma=(\gamma_{1}(t),\ldots,\gamma_{n}(t))_{t\in [0,1]}$ is such a braid then $\beta=(\beta_{1},\ldots,\beta_{n})$ is a geometric $n$-braid, where for all $i=1,\ldots,n$ and $t\in [0,1]$, $\beta_{i}(t)= (\gamma_{i}(t),t)$. Conversely, we may obtain the `particle' notion of braid by reparametrising each string $\beta= (\beta_1,\ldots,\beta_n)$ of a geometric $n$-braid so that $\beta_i(t)$ is of the form $(\gamma_i(t),t)$ for $i=1,\ldots,n$ and $t\in [0,1]$, where $\gamma=(\gamma_{1}(t),\ldots,\gamma_{n}(t))_{t\in [0,1]}$ satisfies conditions~(\ref{it:traja}) and~(\ref{it:trajb}). The transition from a geometric $n$-braid to the `particle notion' may thus be realised geometrically by projecting the strings lying in $M\times [0,1]$ onto the surface $M$.

It is easy to check that two geometric braids are homotopic (in the sense of \resec{defstring}) if and only if the braids defined in terms of trajectories are homotopic. It thus follows that the set of homotopy classes of the latter class of braids may be equipped with a group structure, and that the group thus obtained is isomorphic to $B_{n}(M)$. In this setting, the identity braid is represented by the configuration where all particles are stationary, and the inverse of a braid is given by running through the trajectories in reverse. This point of view proves to be useful when working with braid groups of surface of higher genus, notably in determining presentations~\cite{Be,Bi1,GG1,GM1,S}.

\subsection{Surface braid groups as the fundamental group of configuration spaces}\label{sec:defconfig}

Configuration spaces are important and interesting in their own right~\cite{Coh,FH1}, and have many applications, for example to the study of polynomials in $\C[X]$~\cite{Ha}. The following definition is due to Fox~\cite{FoN} (according to Magnus~\cite{Mag}, the idea first appeared in the work of Hurwitz), and has very important consequences. The motivation for the definition emanates from condition~(\ref{it:geombraid3}) of the definition of geometric $n$-braid given in \resec{defstring}, and is illustrated by Figure~\ref{fig:braid2}.
\begin{defn}
Let $F_n(M)$ denote the \emph{$n\th$ configuration space} of $M$ defined by:
\begin{equation*}
F_n(M)=\set{(p_1,\ldots,p_n)\in M^n}{\text{$p_i\neq p_j$ for all $i,j\in\brak{1,\ldots, n}$, $i\neq j$}}.
\end{equation*}
We equip $F_n(M)$ with the topology induced by the product topology on $M^n$. A transversality argument shows that $F_n(M)$ is a connected $2n$-dimensional open manifold. There is a natural free action of the symmetric group $\sn$ on $F_n(M)$ by permutation of coordinates. The resulting orbit space $F_n(M)/\sn$ shall be denoted by $D_n(M)$, the \emph{$n\th$ permuted configuration space} of $M$, and may be thought of as the configuration space of $n$ \emph{unordered} points. The associated canonical projection $\map{\widehat{\rho}_{n}}{F_n(M)}[D_n(M)]$ is thus a regular $n!$-fold covering map.
\end{defn}

We may thus describe $F_n(M)$ as $M^n\setminus\Delta$, where $\Delta$ denotes the `fat diagonal' of $M^n$:
\begin{equation*}
\Delta=\set{(p_1,\ldots,p_n)\in M^n}{\text{$p_i=p_j$ for some $1\leq i< j\leq n$}}.
\end{equation*}
If $M=\R^2$ then $\displaystyle \Delta= \bigcup_{1\leq i<j\leq n} H_{i,j}$, where $H_{i,j}$ is the hyperplane defined by:
\begin{equation*}
H_{i,j}=\set{(p_1,\ldots,p_n)\in (\R^2)^n}{p_i=p_j}.
\end{equation*}

The following theorem is fundamental, and brings in to play a topological definition of the braid groups that will be very important in what follows. The proof is a good illustration of the use of the short exact sequence~\reqref{permseq}.
\begin{thm}[Fox and Neuwirth~\cite{FoN}]\label{th:pnconfig}
Let $n\in \N$. Then $P_n(M)\cong \pi_1(F_n(M))$ and $B_n(M) \cong \pi_1(D_n(M))$.
\end{thm}

\begin{rems}\mbox{}
\begin{enumerate}[(a)]
\item Since $F_{1}(M)=M$, we have that $B_{1}(M)\cong P_{1}(M)\cong \pi_{1}(M)$. The braid groups of $M$ may thus be seen as generalisations of its fundamental group.
\item The fact that $F_{n}(M)$ (resp.\ $D_{n}(M)$) is connected implies that the isomorphism class of $\pi_1(F_n(M))$ (resp.\ $\pi_1(D_n(M))$) does not depend on the choice of basepoint. We thus have two finite-dimensional topological spaces $F_n(M)$ (resp.\ $D_n(M)$) whose fundamental groups are $P_n(M)$ (resp.\ $B_n(M)$). As we shall see in \resec{fnsequence}, the relations between configuration spaces and braid groups play a fundamental rôle in the study of the latter, notably via the fact that we may form certain natural fibre spaces of the former.
\item The definitions of surface braid groups given in Sections~\ref{sec:defstring}--\ref{sec:defconfig} generalise to any topological space. It was shown in~\cite[Theorem~9]{FaN} that for connected manifolds of dimension $r\geq 3$, there is no braid theory, as it is formulated here.
\end{enumerate}
\end{rems}

The natural inclusion $\altarrow{\iota}{F_{n}(M)}{\lhra}{M^n}$ induces a homomorphism of the corresponding fundamental groups:
\begin{equation*}
\map{\iota_{\#}}{P_{n}(M)}[(\pi_{1}(M))^n],
\end{equation*}
and the inclusion $\altarrow{j}{\dt}{\lhra}{\Int{M}}$ of a topological disc $\dt$ in the interior of $M$ induces a homomorphism $\map{j_{\#}}{P_{n}}[P_{n}(M)]$ that is an embedding for most surfaces:
\begin{prop}[\cite{Bi1}]\label{prop:birembed}
Let $M$ be a compact, orientable surface different from $\St$. Then the inclusion $\altarrow{j}{\dt}{\lhra}{M}$ induces an embedding $P_{n}\lhra P_{n}(M)$.
\end{prop}
\repr{birembed} extends first to the non-orientable case~\cite{Go}, with the exception of $M=\rp$, and secondly, to the full braid groups by applying \req{permseq}. If $M$ is different from $\St$ and $\rp$ then Goldberg showed that the following sequence is short exact~\cite{Go}:
\begin{equation}\label{eq:goldberg}
1 \to \normcl{\im{j_{\#}}}_{P_{n}(M)} \lhra P_{n}(M) \stackrel{\iota_{\#}}{\to} (\pi_{1}(M))^n\to 1,
\end{equation}
where $\normcl{H}_{G}$ denotes the normal closure of a subgroup $H$ in a group $G$. This sequence was analysed in~\cite{GMP} in order to study Vassiliev invariants of braid groups of orientable surfaces. In the case of $\rp$, $\ker{\iota_{\#}}$ was computed and the homotopy fibre of $\iota$ was determined in~\cite{GG15}. 

\subsection{Relationship between braid and mapping class groups}\label{sec:mcg}

Let $M$ be a compact, connected, orientable (resp.\ non-orientable) surface, possibly with boundary $\partial M$, and for $n\geq 0$, let $\mathcal{Q}_{n}$ be a finite subset of $\Int{M}$ consisting of $n$ distinct points (so $\mathcal{Q}_{0}=\vide$). Let $\mathcal{H}(M,\mathcal{Q}_{n})$ denote the group $\operatorname{\text{Homeo}^+}(M,\mathcal{Q}_{n})$ (resp.\ $\operatorname{\text{Homeo}}(M,\mathcal{Q}_{n})$) of orientation-preserving homeomorphisms (resp.\ of homeomorphisms) of $M$ under composition that leave $\mathcal{Q}_{n}$ invariant (so we allow the points of $\mathcal{Q}_{n}$ to be permuted), and that fix $\partial M$ pointwise. We equip $\mathcal{H}(M,\mathcal{Q}_{n})$ with the compact-open topology. Let $\mathcal{H}_{0}(M,\mathcal{Q}_{n})$ denote the path component of $\id_{M}$ in $\mathcal{H}(M,\mathcal{Q}_{n})$. The \emph{$n\th$ mapping class group} of $M$, denoted by $\mcggen{M}$, is defined to be the set of isotopy classes of the elements of $\operatorname{\text{Homeo}^+}(M,\mathcal{Q}_{n})$ (resp.\ $\operatorname{\text{Homeo}}(M,\mathcal{Q}_{n})$), in other words,
\begin{equation*}
\mcggen{M}= \mathcal{H}(M,\mathcal{Q}_{n})/\mathcal{H}_{0}(M,\mathcal{Q}_{n})=\pi_{0}(\mathcal{H}(M,\mathcal{Q}_{n})).
\end{equation*}
It is straightforward to check that $\mcggen{M}$ is indeed a group whose isomorphism class does not depend on the choice of $\mathcal{Q}_{n}$. If $n=0$ then we write simply $\mathcal{H}(M)$ and $\mcgzero{M}$ for the corresponding groups. The mapping class groups have been widely studied and play an important rôle in low-dimensional topology. Some good general references are~\cite{Bi2,FM,I}. 

The mapping class groups are closely related to braid groups. If $M=\dt$ then it is well known that they coincide:
\begin{thm}[\cite{Bi2,Ha,KT}]\label{th:mcgdisc}
$B_n\cong \mathcal{MCG}(\dt,n)$.
\end{thm}
The proof of \reth{mcgdisc} makes use of Artin's representation of $B_{n}$ as a subgroup of the automorphism group of the free group $\F[n]$ of rank $n$, the free group in question being identified with $\pi_{1}(\dt \setminus \mathcal{Q}_{n})$. In the general case, the relationship between $\mcggen{M}$ and $B_{n}(M)$ arises in a topological setting as follows~\cite{Bi1a,Bi2,S}. Let $n\geq 1$, and fix a basepoint $\mathcal{Q}_{n}\in D_{n}(M)$. Then the map $\map{\Psi}{\mathcal{H}(M)}[D_{n}(M)]$ defined by $\Psi(f)=f(\mathcal{Q}_{n})$ is a locally-trivial fibre bundle~\cite{Bi1a,McC}, whose fibre over $\mathcal{Q}_{n}$ is equal to $\mathcal{H}(M,\mathcal{Q}_{n})$. Taking the long exact sequence in homotopy of this fibration yields:
\begin{multline}\label{eq:leshomotopy}
\cdots \to \pi_{1}(\mathcal{H}(M,\mathcal{Q}_{n})) \to \pi_{1}(\mathcal{H}(M)) \to
\pi_{1}(D_{n}(M))\to\\ \pi_{0}(\mathcal{H}(M,\mathcal{Q}_{n})) \to \pi_{0}(\mathcal{H}(M)) \to 1.
\end{multline}
If $M$ is different from $\St$, $\rp$, the torus or the Klein bottle then $\pi_{1}(\mathcal{H}(M))=1$~\cite{Ham}, from which we deduce a short exact sequence of the form:
\begin{equation}\label{eq:mcggen}
1\to B_{n}(M)\to  \mcggen{M}\to \mcgzero{M} \to 1.
\end{equation}
The braid group $B_{n}(M)$ is thus isomorphic to the kernel of the homomorphism that corresponds geometrically to forgetting the marked points. We recover \reth{mcgdisc} by noting that $\mcgzero{\dt}=\brak{1}$ using the Alexander trick. If $M=\St$ (resp.\ $\rp$) and $n\geq 3$ (resp.\ $n\geq 2$) then $\pi_{1}(\mathcal{H}(M,\mathcal{Q}_{n}))=1$~\cite{Ham,McC}, but $\pi_{1}(\mathcal{H}(M))\cong \Z_{2}$~\cite{Ham0,Ham}, which is a manifestation of the fact that the fundamental group of $\operatorname{SO}(3)$ is isomorphic to $\Z_{2}$~\cite{Fa,Ha,Ne}. In this case, we obtain the following short exact sequence:
\begin{equation}\label{eq:mcg}
1 \to \Z_{2}\to B_{n}(M)\to \mcggen{M} \to 1.
\end{equation}
As we shall see in \resec{sphere}, viewed as an element of $B_{n}(M)$, the generator of the kernel is the full twist braid $\ft$~\cite{FvB,vB}. In particular, $B_{n}(M)/\ang{\ft}\cong \mcggen{M}$. In the case of $\St$, the short exact sequence~\reqref{mcg} may be obtained by combining the presentation of $\mcggen{\St}$ due to Magnus~\cite{Mag0,MKS} with  Fadell and Van Buskirk's presentation of $B_{n}(\St)$ (see~\reth{presfvb}). It plays an important part, notably in the study of the centralisers and conjugacy classes of the finite order elements, and of the finite subgroups of $B_{n}(M)$ (see \resec{finitesubgp}).  Finally, if $M$ is the torus $\mathbb{T}^{2}$ or the Klein bottle then~\reqref{leshomotopy} yields a six-term exact sequence starting and ending with $1$. In the case of $\mathbb{T}^{2}$, this exact sequence involves $\mcgzero{\mathbb{T}^{2}}$, which is isomorphic to $\operatorname{SL}(2,\Z)$.

\section{Some properties of surface braid groups}\label{sec:properties}

In this section, we describe various properties of surface braid groups. We start with one of the most important, that makes use of the definition of \resec{defconfig} in terms of configuration spaces.

\subsection{Exact sequences of braid groups}\label{sec:fnsequence}

Let $M$ be a connected surface. For $n\in \N$, we equip $F_{n}(M)$ with the topology induced by the product topology on the $n$-fold Cartesian product $M^n$. For $m\geq 0$, let $\mathcal{Q}_{m}$ be as in \resec{mcg}, and set $F_{m,n}(M)=F_{n}(M\setminus \mathcal{Q}_{m})$ and $D_{m,n}(M)$ to be the quotient space of $F_{m,n}(M)$ by the free action of $S_{n}$, so that the projection $F_{m,n}(M)\to D_{m,n}(M)$ is a covering map. Note that the topological type of $F_{m,n}(M)$ does not depend on the choice of $\mathcal{Q}_{m}$, and that as special cases, we obtain $F_{0,n}(M)=F_{n}(M)$ and $F_{m,1}(M)=M\setminus \mathcal{Q}_{m}$. We have the following important result concerning the topological structure of the spaces $F_{m,n}(M)$.

\begin{thm}[Fadell and Neuwirth~\cite{FaN,Ha,KT}]\label{th:fnfib}
Let $1\leq r<n$ and $m\geq 0$. Suppose that $M$ is a surface with empty boundary. Then the map
\begin{equation}\label{eq:fnles}
\left\{\begin{aligned}
p_{n,r} \colon\thinspace F_{m,n}(M) &\to F_{m,r}(M)\\
(x_{1},\ldots,x_{n}) &\mapsto (x_{1},\ldots,x_{r}) 
\end{aligned}\right.
\end{equation}
is a locally-trivial fibration with fibre $F_{m+r,n-r}(M)$.
\end{thm}

One may then take the long exact sequence in homotopy of the fibration~\reqref{fnles}:
\begin{multline}\label{eq:fnleslong}
\cdots \to \pi_{k}(F_{m+r,n-r}(M)) \to \pi_{k}(F_{m,n}(M)) \to \pi_{k}(F_{m,r}(M))\to \\
\pi_{k-1}(F_{m+r,n-r}(M)) \to \cdots\to \pi_{2}(F_{m+r,n-r}(M)) \to \pi_{2}(F_{m,n}(M)) \to \pi_{2}(F_{m,r}(M))\to \\
\pi_{1}(F_{m+r,n-r}(M)) \to \pi_{1}(F_{m,n}(M)) \to \pi_{1}(F_{m,r}(M))\to 1.
\end{multline}
Since $F_{m+n+i-1,1}(M)$ has the homotopy type of a bouquet of circles for all $0\leq i\leq n-2$, it follows that:
\begin{equation*}
\pi_{k}(F_{m,n}(M))\cong \pi_{k}(F_{m,n-1}(M))\cong \cdots \cong \pi_{k}(F_{m,1}(M))=\pi_{k}(M\setminus \mathcal{Q}_{m})\quad \text{for all $k\geq 3$,}
\end{equation*}
and that the homomorphism $\pi_{2}(F_{m,n-i}(M)) \to \pi_{2}(F_{m,n-i-1}(M))$ is injective for all such $i$. Thus $\pi_{2}(F_{m,n}(M))$ is isomorphic to a subgroup of $\pi_{2}(F_{m,1}(M))$, which is in turn isomorphic to $\pi_{2}(M\setminus \mathcal{Q}_{m})$. Since $\pi_{1}(F_{m,n}(M))\cong P_{n}(M\setminus \mathcal{Q}_{m})$ by \reth{pnconfig}, we recover the following result:

\begin{thm}[\cite{Fa,FaN,FvB,vB}]\label{th:fnses}\mbox{}
\begin{enumerate}[(a)]
\item\label{it:eilenmac} Let $n\in \N$ and $m\geq 0$. We suppose additionally that $M$ is different from $\St$ and $\rp$ if $m=0$. Then the spaces $F_{m,n}(M)$ and $D_{m,n}(M)$ are Eilenberg-Mac~Lane spaces of type $K(P_n(M\setminus \mathcal{Q}_{m}),1)$ and $K(B_n(M\setminus \mathcal{Q}_{m}),1)$ respectively.

\item\label{it:fvb} If $n\geq 3$ (resp.\ $n\geq 2$) then $\pi_{2}(F_{n}(\St))=0$ and $\pi_{2}(F_{n}(\rp))=0$.

\item Let $1\leq r<n$ and $m\geq 0$. If $m=0$ then we suppose that $r\geq 3$ if $M=\St$, and that $r\geq 2$ if $M=\rp$. Then the Fadell-Neuwirth fibration~\reqref{fnleslong} induces a short exact sequence:
\begin{equation}\label{eq:fnses}
1\to P_{n-r}(M\setminus \mathcal{Q}_{m+r}) \to P_{n}(M\setminus \mathcal{Q}_{m})\xrightarrow{(p_{n,r})_{\#1}} P_{r}(M\setminus \mathcal{Q}_{m})\to 1.
\end{equation}
\end{enumerate}
\end{thm}

\begin{rems}\mbox{}\label{rem:fnses}
\begin{enumerate}[(a)]
\item The short exact sequence~\reqref{fnses} is known as the \emph{Fadell-Neuwirth short exact sequence of surface braid groups.} It plays a central rôle in the study of surface (pure) braid groups. It was used to study mapping class groups in~\cite{PR}, and in work on Vassiliev invariants for braid groups~\cite{GMP}.
\item \reth{fnses}(\ref{it:fvb}) was proved in~\cite{Fa,FvB,vB} by showing that $\pi_{2}(F_{3}(\St))=\pi_{2}(F_{2}(\rp))=0$ and using induction.
\item The projection $P_{n}(M\setminus \mathcal{Q}_{m})\to P_{r}(M\setminus \mathcal{Q}_{m})$ may be interpreted geometrically as the epimorphism that `forgets' the last $n-r$ strings.
\item\label{it:fnsesb} In order to prove that~\reqref{fnles} is a locally-trivial fibration, one needs to suppose that $M$ is without boundary. However, the long exact sequence~\reqref{fnleslong} exists even if $M$ has boundary, and thus \reth{fnses} holds for any connected surface. To see this, let $M$ be a surface with boundary, and let $M'=M\setminus \partial M$. Then $M'$ is a surface with empty boundary, and so Theorems~\ref{th:fnfib} and~\ref{th:fnses} hold for $M'$. The inclusion of $M'$ in $M$ is not only a homotopy equivalence between $M'$ and $M$, but it also induces a homotopy equivalence between their $n\up{th}$ configuration spaces. In particular,~\reqref{fnleslong} and~\reth{fnses} are valid also for $M$, and the $n\up{th}$ (pure) braid groups of $M'$ and $M$ are isomorphic.

\item Let $n\geq 4$ if $M=\St$, $n\geq 3$ if $M=\rp$, and $n\geq 2$ otherwise. Two special cases to which we will refer frequently are:
\begin{enumerate}[(i)]
\item $m=0$, in which case the short exact sequence~\reqref{fnses} becomes:
\begin{equation}\label{eq:sesgennr}
1 \to P_{n-r}(M\setminus \mathcal{Q}_{r}) \to P_{n}(M) \xrightarrow{(p_{n,r})_{\#1}} P_{r}(M) \to 1.
\end{equation}
\item $m=0$ and  $r=n-1$, in which case the short exact sequence~\reqref{fnses} becomes:
\begin{equation}\label{eq:sesgen}
1 \to \pi_{1}(M\setminus \mathcal{Q}_{n-1}) \to P_{n}(M) \xrightarrow{(p_{n,n-1})_{\#1}} P_{n-1}(M) \to 1.
\end{equation}
In particular, each element of $\ker{(p_{n,n-1})_{\#1}}$ may be interpreted as an $n$-string braid whose first $n-1$ strings are vertical. This short exact sequence lends itself naturally to induction, and may be used for example to solve the word problem in surface braid groups~\cite{A2,GVB,S}, and to obtain presentations (see \resec{present}).
\end{enumerate}
\end{enumerate}
\end{rems}

By a theorem of P.~A.~Smith (see~\cite[page~149]{Hat} or~\cite[page~287]{Hu}), the fundamental groups of finite-dimensional Eilenberg-Mac~Lane spaces of type $K(\pi,1)$ are torsion free. This implies immediately the sufficiency of the following assertions:
\begin{cor}[\cite{FaN,FvB,vB}]\label{cor:torsion}
Let $M$ be a surface. Then the braid groups $P_n(M\setminus \mathcal{Q}_{m})$ and $B_n(M\setminus \mathcal{Q}_{m})$ are torsion free if and only if either:
\begin{enumerate}[(a)]
\item $m\geq 1$, or
\item $m=0$ and $M$ is a surface different from $\St$ and $\rp$.
\end{enumerate}
\end{cor}
As for the necessity of the conditions, we already mentioned in \resec{mcg} that the full twist $\ft$ is an element of $P_{n}(M)$ of order~$2$ if $M=\St$ or $\rp$. The existence of torsion in the braid groups of $\St$ and $\rp$ is a fascinating phenomenon to which we shall return in Sections~\ref{sec:sphere} and~\ref{sec:finitesubgp}, and makes for interesting and intricate $K$-theoretical structure (see \resec{ktheoryresults}). More will be said about the homotopy groups of the configuration spaces of the exceptional surfaces, $\St$ and $\rp$, in \resec{homotype}. 

We remark that a purely algebraic proof of the fact that the Artin braid groups are torsion free was given later by Dyer~\cite{Dy}. We shall see another proof in \resec{orderable}.

The short exact sequences~\reqref{fnses}--\reqref{sesgen} do not extend directly to the full braid groups, but may be generalised  as follows to certain subgroups that lie between $P_{n}(M)$ and $B_{n}(M)$. Once more, let $1\leq r<n$, and suppose that $r\geq 3$ if $M=\St$ and $r\geq 2$ if $M=\rp$. We consider the space obtained by taking the quotient of $F_{n}(M)$ by the subgroup $\sn[r]\times \sn[n-r]$ of $\sn[n]$. If $M$ is without boundary then as in \reth{fnfib} we obtain a locally-trivial fibration $\map{q_{n,r}}{F_{n}(M)/(\sn[r]\times \sn[n-r])}[D_r(M)]$, defined by forgetting the last $n-r$ coordinates. We set $B_{r,n-r}(M)=\pi_1\bigl(F_{n}(M)/(\sn[r]\times \sn[n-r])\bigr)$, which is often termed a \emph{`mixed'} braid group, and is defined whether or not $M$ has boundary. As in the pure braid group case, we obtain the following generalisation of~\reqref{sesgennr}~\cite{GG2}:
\begin{equation}\label{eq:fnsurface}
1\to B_{n-r}(M\setminus\mathcal{Q}_{r}) \to B_{r,n-r}(M)\xrightarrow{(q_{n,r})_{\#1}} B_r(M)\to 1,
\end{equation}
known as a \emph{generalised Fadell-Neuwirth short exact sequence} of mixed braid groups. Such braid groups are very useful, and have been studied in~\cite{BGG2,GG2,GG4,GG11,Lam,Man,PR} for example. Further generalisations are possible by taking quotients by direct products of the form $\sn[i_{1}]\times \cdots \times \sn[i_{r}]$, where $\displaystyle \sum_{j=1}^{r}\, i_{j}=n$.

\subsection{Presentations of surface braid groups}\label{sec:present}

We recall the classical presentation of the Artin braid groups:
\begin{thm}[Artin, 1925~\cite{A1}]\label{th:presbn}
For all $n\geq 1$, the braid group $B_n$ admits the following presentation
\begin{enumerate}
\item[\underline{\textbf{generators}:}] $\sigma_1,\ldots, \sigma_{n-1}$ (known as the \emph{Artin generators}).
\item[\underline{\textbf{relations}:}] (known as the \emph{Artin relations}) 
\begin{gather}
\text{$\sigma_{i}\sigma_{j}=\sigma_{j}\sigma_{i}$ if $\lvert i-j\rvert\geq 2$
and $1\leq i,j\leq n-1$}\label{eq:Artin1}\\
\text{$\sigma_{i}\sigma_{i+1}\sigma_{i}=\sigma_{i+1}\sigma_{i}\sigma_{i+1}$ for
all $1\leq i\leq n-2$.}\label{eq:Artin2}
\end{gather}
\end{enumerate} 
\end{thm}

\begin{figure}[h]
\hfill
\begin{tikzpicture}[scale=0.75, very thick]

\foreach \k in {5}
{\draw (\k,3) .. controls (\k,2) and (\k-1,2) .. (\k-1,1);};

\foreach \k in {4}
{\draw[white,line width=6pt] (\k,3) .. controls (\k,2) and (\k+1,2) .. (\k+1,1);
\draw (\k,3) .. controls (\k,2) and (\k+1,2) .. (\k+1,1);};

\foreach \k in {15}
{\draw (\k,3) .. controls (\k,2) and (\k+1,2) .. (\k+1,1);};

\foreach \k in {16}
{\draw[white,line width=6pt] (\k,3) .. controls (\k,2) and (\k-1,2) .. (\k-1,1);
\draw (\k,3) .. controls (\k,2) and (\k-1,2) .. (\k-1,1);};

\foreach \k in {1,3,6,8,12,14,17,19}
{\draw (\k,1)--(\k,3);};

\foreach \k in {2,7,13,18}
{\node at (\k,2) {$\cdots$};};

\foreach \k in {1,12}
{\node at (\k,3.5) {$1$};
\node at (\k+1.9,3.5) {$i-1$};
\node at (\k+2.9,3.52) {$i$};
\node at (\k+3.85,3.5) {$i+1$};
\node at (\k+5.2,3.5) {$i+2$};
\node at (\k+7,3.5) {$n$};};

\node at (4.5,0.25) {$\sigma_{i}$};
\node at (15.5,0.25) {$\sigma_{i}^{-1}$};

\end{tikzpicture}
\hspace*{\fill}
\caption{The braid $\sigma_{i}$ and its inverse.}\label{fig:braid3}
\end{figure}
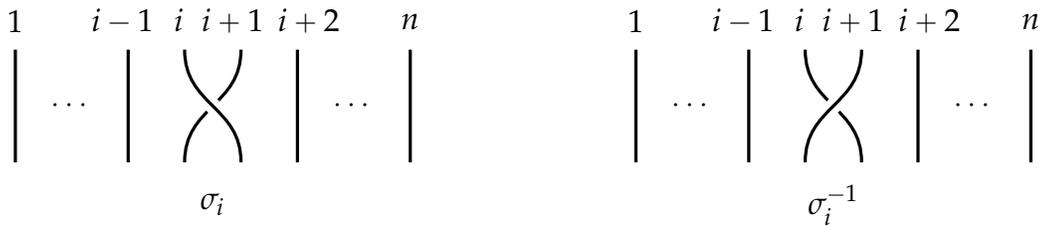

The generator $\sigma_i$ may be regarded geometrically as the braid with a single positive crossing of the $i\th$ string with the $(i+1)\textsuperscript{st}$ string, while all other strings remain vertical (see Figure~\ref{fig:braid3}). 
It follows from this presentation that $B_1=\brak{1}$ and $B_2=\ang{\sigma_1}\cong \Z$. Adding the relations $\sigma_i^2=1$, $i=1,\ldots,n-1$, to those  of \reth{presbn} yields the Coxeter presentation of $\sn$. If $1\leq i<j\leq n$, the pure braid defined by:
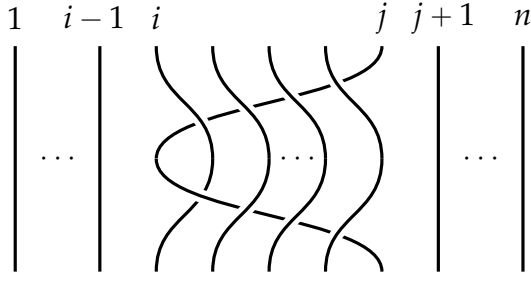
\begin{figure}[h]
\hfill
\begin{tikzpicture}[scale=0.75, very thick]
\draw (7,6.5) .. controls (7,5.5) and (3,5.5) .. (3,4.5);
\foreach \j in {3,...,6} 
{\draw[white,line width=5pt] (\j,6.5).. controls (\j,5.5) and (\j+1,5.5) .. (\j+1,4.5);
\draw (\j,6.5).. controls (\j,5.5) and (\j+1,5.5) .. (\j+1,4.5);
\draw[white,line width=5pt] (3,4.5) .. controls (3,3.5) and (7,3.5) .. (7,2.5);};
\draw (4,4.5) .. controls (4,3.5) and (3,3.5) .. (3,2.5);
\draw[white,line width=5pt] (3,4.5) .. controls (3,3.5) and (7,3.5) .. (7,2.5);
\draw (3,4.5) .. controls (3,3.5) and (7,3.5) .. (7,2.5);
\foreach \j in {4,...,6} 
{\draw[white,line width=5pt] (\j+1,4.5).. controls (\j+1,3.5) and (\j,3.5) .. (\j,2.5);
\draw (\j+1,4.5).. controls (\j+1,3.5) and (\j,3.5) .. (\j,2.5);
};
\foreach \j in {0.5,2,8,9.5}
{\draw (\j,6.5)-- (\j,2.5);};
\foreach \k in {0.75,5,8.25}
{\node at (\k+0.55,4.5) {$\cdots$};};
\node at (0.5,7) {$1$};
\node at (1.9,7) {$i-1$};
\node at (3,7) {$i$};
\node at (7,7) {$j$};
\node at (8.1,7) {$j+1$};
\node at (9.5,7) {$n$};
\end{tikzpicture}
\hspace*{\fill}
\caption{The element $A_{i,j}$ of $B_{n}$.}\label{fig:aij}
\end{figure}
\begin{equation}\label{eq:defaij}
A_{i,j}=\sigma_{j-1}\cdots \sigma_{i+1}\sigma_{i}^2\sigma_{i+1}^{-1}\cdots \sigma_{j-1}^{-1},
\end{equation}
may be represented geometrically by the braid all of whose strings are vertical, with the exception of the $j\textsuperscript{th}$ string, that wraps around the $i\textsuperscript{th}$ string (see Figure~\ref{fig:aij}). 
Such elements generate $P_{n}$:
\begin{prop}[\cite{Ha}]\label{prop:hansen}
For all $n\geq 1$, $P_n$ is generated by  $\setl{A_{i,j}}{1\leq i<j\leq n}$ whose elements are subject to the following relations:
\begin{equation*}
A_{r,s}^{-1}A_{i,j}A_{r,s} =
\begin{cases}
A_{i,j} & \text{if $ i < r < s < j$ or $r < s < i < j$}\\
A_{r,j}A_{i,j}A_{r,j}^{-1} & \text{if $r < i = s < j$}\\
A_{r,j}A_{s,j}A_{i,j}A_{s,j}^{-1}A_{r,j}^{-1} & \text{if $i = r < s < j$}\\
A_{r,j}A_{s,j}A_{r,j}^{-1}A_{s,j}^{-1}A_{i,j}A_{s,j}A_{r,j}A_{s,j}^{-1}A_{r,j}^{-1} & \text{if $ r < i < s < j$.}
\end{cases}
\end{equation*}
\end{prop}

One interesting fact that may be deduced immediately from the presentation of \repr{hansen} is that the action by conjugation of $P_{n}$ on itself induces the identity on the Abelianisation of $P_{n}$, and via the short exact sequence~\reqref{sesgen} in the case where $M=\R^2$, implies that $P_{n}$ is an \emph{almost-direct} product of $\F[n-1]$ and $P_{n-1}$. This plays an important rôle in various aspects of the theory, for example in the proof of the fact that $P_{n}$ is residually nilpotent~\cite{FaRa,FaRa2}.

A number of presentations are known for surface (pure) braid groups~\cite{Be,Bi1,GG1,GG7,GM1,Lad,Lam,S,Z1,Z2}, the first being due to Birman and Scott. We recall those due to Bellingeri for $B_{n}(N)$, where $N$ is a connected surface of the form $M\setminus \mathcal{Q}_{m}$, $M$ being compact and without boundary, and orientable in the first case, and non-orientable in the second. One way to find such presentations is to apply standard techniques to obtain presentations of group extensions~\cite{J}. One first uses induction and the short exact sequence~\reqref{sesgen} to obtain presentations of the pure braid groups, and then~\reqref{permseq} yields presentations of the full braid groups.

\enlargethispage{2mm}
\begin{thm}[\cite{Be}]\label{th:belli1}
Let $M$ be a compact, connected, orientable surface without boundary of genus $g$, where $g\geq 1$, and let $m\geq 0$. Then $B_{n}(M\setminus \mathcal{Q}_{m})$ admits the following presentation:
\begin{enumerate}
\item[Generators:] $\sigma_{1},\ldots, \sigma_{n-1}, a_{1},\ldots, a_{g}, b_{1},\ldots, b_{g}, z_{1},\ldots,z_{m-1}$.
\item[Relations:]\mbox{}
\begin{enumerate}[(a)] 
\item the Artin relations~\reqref{Artin1} and~\reqref{Artin2}.
\item $a_{r}\sigma_{i}=\sigma_{i}a_{r}$, $b_{r}\sigma_{i}=\sigma_{i}b_{r}$ and $z_{j}\sigma_{i}=\sigma_{i}z_{j}$ for all $1\leq r\leq g$, $2\leq i\leq n-1$ and $1\leq j\leq m-1$.
\item $(\sigma_{1}^{-1} a_{r})^{2} =(a_{r}\sigma_{1}^{-1})^{2}$, $(\sigma_{1}^{-1} b_{r})^{2} =(b_{r}\sigma_{1}^{-1})^{2}$ and $(\sigma_{1}^{-1} z_{j})^{2}= (z_{j}\sigma_{1}^{-1})^{2}$ for all $1\leq r\leq g$ and $1\leq j\leq m-1$.
\item $\sigma_{1}^{-1} a_{s}\sigma_{1} a_{r}=a_{r}\sigma_{1}^{-1} a_{s}\sigma_{1}$, $\sigma_{1}^{-1} b_{s}\sigma_{1} b_{r}=b_{r}\sigma_{1}^{-1} b_{s}\sigma_{1}$, $\sigma_{1}^{-1} a_{s}\sigma_{1} b_{r}=b_{r}\sigma_{1}^{-1} a_{s}\sigma_{1}$ and $\sigma_{1}^{-1} b_{s}\sigma_{1} a_{r}=a_{r}\sigma_{1}^{-1} b_{s}\sigma_{1}$ for all $1\leq s<r\leq g$.

\item if $n\geq 2$, $\sigma_{1}^{-1} z_{j}\sigma_{1} a_{r}=a_{r}\sigma_{1}^{-1} z_{j}\sigma_{1}$ and $\sigma_{1}^{-1} z_{j}\sigma_{1} b_{r}=b_{r}\sigma_{1}^{-1} z_{j}\sigma_{1}$ for all $1\leq r\leq g$ and $1\leq j\leq m-1$. 

\item $\sigma_{1}^{-1} z_{j}\sigma_{1} z_{l}=z_{l}\sigma_{1}^{-1} z_{j}\sigma_{1}$ for all $1\leq j<l\leq m-1$.

\item $\sigma_{1}^{-1} a_{r}\sigma_{1}^{-1} b_{r}=b_{r}\sigma_{1}^{-1} a_{r}\sigma_{1}$ for all $1\leq r\leq g$.

\item if $m=0$ then $[a_{1},b_{1}^{-1}]\cdots [a_{g},b_{g}^{-1}]=\sigma_{1}\cdots \sigma_{n-2}\sigma_{n-1}^{2} \sigma_{n-2} \cdots \sigma_{1}$, where $[a,b]=aba^{-1}b^{-1}$.
\end{enumerate}
\end{enumerate}
\end{thm}

\begin{thm}[\cite{Be}]\label{th:belli2}
Let $M$ be a compact, connected, non-orientable surface without boundary of genus $g$, where $g\geq 2$, and let $m\geq 0$. Then $B_{n}(M\setminus \mathcal{Q}_{m})$ admits the following presentation:
\begin{enumerate}
\item[Generators:] $\sigma_{1},\ldots, \sigma_{n-1}, a_{1},\ldots, a_{g}, z_{1},\ldots,z_{m-1}$.
\item[Relations:]\mbox{}
\begin{enumerate}[(a)] 
\item the Artin relations~\reqref{Artin1} and~\reqref{Artin2}.
\item $a_{r}\sigma_{i}=\sigma_{i}a_{r}$ for all $1\leq r\leq g$ and $2\leq i\leq n-1$.

\item $(\sigma_{1}^{-1} a_{r})^{2} =a_{r}\sigma_{1}^{-1}a_{r}\sigma_{1}$ and $(\sigma_{1}^{-1} z_{j})^{2}= (z_{j}\sigma_{1}^{-1})^{2}$ for all $1\leq r\leq g$ and $1\leq j\leq m-1$.

\item $\sigma_{1}^{-1} a_{s}\sigma_{1} a_{r}=a_{r}\sigma_{1}^{-1} a_{s}\sigma_{1}$ for all $1\leq s<r\leq g$.

\item $z_{j}\sigma_{i}=\sigma_{i}z_{j}$ for all $2\leq i\leq n-1$ and $1\leq j\leq m-1$.

\item if $n\geq 2$, $\sigma_{1}^{-1} z_{j}\sigma_{1} a_{r}=a_{r}\sigma_{1}^{-1} z_{j}\sigma_{1}$ for all $1\leq r\leq g$ and $1\leq j\leq m-1$. 

\item $\sigma_{1}^{-1} z_{j}\sigma_{1} z_{l}=z_{l}\sigma_{1}^{-1} z_{j}\sigma_{1}$ for all $1\leq j<l\leq m-1$.

\item if $m=0$ then $a_{1}^{2}\cdots a_{g}^{2}=\sigma_{1}\cdots \sigma_{n-2}\sigma_{n-1}^{2} \sigma_{n-2} \cdots \sigma_{1}$.
\end{enumerate}
\end{enumerate}
\end{thm}

\begin{rems}\mbox{}
\begin{enumerate}[(a)]
\item Geometrically, the generators $a_{1},\ldots, a_{g}, b_{1},\ldots, b_{g}$ (resp.\ $a_{1},\ldots, a_{g}$) of $B_{n}(M)$ given in \reth{belli1} (resp.\ \reth{belli2}) correspond to a standard set of generators of $\pi_{1}(M)$ based at the first basepoint of the braid, and in both cases, the generator $z_{i}$, $i\in \brak{1,\ldots,m-1}$, corresponds to the braid all of whose strings are vertical, with the exception of the first string that wraps around the $i\up{th}$ puncture. 
\item By \rerems{fnses}(\ref{it:fnsesb}), it follows that we may also take some or all of the punctures to be boundary components. In other words, Theorems~\ref{th:belli1} and~\ref{th:belli2} yield presentations of the braid groups of any surface as defined at the beginning of \resec{basic}.
\end{enumerate}
\end{rems}

Presentations for $B_{n}(\St)$ and $B_{n}(\rp)$ will be given in \resec{sphere}. Results on the minimal cardinality of different types of generating sets of $B_{n}(M)$, where $M=\dt$, $\St$ or $\rp$, are given in~\cite{GG13}. Positive presentations of braid groups of orientable surfaces were obtained in~\cite{BG}. Braid groups of the annulus, which are Artin-Tits groups of type $B_{n}$, were studied in~\cite{Cr,GGjktr2,KP,Lam,Man,PR}.

\subsection{The centre of surface braid groups}\label{sec:centre}

\begin{figure}[h]
\hfill
\begin{tikzpicture}[scale=0.75, very thick]
\foreach \j in {2,3,...,6} 
{\foreach \k in {2,...,\j}
{\draw (\k,\j) .. controls (\k,\j-0.5) and (\k-1,\j-0.5) .. (\k-1,\j-1);};
{\draw[white,line width=5pt] (1,\j) .. controls (1,\j-0.5) and (\j,\j-0.5) .. (\j,\j-1);
\draw (1,\j) .. controls (1,\j-0.5) and (\j,\j-0.5) .. (\j,\j-1);};
if \j>2 then \draw (\j,1) -- (\j,\j-1);
\draw (\j,6) -- (\j,6.5);
\draw (1,6) -- (1,6.5);
};
\end{tikzpicture}
\hspace*{\fill}
\caption{The Garside element $\garside[6]$ of $B_{6}$.}\label{fig:garside}
\end{figure}
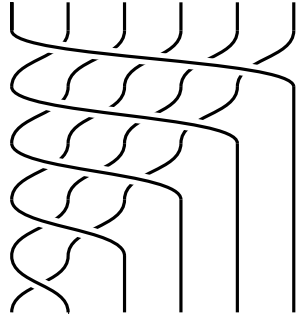

In terms of the presentation of \reth{presbn}, the `full twist' braid $\ft$ of $B_{n}$ is defined by:
\begin{equation}\label{eq:fulltwist}
\ft=(\sigma_1\cdots\sigma_{n-1})^n\in B_n.
\end{equation}
It has a special rôle in the theory of Artin braid groups. Since $\tau_{n}(\sigma_1\cdots\sigma_{n-1})=(1,n,n-1,\ldots,2)$, we see that $\ft$ belongs to $P_n$, and in terms of the generators of $P_n$ of \req{defaij}, one may check that:
\begin{equation*}
\ft=(A_{1,2})(A_{1,3}A_{2,3})\cdots (A_{1,n}A_{2,n}\cdots A_{n-1,n}).
\end{equation*}
The parenthesised terms in this expression commute pairwise~--~geometrically, this is obvious. This braid is the square of the well-known \emph{Garside element} (or `half-twist') $\garside$ of $B_n$ (see Figure~\ref{fig:garside}), defined by:
\begin{equation*}
\garside=(\sigma_1 \sigma_2\cdots \sigma_{n-1})(\sigma_1 \sigma_2\cdots \sigma_{n-2})\cdots (\sigma_1 \sigma_2)(\sigma_1).
\end{equation*}
The notion of Garside element is important in the study of braid groups, notably in the resolution of the conjugacy problem in $B_n$~\cite{Bi2,Ga}, and in a more general setting, in the theory of Garside groups and monoids~\cite{DDGKM,KT}. By~\cite[Lemma~2.5.1]{Bi2}, we have:
\begin{equation}\label{eq:conjgarside}
\garside \sigma_{i} \garside^{-1}= \sigma_{n-i} \quad \text{for all $1\leq i\leq n-1$,}
\end{equation}
and the $n\up{th}$ root $\sigma_1\cdots\sigma_{n-1}$ of $\ft$ that appears in \req{fulltwist} satisfies~\cite{A1,GG12,Mo}:
\begin{gather}
(\sigma_1\cdots\sigma_{n-1}) \sigma_{i} (\sigma_1\cdots\sigma_{n-1})^{-1}= \sigma_{i+1} \quad\text{for all $1\leq i\leq n-2$, and}\label{eq:alpha0a}\\
(\sigma_1\cdots\sigma_{n-1})^2 \sigma_{n-1} (\sigma_1\cdots\sigma_{n-1})^{-2}.\label{eq:alpha0b}
\end{gather}
From equations~\reqref{alpha0a}--\reqref{alpha0b}, it follows that $\ft$ commutes with all of the generators of \reth{presbn}, and so belongs to the centre of $B_{n}$ and of $P_{n}$. A straightforward argument using the short exact sequences~\reqref{permseq} and~\reqref{sesgen} with $M=\R^2$ enables one to show that $\ft$ generates the centre of the Artin braid groups. 

\begin{thm}[Chow~\cite{Ch}]\label{th:centre}
Let $n\geq 3$. Then $Z(B_n)=Z(P_n)=\ang{\ft}$.
\end{thm}

\begin{rem}\label{rem:centriv}
From \resec{present}, we know that $B_1=P_1=\brak{1}$, and $B_2$ and $P_2$ are infinite cyclic: $Z(B_2)=\ang{\sigma_1}$, and $Z(P_2)=\ang{A_{1,2}}=\ang{\ft[2]}$.
\end{rem}

A small number of surface braid groups possess non-trivial centre:
\begin{enumerate}[(a)]
\item If $M=\St$ (resp.\ $M=\rp$) and $n\geq 2$ then $Z(B_{n}(M))$ is cyclic of order $2$~\cite{GVB,M,vB} (see \resec{sphere}).
\item Let $\mathbb{T}^2$ denote the $2$-torus. Then $Z(B_n(\mathbb{T}^2))$ is free Abelian of rank~$2$~\cite{Bi2,PR}.
\end{enumerate}

Apart from these cases and a few other exceptions, most surface braid groups have trivial centre. With the aid of \reco{torsion}, one may once more use the short exact sequences~\reqref{permseq} and~\reqref{sesgen} to prove the following: 
\begin{prop}[\cite{GG2,PR}]\label{prop:nocentre}
Let $M$ be a compact surface different from the disc and the sphere whose fundamental group has trivial centre. Then for all $n\geq 1$, $Z(B_n(M))$ is trivial.
\end{prop}

\begin{rem}
The only compact surfaces whose fundamental group does not have trivial centre are the real projective plane, the annulus, the torus, the M\"obius band and the Klein bottle. 
\end{rem}

\subsection{Embeddings of surface braid groups}\label{sec:embeddings}

One possible approach in the study of surface braid groups is to determine relationships between braid groups of different surfaces. The first result in this direction is the embedding of $P_{n}$ in $P_{n}(M)$ given by \repr{birembed} and its extensions to non-orientable surfaces and to the full braid groups. The proof of \repr{birembed} uses induction and the Fadell-Neuwirth short exact sequences~\reqref{fnses} and~\reqref{sesgen}. 

Let $N$ be a subsurface of $M$, and let $m\geq 0$. Paris and Rolfsen studied the homomorphism $B_{n}(N)\to B_{n+m}(M)$ of braid groups induced by inclusion of $N$ in $M$, and gave necessary and sufficient conditions for it to be injective~\cite{PR}. In another direction, it is reasonable to ask whether it is possible to obtain embeddings of braid groups of surfaces that are not induced by inclusions (see \cite[page 216, Problem 1]{Bi2} for example). The answer is affirmative in the case of covering spaces:
\begin{thm}[\cite{GG11}]\label{th:covering}
Let $M$ be a compact, connected surface, possibly with a finite set of points removed from its interior. Let $d,n\in \N$, and let $\widetilde{M}$ be a $d$-fold covering space of $M$.  Then the covering map induces an embedding of the $n\up{th}$ braid group $B_{n}(M)$ of $M$ in the $dn\up{th}$ braid group $B_{dn}(\widetilde{M})$ of $\widetilde{M}$.
\end{thm}

To prove \reth{covering}, note that the inverse image of the covering map induces a map between the permuted configuration spaces of $M$ and $\widetilde{M}$. By restricting first to $F_{n}(M)$, one shows that this map induces the embedding mentioned in the statement of \reth{covering}. Note however that the embedding does not restrict to the corresponding pure braid groups: the image of $P_{n}(M)$ is a subgroup of the `mixed' subgroup $\pi_{1}\bigl( F_{dn}(\widetilde{M})/(S_{2}\times \cdots \times S_{2})\bigr)$ that is not contained in $P_{dn}(\widetilde{M})$. Although the map in question appears at first sight to be very natural, to our knowledge, it does not seem to have been studied previously in the literature. \reth{covering} should prove to be useful in the analysis of the structure of surface braid groups. As examples of this, one may deduce the linearity of the braid and mapping class groups of $\rp$ (see \resec{linear}), and one may classify their finite subgroups (see \resec{finitesubgp}). The following is an immediate consequence of \reth{covering}:
\begin{cor}\label{cor:embed}
Let $n\in \N$. The $n\up{th}$ braid group of a non-orientable surface embeds in the $2n\up{th}$ braid group of its orientable double covering. In particular, $B_{n}(\rp)$ embeds in $B_{2n}(\St)$.
\end{cor}
Using the covering map, one may write down explicitly the images in $B_{2n}(\St)$ of elements of $B_{n}(\rp)$. In this case, we see once more that such an embedding does not restrict to an embedding of the corresponding pure braid subgroups since if $n\geq 2$, $P_{n}(\rp)$ has torsion $4$ (see \repr{agt}(\ref{it:agt3})), while $P_{2n}(\St)$ has torsion $2$ (see \repr{splitpns}). \reco{embed} (and \reth{covering} in a more general context) would seem to be a significant step towards the resolution of the problem of Birman mentioned above concerning the relationship between the braid groups of a non-orientable surface and those of its orientable double covering.

\subsection{Braid combing and the splitting problem}\label{sec:combing}

Let $n\geq 2$ and $M=\R^2$, and consider the short exact sequence~\reqref{sesgen}:
\begin{equation}\label{eq:sesartin}
1 \to \F[n-1] \to P_n \stackrel{p_{n\#}}{\to} P_{n-1}\to 1,
\end{equation}
where we set $p_{n\#}=(p_{n,n-1})_{\#1}$ and we identify $\F[n-1]$ naturally with the free group $\ker{p_{n\#}}\cong\pi_1(\R^2\setminus\mathcal{Q}_{n-1},x_{n})$, where $\brak{x_{n}}=\mathcal{Q}_{n} \setminus \mathcal{Q}_{n-1}$. Recall that geometrically, $p_{n\#}$ `forgets' the $n\th$ string of a braid in $P_n$, and using \repr{hansen}, it may be seen easily that $p_{n\#}$ admits a section $\map{s_{n\#}}{P_{n-1}}[P_n]$ given geometrically by adding a vertical string (in terms of the generators of \repr{hansen}, $s_{n\#}$ maps $A_{i,j}$, $1\leq i<j\leq n-1$, considered as an element of $P_{n-1}$ to $A_{i,j}$, considered as an element of $P_n$). It follows that $P_n$ is isomorphic to the semi-direct product $\F[n-1]\rtimes_{\phi} P_{n-1}$, where the action $\phi$ is given by conjugation via $s_{n\#}$. By induction on $n$, $P_n$ may be written as an iterated semi-direct product of free groups, known as the \emph{Artin normal form}:
\begin{equation}\label{eq:combing}
P_n \cong \F[n-1]\rtimes \F[n-2]\rtimes \cdots \rtimes \F[2]\rtimes \F[1].
\end{equation}
The procedure for obtaining the Artin normal form of a pure braid $\beta$ is known as \emph{Artin combing}, and involves writing $\beta$ in the form $\beta=\beta_{n-1}\cdots \beta_1$, where $\beta_i\in \F[i]$. Since this expansion is unique and the word problem in free groups is soluble, this yields a (finite) algorithm to solve the word problem in $P_n$. Furthermore, $P_n$ is of finite index in $B_n$, and it is then an easy matter to solve the word problem in $B_n$ also. The decomposition~\reqref{combing} is one of the fundamental results in classical braid theory, and is frequently used to prove assertions about $P_{n}$ by induction, such as the study of the lower central series and the residual nilpotence of $P_n$~\cite{FR}, the bi-orderability of $P_n$ (see \reth{pnbiord}) and the fact that $P_{n}$ is poly-free (see \resec{fjconjasp}). Another application is obtained by taking $M=\R^2$ and $r=2$ in the short exact sequence~\reqref{sesgennr}, and using \reth{centre} and the fact that the projection $\map{(p_{n,2})_{\#}}{P_{n}}[P_{2}]$ sends $\ft$ to the generator $\ft[2]$ of $P_{2}$:

\begin{prop}[\cite{GG2}]\label{prop:splitpn}
Let $n\geq 3$. Then $P_n\cong P_{n-2}(\R^2\setminus \mathcal{Q}_{2}) \times \Z$.
\end{prop}

The problem of deciding whether a decomposition of the form~\reqref{combing} exists for surface braid groups is thus fundamental. This was indeed a recurrent and central question during the foundation of the theory and its subsequent development during the 1960's~\cite{Bi1,Fa,FaN,FvB,vB}. An interesting and natural question, to which we shall refer henceforth as \emph{the splitting problem}, is that of whether the short exact sequences~\reqref{fnses}--\reqref{fnsurface} split. Clearly, the existence of a geometric cross-section on the level of configuration spaces implies that of a section on the algebraic level, and in most cases the converse is true. Indeed, if $M$ is aspherical, this follows from~\cite{Baue,Wh1}, while if $M=\St$ or $\rp$, one may consult~\cite{GG3,GG4}. We sum up the situation as follows.
\begin{prop}\label{prop:xsplit}
Let $M$ be a compact, connected surface (so $m=0$ in \req{fnses}). Let $1\leq r<n$, and suppose that $r\geq 3$ if $M=\St$ and $r\geq 2$ if $M=\R P^2$. Then the Fadell-Neuwirth fibration $\map{p_{n,r}}{F_{n}(M)}[F_r(M)]$ (resp.\ $\map{q_{n,r}}{F_{n}(M)/(\sn[r]\times \sn[n-r])}[D_r(M)]$)
admits a cross-section if and only if the short exact sequence~\reqref{fnses} (resp.~\reqref{fnsurface}) splits.
\end{prop}


In the case of the pure braid groups, the splitting problem for~\reqref{fnses} has been studied for other surfaces besides the plane. Fadell and Neuwirth gave various sufficient conditions for the existence of a geometric section for $p_{n,r}$~\cite{FaN}. If $m\geq 1$ (or if $\partial M\neq \vide$) then $p_{n,r}$ always admits a cross-section, and hence $(p_{n,r})_{\#1}$ does too~\cite{GG1,GG7}. So suppose that $m=0$. If $M=\St$ and $r\geq 3$, $p_{n,r}$ admits a cross-section~\cite{FvB}, and thus the short exact sequence~\reqref{sesgennr} splits. In the case $M=\rp$, Van Buskirk showed that the fibration $p_{3,2}$ admits a cross-section~\cite{vB} (and hence so does the corresponding homomorphism $(p_{3,2})_{\#1}$), but that for $n\geq 2$, neither the fibration $p_{n,1}$ nor the homomorphism $(p_{n,1})_{\#1}$ admit a section (this is one of the cases not covered by \repr{xsplit}), this being a consequence of the fact that $\rp$ has the fixed point property. If $M$ is the $2$-torus then Birman exhibited an explicit cross-section for $p_{n,n-1}$ if $n\geq 2$~\cite{Bi1}, which implies that the short exact sequence~\reqref{sesgen} splits for all $n$. This implies that~\reqref{sesgennr} splits for all $1\leq r<n$. In the case of orientable surfaces without boundary of genus at least two, the question of the splitting of~\reqref{sesgen} was posed explicitly by Birman in 1969~\cite{Bi1}, and was finally answered in~\cite{GG1}:
\begin{thm}[\cite{GG1}]\label{th:fnsplit}
If $M$ is a compact orientable surface without boundary of genus $g\geq 2$, the short exact sequence~\reqref{sesgennr} splits if and only if $r=1$.
\end{thm}

For the remaining cases, the problem was studied in a series of papers~\cite{GG3,GG4,GGgeom}, and a complete solution to the splitting problem for~\reqref{sesgennr} was given in~\cite{GG7}:
\begin{thm}[\cite{GG7}]\label{th:complete}
Let $1\leq r<n$ and $m\geqslant 0$, and let $M$ be a connected surface.
\begin{enumerate}[(a)]
\item If $m>0$ or if $M$ has non-empty boundary then $(p_{n,r})_{\#1}$ admits a section.
\item Suppose that $m=0$ and that $\partial M=\vide$. Then $(p_{n,r})_{\#1}$ admits a section if and only if one of the following conditions holds:
\begin{enumerate}[(i)]
\item $M=\St$, the $2$-torus $\mathbb{T}^2$ or the Klein bottle $\mathbb{K}^2$.
\item $M=\rp$, $n=3$ and $r=2$.
\item $M\neq \rp,\St,\mathbb{T}^2,\mathbb{K}^2$ and $r=1$.
\end{enumerate}
\end{enumerate}
\end{thm}

To obtain a positive answer to the splitting problem, it suffices of course to exhibit an explicit section. However, in general it is very difficult
to prove directly that the (generalised) Fadell-Neuwirth short exact sequences do not split. One of the principal methods that was used in the proof of \reth{complete} is based on the following observation: let $G$ be a group, and let $K,H$ be normal subgroups of $G$ such that $H$ is contained in $K$. If the extension $1\to K\to G\to Q\to 1$ splits then so does the extension $1\to K/H\to G/H\to Q\to 1$. The condition on $H$ is satisfied for example if $H$ is an element of either the lower central series $(\Gamma_{i}(K))_{i\in \N}$ or of the derived series of $K$. In many parts of the proof of \reth{complete}, it suffices to take $H=\Gamma_{2}(K)$, in which case $K/H$ is the Abelianisation of $K$, to show that this second extension does not split, which then implies that the first extension does not split.

From the point of view of the splitting problem, it is thus helpful to know the lower central and derived series of the braid groups occurring in these group extensions. These series have been calculated in many cases~\cite{BGG,BGG2,GG6,GGjktr2,GG10,GL}. The splitting problem for the generalised Fadell-Neuwirth short exact sequence~\reqref{fnsurface} has been studied in the case $M=\St$~\cite{GG4}.

\subsection{Homotopy type of the configuration spaces of $\St$ and $\rp$ and periodicity}\label{sec:homotype}

As we saw in \reth{fnses}(\ref{it:eilenmac}), the configuration spaces of surfaces different from $\St$ and $\rp$ are Eilenberg-Mac~Lane spaces of type $K(\pi,1)$. For the two exceptional cases of $\St$ and $\rp$, the situation is very different, and in view of the relation with the homotopy groups of $\St$ (and $\St[3]$), motivates the study of their configuration spaces. In the case of $\St$, the following proposition may be found in~\cite{BCP,FZ}. An alternative proof was given in~\cite{GG12}.
\begin{prop}[\cite{BCP,FZ}]\label{prop:homot}\mbox{}
\begin{enumerate}[(a)]
\item\label{it:homottype2} The space $F_2(\St)$ (resp.\ $D_2(\St)$) has the homotopy type of $\St$ (resp.\ of $\rp$). Hence the universal covering space of $D_2(\St)$ is $F_2(\St)$.
\item If $n\geq 3$, the universal covering space of $F_n(\St)$ or of $D_n(\St)$ has the homotopy type of the $3$--sphere $\St[3]$.
\end{enumerate}
\end{prop}

A similar result holds for the configuration spaces of $\rp$:
\enlargethispage{3mm}

\begin{prop}[{\cite{GG3}}]\mbox{}\label{prop:homotopy}
\begin{enumerate}
\item \label{it:f1rp} The universal covering of $F_1(\rp)$ is $\St$.
\item \label{it:fnrp} For $n\geq 2$, the universal covering space of $F_n(\rp)$ or of $D_{n}(\rp)$ has the homotopy type of $\St[3]$. 
\end{enumerate}
\end{prop}

Suppose that $n\geq 3$ if $M=\St$ and that $n\geq 2$ if $M=\rp$. From Propositions~\ref{prop:homot} and~\ref{prop:homotopy}, the universal covering space $X$ of $F_{n}(M)$ is a finite-dimensional complex that has the homotopy type of $\St[3]$. Thus any finite subgroup of $B_n(M)$ acts freely on $X$, and so has period $2$ or $4$ by \cite[Proposition~10.2, Section~10, Chapter~VII]{Br}. It thus follows that such a subgroup must be one of the subgroups that appear in the Suzuki-Zassenhaus classification of periodic groups~\cite{AM}. We shall come back to the finite subgroups of $B_{n}(M)$ in \resec{finitesubgp}. Using results of~\cite[Section~2]{AS} allows one to obtain a periodicity result for \emph{any} subgroup of $B_{n}(M)$:
\begin{prop}[\cite{GG12}]\label{prop:per24}
Let $M=\St$ or $\rp$, let $n\geq 3$ if $M=\St$ and $n\geq 2$ if $M=\rp$, and let $G$ be a group abstractly isomorphic to a subgroup of $B_{n}(M)$. Then there exists $r_{0}\geq 1$ such that $H^r(G; \Z) \cong H^{r+4}(G; \Z)$ for all $r\geq r_{0}$. 
\end{prop} 

The connections between surface braid groups and the homotopy groups of $\St$ do not end there. If $M$ is a surface, recall that an element of $P_{n}(M)$ is said to be \emph{Brunnian} if it becomes trivial after removing any one of its $n$ strings. The subgroup $\operatorname{Brun}_{n}(M)$ of Brunnian braids may thus be seen to be the intersection $\bigcap_{i=1}^n \ker{\map{d_{i}}{P_{n}(M)}[P_{n-1}(M)]}$, where $d_{i}$ corresponds geometrically to removing the $i\th$ string. The study of the homomorphisms $d_{i}$ allows one to introduce a simplicial structure on the pure braid groups of $M$. In this way, the following result was proved in~\cite{BCWW}:
\begin{thm}[\cite{BCWW}]\label{th:bcww}
Let $n\geq 4$. Then there is an exact sequence of the form:
\begin{equation*}
1\to  \operatorname{Brun}_{n+2}(\St)\to \operatorname{Brun}_{n+1}(\dt) \to \operatorname{Brun}_{n+1}(\St) \to \pi_{n}(\St)\to 1.
\end{equation*}
\end{thm}

\reth{bcww} has been generalised in some sense to $\rp$ in~\cite{bmvw}, and to other surfaces in~\cite{Oc}. The hope is that one might understand better the homotopy groups of $\St$ using the structure of Brunnian braid groups.

\subsection{Orderability}\label{sec:orderable}

A group $G$ is said to be \emph{left orderable} (resp.\ \emph{right orderable}) if it admits a total ordering $<$ that is invariant under left (resp.\ right) multiplication in $G$. In other words,
\begin{equation*}
\forall x,y,z \in G, \quad x<y \Longrightarrow zx<zy \quad \text{(resp.\ $x<y \Longrightarrow xz<yz$).}
\end{equation*}
Any left ordering may be converted into a right ordering by considering inverses of elements, but the two orderings will in general be different. A group is said to be \emph{biorderable} if there exists a total ordering $<$ for which $G$ is both right and left orderable. The classes of left orderable and biorderable groups are closed under subgroups, direct products and free products (so free groups are biorderable), and that the class of left orderable groups is also closed under extensions. It is an easy exercise to show that a left orderable group is torsion free. Further, a biorderable group has no generalised torsion (a group $G$ is said to have \emph{generalised torsion} if there exist $g,h_1,\ldots,h_k\in G$, $g\neq 1$, such that $h_1gh_1^{-1}\cdots h_kgh_k^{-1}=1$). 

One of the most exciting developments over the past twenty years in the theory of braid groups is the discovery of Dehornoy~\cite{Deh}, using some deep results in set theory, that $B_{n}$ is left orderable:
\begin{thm}[Dehornoy~\cite{Deh,DDRW,Ka,KT}]\label{th:order}
$B_n$ is left orderable.
\end{thm} 

\reth{order} thus yields an alternative proof of \reco{torsion}, that is, $B_n$ is torsion free. In the wake of Dehornoy's paper, a group of topologists came up with a different way of interpreting his ordering of $B_n$ in terms of $\mathcal{MCG}(\dt,n)$~\cite{FGRW}. Short and Wiest described another approach due to Thurston using the action of the mapping class group on the hyperbolic plane which in fact defines uncountably many different orderings on $B_n$~\cite{SW}. The reader is referred to the monograph~\cite{DDRW} for a full description of these different points of view, as well as to~\cite[Chapter~7]{KT}. These results have led to renewed interest in orderable groups, notably in the case of $3$-manifold groups~\cite{BRW}.

If $n\geq 3$ then $B_n$ is not biorderable since it has generalised torsion. Indeed, by \req{conjgarside}, we have $\garside (\sigma_{n-1}^{-1}\sigma_{1}) \garside^{-1}=(\sigma_{n-1}^{-1}\sigma_{1})^{-1}$. However:
\begin{thm}[Falk and Randell, Kim and Rolfsen~\cite{DDRW,FaRa,KT,KR,RZ}]\label{th:pnbiord}
$P_n$ is biorderable.
\end{thm}

Falk and Randell's result is a consequence of the residual nilpotence of $P_{n}$, and the fact that its lower central series quotients are torsion free. Kim and Rolfsen's proof gives an explicit biordering, and makes use of \req{sesartin} and an ordering emanating from the Magnus expansion of free groups.

Theorems~\ref{th:order} and~\ref{th:pnbiord} motivated the study of the (bi)orderability of surface braid groups. We summarise the known results as follows.

\begin{enumerate}[(a)]
\item Since the braid groups of the $\St$ and $\rp$ have torsion (see \rerem{centriv} and \resec{sphere}), they are not left orderable. 
\item As was pointed out in~\cite{RolW}, the short exact sequence~\reqref{sesgen} implies that the braid groups of any compact surface different from $\St$ and $\rp$ are left orderable. Pure braid groups of compact, orientable surfaces without boundary of genus $g\geq 1$ are biorderable~\cite{GM2}: the proof makes use of the short exact sequence~\reqref{goldberg}. Pure braid groups of compact, non-orientable surfaces without boundary of genus $g\geq 2$ have generalised torsion, and so are not biorderable, but are left orderable~\cite{GM2}. 
\item If $n\geq 3$ and $M$ is a compact surface different from $\St$ and $\rp$ then the generalisation of \repr{birembed} to $B_{n}(M)$ and the fact that $B_{n}$ is not biorderable imply that $B_n(M)$ is not biorderable.  Using \req{mcggen} and the fact that mapping class groups of surfaces with non-empty boundary are left orderable~\cite{RouW}, it follows that the braid groups of any surface with boundary are left orderable. If $M$ is without boundary and $n\geq 2$ then it seems to be an open question as to whether $B_n(M)$ is left orderable.
\end{enumerate}

\subsection{Linearity}\label{sec:linear}

A group is said to be \emph{linear} if it admits a faithful representation in a multiplicative group of matrices over some field. The linearity of the braid groups is a classical problem (see \cite[page~220, Problem~30]{Bi1} and \cite[Question~1]{Bar2} for example). Krammer~\cite{Kr1,Kr2} and Bigelow~\cite{Big2} showed that $B_{n}$ is linear. The question of the linearity of surface braid groups has been the subject of various papers during the last few years~\cite{Bar,Bar2,Big2,BB,Kor}.  The linearity of $\operatorname{\mathcal{MCG}}(\St,n)$ was proved in~\cite{Bar,Bar2,BB,Kor}, and that of $B_{n}(\St)$ was obtained in~\cite{Bar,Bar2,BB}. If $n=1$ then we are in the case of surface groups, which are known to be linear for any surface $M$. If $n\leq 2$ then $B_{n}(\rp)$ is linear because it is finite, while $B_{3}(\rp)$ is known to be isomorphic to a subgroup of $\operatorname{\mathrm{GL}}(96,\Z)$~\cite{Bar2}. With the help of \reco{embed} and the short exact sequence~\reqref{mcg}, we have the following results.

\pagebreak
\begin{thm}[\cite{GG11}]\label{th:linear}
Let $n\in \N$.
\begin{enumerate}[(a)]
\item Let $M$ be a compact, connected surface, possibly with boundary, of genus zero if $M$ is orientable, and of genus one if $M$ is non-orientable. Then $B_{n}(M)$ is linear. 
\item The mapping class groups $\operatorname{\mathcal{MCG}}(\rp,n)$ are linear.
\item Let $\mathbb{T}^2$ denote the $2$-torus, and let $x\in \mathbb{T}^2$. Then $B_{n+1}(\mathbb{T}^2)$ is linear if and only if $B_{n}(\mathbb{T}^2\setminus\brak{x})$ is linear. Consequently, $B_{2}(\mathbb{T}^2)$ is linear.
\end{enumerate}
\end{thm}

In particular, the braid groups of $\rp$ and the M\"obius band are linear. To our knowledge, very little is known about the linearity of braid groups of other surfaces.
 
\section{Braid groups of the sphere and the projective plane}\label{sec:virtually}

Together with the braid groups of $\rp$, the braid groups of $\St$ are of particular interest, notably because they have non-trivial centre (see \repr{nocentres}), and torsion (see \reth{murasugi}).  In \resec{sphere}, we begin by recalling some of their basic properties, including the characterisation of their torsion elements. In \resec{finitesubgp}, we give the classification of the isomorphism classes of their finite subgroups, and in \resec{vcyclics}, this is extended to the isomorphism classes of the virtually cyclic subgroups of their pure braid groups and of $B_{n}(\St)$. As well as being interesting in their own right, these results play an important rôle in the determination of the lower algebraic $K$-theory of the group rings of the braid groups of these two surfaces (see~\resec{ktheory}). From this point of view, it is also necessary to have a good understanding of the conjugacy classes of the finite order elements and the finite subgroups.

\subsection{Basic properties}\label{sec:sphere}

In this section, we recall briefly some of the basic properties of the braid groups of $\St$ and $\rp$. We first consider $B_{n}(\St)$. The reader may consult~\cite{Fa,FvB,GVB,GG2,vB} for more details. Consider the group homomorphism $\map {j_{\#}}{B_n}[B_n(\St)]$ of \resec{defconfig} induced by an inclusion $\map{j}{\dt}[\St]$. If $\beta\in B_n$ then we shall denote its image $j_{\#}(\beta)$ simply by $\beta$. A presentation of $B_n(\St)$ is as follows:
\begin{thm}[\cite{FvB}]\label{th:presfvb}
The following constitutes a presentation of the group $B_n(\St)$:
\begin{enumerate}
\item[\underline{generators:}] $\sigma_1,\ldots,\sigma_{n-1}$.
\item[\underline{relations:}]\mbox{}
\begin{enumerate}[(i)]
\item relations~\reqref{Artin1} and~\reqref{Artin2}.
\item the `surface relation' of $B_n(\St)$:
\begin{equation}\label{eq:surface}
\text{$\sigma_1\cdots \sigma_{n-2}\sigma_{n-1}^2 \sigma_{n-2}\cdots \sigma_1=1$.}
\end{equation}
\end{enumerate}
\end{enumerate}
\end{thm}

The surface relation may be seen geometrically to indeed represent the trivial element of $B_{n}(\St)$ (see~\cite[page~194]{MK} for example). It follows from \reth{presfvb} that $B_n(\St)$ is a quotient of $B_n$, and that its Abelianisation is isomorphic to $\Z_{2(n-1)}$. The first three sphere braid groups are finite: $B_1(\St)$ is trivial, $B_2(\St)$ is cyclic of order~$2$, and $B_3(\St)$ is isomorphic to $\Z_3 \rtimes \Z_4$, the action being the non-trivial one. For $n\geq 4$, $B_n(\St)$ is infinite. Just as for the Artin braid groups, the full twist braid of $B_n(\St)$ plays an important part, and has some interesting additional properties.

\pagebreak
\begin{prop}[\cite{GVB,GG2}]\label{prop:nocentres}
Let $n\geq 3$. Then:
\begin{enumerate}[(a)]
\item\label{it:propsbns2} $\ft$ is the unique element in $P_n(\St)$ of finite order, and is the unique element of $B_n(\St)$ of order $2$.
\item\label{it:propsbns1} $\ft$ generates the centre $Z(B_n(\St))$ of $B_n(\St)$.
\end{enumerate}
\end{prop}

Taking $M=\St$, $m=0$ and $r=3$ in the short exact sequence~\reqref{fnses} and applying an argument similar to that used in the proof of \repr{splitpn} yields:
\begin{prop}[\cite{GG2}]\label{prop:splitpns}
Let $n\geq 4$. Then $P_n(\St)\cong P_{n-3}(\St\setminus \mathcal{Q}_{3}) \times \Z_2$.
\end{prop}

From this and \repr{nocentre}, it follows that $Z(P_n(\St))=\ang{\ft}$ for all $n\geq 4$.

Let $n\geq 3$. Fadell and Van Buskirk showed that the element $\alpha_{0}= \sigma_1\cdots \sigma_{n-2} \sigma_{n-1}$ is of order $2n$ in $B_{n}(\St)$~\cite{FvB}. Gillette and Van Buskirk later proved that if $k\in \N$ then $B_n(\St)$ has an element of order $k$ if and only if $k$ divides one of $2n$, $2(n-1)$ or $2(n-2)$~\cite{GVB}. Using Seifert fibre space theory, Murasugi characterised the finite order elements of $B_n(\St)$ and $B_n(\rp)$. In the case of the sphere, $B_{n}(\St)$, up to conjugacy and powers, there are precisely three torsion elements:
\begin{thm}[\cite{M}]\label{th:murasugi}
Let $n\geq 3$. Then the torsion elements of $B_n(\St)$ are precisely the conjugates of powers of the three elements $\alpha_0$, $\alpha_1=\sigma_1\cdots \sigma_{n-2} \sigma_{n-1}^2$ and $\alpha_2=\sigma_1\cdots \sigma_{n-3} \sigma_{n-2}^2$, which are of order $2n$, $2(n-1)$ and $2(n-2)$ respectively. 
\end{thm}
\reth{murasugi} implies Gillette and Van Buskirk's result, and in conjunction with \repr{nocentres}(\ref{it:propsbns2}), yields the useful relation:
\begin{equation}\label{eq:uniqueorder2}
\ft=\alpha_{i}^{n-i} \quad \text{for all $i\in\brak{0,1,2}$,}
\end{equation}
which implies that $\alpha_{i}$ is an $(n-i)\th$ root of $\ft$.  Since the permutation $\tau_{n}(\alpha_{i})$ consists of an $(n-i)$-cycle and $i$ fixed elements, we see that the $\alpha_{i}$ are pairwise non conjugate. One interesting fact about the group $B_n(\St)$ is that it is generated by $\alpha_0$ and $\alpha_1$~\cite{GG4}, and so is torsion generated in the sense of~\cite{GG13}. Equations~\reqref{alpha0a}--\reqref{alpha0b} also hold in $B_{n}(\St)$, and more generally, for $i\in \brak{0,1,2}$ we have~\cite{GG12}:
\begin{gather}
\alpha_{i}^l \sigma_{j} \alpha_{i}^{-l}=\sigma_{j+l} \quad \text{for all $j,l\in \N$ satisfying $j+l\leq n-i-1$,}\label{eq:fundaa}\\
\sigma_{1}=\alpha_{i}^2 \sigma_{n-i-1} \alpha_{i}^{-2}\label{eq:fundab}
\end{gather}
in $B_{n}$ and so also in $B_{n}(\St)$, in other words, conjugation by $\alpha_{i}$ permutes the $n-i$ elements $\sigma_{1},\ldots, \sigma_{n-i-1}, \alpha_{i} \sigma_{n-i-1} \alpha_{i}^{-1}$ cyclically. These relations prove to be very useful in the study of the finite and virtually cyclic subgroups of $B_{n}(\St)$.

We now turn to the braid groups of the projective plane. Some basic references are~\cite{GG2,GG4,GG9,GG10,vB}. We first recall a presentation of $B_{n}(\rp)$ due to Van Buskirk~\cite{vB}:
\enlargethispage{3mm}
\begin{thm}[\cite{vB}]\label{th:presvb}
The following constitutes a presentation of the group $B_n(\rp)$:
\begin{enumerate}
\item[\underline{generators:}] $\sigma_1,\ldots,\sigma_{n-1},\rho_1,\ldots,\rho_n$.
\item[\underline{relations:}]\mbox{}
\begin{enumerate}[(i)]
\item relations~\reqref{Artin1} and~\reqref{Artin2}.
\item $\sigma_i\rho_j=\rho_j\sigma_i$ for $j \ne i, i+1$.
\item\label{it:rel3rp2} $\rho_{i+1}=\sigma_i^{-1}\rho_i\sigma_i^{-1}$ for $1\leq i \leq n-1$.
\item $\rho_{i+1}^{-1}\rho_i^{-1}\rho_{i+1}\rho_i=\sigma_i^2$  for $1\leq i \leq n-1$.
\item $\rho_1^2=\sigma_1\sigma_2\cdots\sigma_{n-2}\sigma_{n-1}^2\sigma_{n-2}\ldots\sigma_2\sigma_1$.
\end{enumerate}
\end{enumerate}
\end{thm}

Each of the generators $\rho_{i}$ corresponds geometrically to an element of the fundamental group of $\rp$ based at the $i\th$ basepoint. A presentation of $P_{n}(\rp)$ was given in~\cite{GG4}. From these presentations, we see that the first two braid groups of $\rp$ are finite: $B_1(\rp)=P_1(\rp)\cong \Z_{2}$, $P_{2}(\rp)$ is isomorphic to the quaternion group $\quat$ of order~$8$, and $B_{2}(\rp)$ is isomorphic to the generalised quaternion group of order $16$~\cite{vB}. For $n\geq 3$, $B_{n}(\rp)$ is infinite. If $n\geq 2$, the Abelianisation of $B_{n}(\rp)$ is $\Z_{2}^2$, while that of $P_{n}(\rp)$ is $\Z_{2}^n$. If $M=\rp$ and $m=0$, the map $p_{3,2}$ of \req{fnles} admits a geometric section given by taking the vector product of two directions, and so by \req{sesgennr}, $P_{3}(\rp)$ is isomorphic to a semi-direct product of a free group of rank $2$ by $\quat$~\cite{vB}; an explicit action was given in~\cite{GG2,GG10}.

We recall that the \emph{virtual cohomological dimension} of a group is equal to the (common)  cohomological dimension of its torsion-free subgroups of finite index~\cite[page~226]{Br}. As an application of the Fadell-Neuwirth short exact sequence~\reqref{sesgennr}, \repr{splitpns} and the fact that $P_{2}(\rp)\cong \quat$, one may compute the virtual cohomological dimension of the braid groups of $\St$ and $\rp$:
\begin{prop}[\cite{GG15}]
Let $M$ be equal to $\St$ (resp.\ $\rp$), and let $n\geq 3$ (resp.\ $n\geq 2$). Then the virtual cohomological dimension of $B_{n}(M)$ and of $P_{n}(M)$ is equal to $n-3$ (resp.\ $n-2$). 
\end{prop}

For $n\geq 2$, Murasugi showed that $\ft$ generates the centre of $B_n(\rp)$~\cite{M}. The following proposition summarises some other basic results concerning the torsion of the braid groups of $\rp$.

\begin{prop}[\cite{GG2,GG9}]\label{prop:agt}
Let $n\geq 2$. Then:
\begin{enumerate}[(a)]
\item\label{it:agt1} $B_{n}(\rp)$ has an element of order $k$ if and only if $k$ divides either $4n$ or
$4(n-1)$.
\item\label{it:agt3} the (non-trivial) torsion of $P_{n}(\rp)$ is precisely $2$ and $4$.
\item\label{it:agt2} the full twist $\ft$ is the unique element of $B_{n}(\rp)$ of order $2$.
\end{enumerate}
\end{prop}

If $M=\St$ or $\rp$, it follows from Propositions~\ref{prop:nocentres} and~\ref{prop:agt} that the kernel of the short exact sequence~\reqref{mcg} is generated by $\ft$. In~\cite[Proposition~26]{GG2}, it was proved that the following elements of $B_{n}(\rp)$:
\begin{align*}
a &= \sigma_{n-1}^{-1} \cdots \sigma_{1}^{-1} \rho_{1}\\
b &= \sigma_{n-2}^{-1} \cdots \sigma_{1}^{-1} \rho_{1}
\end{align*}
are of order $4n$ and $4(n-1)$ respectively. By~\cite[Remark~27]{GG2}, we have
\begin{equation}
\label{eq:defalpha}
\left\{ \begin{aligned}
\alpha &= a^n= \rho_{n} \cdots \rho_{1}\\
\beta &= b^{n-1}= \rho_{n-1} \cdots \rho_{1}.
\end{aligned} \right.
\end{equation}
It is clear that $\alpha$ and $\beta$ are pure braids of order $4$. The finite order elements of $B_{n}(\rp)$ had previously been characterised in~\cite{M}, but the results are less transparent than in the case of $\St$ given by \reth{murasugi}. For example, it is not clear what the orders of the given torsion elements are, even for elements of $P_{n}(\rp)$. In~\cite{GG8}, Murasugi's characterisation was simplified somewhat as follows.

\begin{thm}[\cite{GG8}]\label{th:murasugi2}
Let $n\geq 2$, and let $x\in B_{n}(\rp)$. Then $x$ is of finite order if and only if there exist $i\in \brak{1,2}$ and $0\leq r\leq n+1-i$ such that $x$ is a power of a conjugate of the following element:
\begin{equation}\label{eq:finiteordrp2}
(\rho_{r}\sigma_{r-1}\cdots \sigma_{1})^{2r/l} (\sigma_{r+1}\cdots \sigma_{n-1}\sigma_{r+1}^{i-1})^{p/l}
\end{equation}
where $p=(n+1-i)-r$ and $l=\gcd(p,2r)$. Further, this element is of order $2l$.
\end{thm}

Using \reth{presvb}, one may check that the element $a$ (resp.\ $b$) is one of the above elements by taking $r=n$ and $i=1$ (resp.\ $r=n-1$ and $i=2$). The permutation and the Abelianisation may be used to distinguish the conjugacy classes of the elements given by \req{finiteordrp2}. The following result gives information about the conjugacy classes  and the centralisers of elements of $P_{n}(\rp)$ of order $4$:
\begin{prop}[\cite{GG8}]\label{prop:order4}
Let $n\geq 2$, and let $x \in P_{n}(\rp)$ be an element of order $4$.
\begin{enumerate}[(a)]
\item\label{it:order4a} In $B_{n}(\rp)$, $x$ is conjugate to an element of $\brak{\alpha,\beta,\alpha^{-1},\beta^{-1}}$.
\item\label{it:order4b} The centraliser $Z_{P_{n}(\rp)}(x)$ of $x$ in $P_{n}(\rp)$ is equal to $\ang{x}$.
\end{enumerate}
\end{prop}

It was shown in~\cite{GG13} that if $n\geq 2$, there are $(n-2)! \, (2n-1)$ conjugacy classes of elements of order $4$ in $P_{n}(\rp)$ (there is a misprint in the statement of~\cite[Proposition~11]{GG13}, $B_{n}(\rp)$ should read $P_{n}(\rp)$). The analysis of the conjugacy classes of finite order elements of $B_{n}(\rp)$ is the subject of work in progress~\cite{GG14}.

The elements $a$ and $b$ have some interesting properties that mirror those of equations~\reqref{fundaa}--\reqref{fundab} that may be used to study the structure of $B_{n}(\rp)$. From~\cite[pages 777--778]{GG2}, conjugation by $a^{-1}$ permutes cyclically the elements of the following sets:
\begin{equation*}
\brak{\sigma_{1},\ldots, \sigma_{n-1}, a^{-1}\sigma_{n-1}a, \sigma_{1}^{-1},\ldots, \sigma_{n-1}^{-1}, a^{-1}\sigma_{n-1}^{-1}a} \quad\text{and}\quad \brak{\rho_{1},\ldots \rho_{n}, \rho_{1}^{-1},\ldots, \rho_{n}^{-1}},
\end{equation*}
and conjugation by $b^{-1}$ permutes cyclically the following elements:
\begin{equation*}
\sigma_{1},\ldots, \sigma_{n-2}, b^{-1}\sigma_{n-2}b, \sigma_{1}^{-1},\ldots, \sigma_{n-2}^{-1}, b^{-1} \sigma_{n-1}^{-1} b.
\end{equation*}
Note that there is a typographical error in line~16 of \cite[page~778]{GG2}: it should
read `\ldots shows that $b^{-2}\sigma_{n-2}b^2=\sigma_{1}^{-1}$ \ldots', and not `\ldots shows that
$b^{-2}\sigma_{n-1}b^2=\sigma_{1}^{-1}$ \ldots'. By~\cite[pages~865--866]{GG11}, we also have that:
\begin{equation}\label{eq:propgarside}
\quad \garside a \garside^{-1}=a^{-1} \quad \text{and} \quad (\garside a^{-1}) b (a \garside^{-1})=b^{-1} \quad\text{for all $i=1,\ldots,n-1$.} 
\end{equation}
As for $\St$, such relations are very useful in the study of the finite and virtually cyclic subgroups of $B_{n}(\rp)$.

\subsection{Finite subgroups of the braid groups of $\St$ and $\rp$}\label{sec:finitesubgp}

We start by considering the pure braid groups of $\St$ and $\rp$. In the case of $P_{n}(\St)$, there are only two finite subgroups for $n\geq 3$ by \repr{splitpns}, the trivial group $\brak{e}$ and that generated by the full twist $\ft$. In the case of $P_{n}(\rp)$, there are more possibilities: 
\begin{prop}[\cite{GG8}]\label{prop:finitepn}
Up to isomorphism, the maximal finite subgroups of $P_{n}(\rp)$ are:
\begin{enumerate}[(a)]
\item $\Z_{2}$ if $n=1$.
\item $\quat$ if $n=2,3$.
\item $\Z_{4}$ if $n\geq 4$.
\end{enumerate}
\end{prop}
As we mentioned above, in \repr{order4}, we know the numebr of conjugacy classes of the elements of $P_{n}(\rp)$ of order $4$, both in $P_{n}(\rp)$ and in $B_{n}(\rp)$.

We now turn to $B_{n}(\St)$ and $B_{n}(\rp)$. The results of \reth{murasugi} and \repr{agt} imply that we know the isomorphism classes of their finite cyclic subgroups. This leads naturally to the question as to which isomorphism classes of finite groups are realised as subgroups of these two groups. From~\cite{GG4}, if $n\geq 3$ then $B_n(\St)$ contains an isomorphic copy of the finite group $B_3(\St)$ of order $12$ if and only if $n\nequiv 1 \bmod 3$. During the study of the lower central series of $B_{n}(\St)$, it was observed that the commutator subgroup $\Gamma_2\left(B_4(\St) \right)$ of $B_{4}(\St)$ is isomorphic to a semi-direct product of $\quat$ by a free group of rank $2$~\cite{GG6} (see also~\cite{GJM}). The question of the realisation of $\quat$ as a subgroup of $B_n(\St)$ was posed explicitly by R.~Brown~\cite{ATD} in connection with the Dirac string problem and the fact that the fundamental group of $\operatorname{SO}(3)$ is isomorphic to $\Z_2$~\cite{Fa,Ha,Ne}. The existence of a subgroup of $B_{4}(\St)$ isomorphic to $\quat$ was studied by J.~G.~Thompson~\cite{Th}. It was shown in~\cite{GG5} that if $n\geq 3$, $B_n(\St)$ contains a subgroup isomorphic to $\quat$ if and only if $n$ is even. The construction of $\quat$ given in~\cite{GG5} may be generalised. If $m\geq 2$, let $\dic{4m}$ denote the \emph{dicyclic group} of order $4m$. It admits a presentation of the form: 
\begin{equation}\label{eq:presdic}
\setangr{x,y}{x^m=y^2,\; yxy^{-1}=x^{-1}}.
\end{equation}
If in addition $m$ is a power of $2$ then we will refer to the dicyclic group of order $4m$ as the \emph{generalised quaternion group} of order $4m$, and denote it by $\quat[4m]$. For example, if $m=2$ then we obtain the usual quaternion group $\quat$. For $i\in\brak{0,2}$, we have:
\begin{equation}\label{eq:basicconj}
\garside \alpha_{i}'\garside^{-1}=\alpha_{i}'^{-1}, \quad\text{where $\alpha_{i}'=\alpha_{0}\alpha_{i}\alpha_{0}^{-1}= \alpha_{0}^{i/2}\alpha_{i}\alpha_{0}^{-i/2}$,}
\end{equation}
and the group $\dic{4(n-i)}$ is realised in terms of the generators of $B_{n}(\St)$ by the subgroup $\ang{\alpha_{i}', \garside}$, which we shall call the \emph{standard copy} of $\dic{4(n-i)}$ in $B_{n}(\St)$~\cite{GG6,GG7}. Let $\tonestar$ (resp.\ $\oonestar$, $\istar$) denote the \emph{binary tetrahedral group} of order $24$ (resp.\ the \emph{binary octahedral group} of order $48$, the \emph{binary icosahedral group} of order $120$). 
The groups $\tonestar,\oonestar$ and $\istar$, to which we refer collectively as the \emph{binary polyhedral groups}, admit presentations of the form~\cite{Cox,CM}:
\begin{equation*}
\ang{p,3,2}=\setangr{A,B}{A^p=B^3= (AB)^2},
\end{equation*}
where $p=3,4,5$ respectively, and the element $A^p$ is central and is the unique element of order~$2$. The group $\tonestar$ also admits the following presentation~\cite[page~198]{Wo}:
\begin{equation}\label{eq:preststar}
\setangr{P,Q,X}{X^3=1,\, P^2=Q^2,\, PQP^{-1}=Q^{- 1},\, XPX^{-1}=Q, \, XQX^{-1}=PQ},
\end{equation}
and thus $\tonestar$ is a semi-direct product of $\ang{P,Q}\cong \quat$ by $\ang{X}\cong \Z_{3}$. Also, $\tonestar$ is abstractly a subgroup of $\oonestar$ and of $\istar$. We refer the reader to~\cite{AM,Cox,CM,GG12,Wo} for more properties of the binary polyhedral groups. One important property that they share with the family of cyclic and dicyclic groups is that they possess a unique element of order~$2$ (except for cyclic groups of odd order), which is a ramification of the fact that they are periodic in the sense of \resec{homotype}, and that in the non-cyclic case, this element generates the centre of the group. Further, the quotient by the unique subgroup of order $2$ induces a correspondence between the family of even-order cyclic, dicyclic and binary polyhedral groups with the finite subgroups of $\operatorname{SO}(3)$, the dicyclic group $\dic{4m}$ being associated with the dihedral group $\dih{2m}$ of order $2m$, and $\tonestar,\oonestar$ and $\istar$ being associated respectively with the polyhedral groups $\an[4]$, $\sn[4]$ and $\an[5]$. Using Kerckhoff's solution to the Nielsen realisation problem, Stukow classified the isomorphism classes of the finite subgroups of $\mcggen{\St}$, showing that they are finite subgroups of $\operatorname{SO}(3)$, with appropriate restrictions on $n$~\cite{St}. The analysis of \req{mcg} then leads to the complete classification of the isomorphism classes of the finite subgroups of $B_{n}(\St)$.

\pagebreak
\begin{thm}[\cite{GGagt2}]\label{th:finitebn}
Let $n\geq 3$. The isomorphism classes of the maximal finite subgroups of $B_n(\St)$ are as follows:
\begin{enumerate}[(a)]
\item\label{it:fina} $\Z_{2(n-1)}$ if $n\geq 5$.
\item\label{it:finb} $\dic{4n}$.
\item\label{it:finc} $\dic{4(n-2)}$ if $n=5$ or $n\geq 7$.
\item\label{it:find} $\tonestar$ if $n\equiv 4 \bmod 6$.
\item\label{it:fine} $\oonestar$ if $n\equiv 0,2\bmod 6$.
\item\label{it:finf} $\istar$ if $n\equiv 0,2,12,20\bmod 30$.
\end{enumerate}
\end{thm}

The geometric realisation of the finite subgroups of $B_{n}(\St)$ may be obtained by letting the corresponding finite subgroup of $\mcggen{\St}$ act by homeomorphisms on $\St$ (see~\cite[Section~3.2]{GGagt2} for more details). Concretely, consider the geometric definition given in \resec{defstring}. We visualise the space $\St \times [0,1]$ as that confined between two concentric spheres (see~\cite[page~41]{Ha} for example). For the (maximal) subgroups $\dic{4n}$, $\Z_{2(n-1)}$ and $\dic{4(n-2)}$, we attach strings, each representing the constant path in terms of the definition of \resec{particles}, to $n$ (resp.\ $n-1$, $n-2$) equally-spaced points on the equator, and $0$ (resp.\ $1$, $2$) points at the poles. For $\tonestar$, $\oonestar$ and $\istar$, the $n$ strings are attached symmetrically with respect to the associated regular polyhedron. We now let the corresponding finite subgroup of $\mcggen{\St}$ act on the inner sphere as a group of homeomorphisms, so that the set of basepoints is left invariant globally. This yields a subgroup of $B_{n}(\St)$, and one may check that it is exactly the given finite subgroup of \reth{finitebn}. In particular, a complete rotation of the inner sphere gives rise to the full twist braid $\ft$, and is a manifestation of the famous `Dirac string trick' (see~\cite[Section~6]{Fa},~\cite[page~43]{Ha} or~\cite[page~628]{Mag0}).

Algebraic representations of some of the binary polyhedral groups have been found: see~\cite[Remarks~3.2 and~3.3]{GGagt2} for realisations of $\tonestar$ in $B_{4}(\St)$ and $B_{6}(\St)$. Note however that in the second case there is a misprint, and the expression for $\delta$ should read
\begin{equation*}
\delta=\sigma_{3}^{-1} \sigma_{4}^{-1}\sigma_{5}^{-1} \sigma_{2}^{-1} \sigma_{1}^{-1}\sigma_{2}^{-1} \sigma_{5}\sigma_{4}\sigma_{5} \sigma_{5}\sigma_{4}\sigma_{3}.
\end{equation*}

By~\cite[Proposition~1.5]{GGagt2}, there are at most two conjugacy classes of each isomorphism class of the finite subgroups of $B_{n}(\St)$, and there is a single conjugacy class for each maximal finite subgroup.

As another application of \reco{embed}, we obtain the classification of  the finite subgroups of $B_{n}(\rp)$.
\begin{thm}[\cite{GG11}]\label{th:finitebnrp2}
Let $n\geq 2$. The isomorphism classes of the finite subgroups of $B_n(\rp)$ are the subgroups of the following groups:
\begin{enumerate}
\item\label{it:finrpa} $\dic{8n}$.
\item\label{it:finrpb} $\dic{8(n-1)}$ if $n\geq 4$.
\item\label{it:finrpc} $\oonestar$ if $n\equiv 0,1\pmod 3$.
\item\label{it:finrpd} $\istar$ if $n\equiv 0,1,6,10\pmod{15}$.
\end{enumerate}
\end{thm}
Although the groups involved in the statements of Theorems~\ref{th:finitebn} and~\ref{th:finitebnrp2} are basically the same, there is a difference in terms of those that are maximal. The finite groups described in \reth{finitebnrp2}(\ref{it:finrpa})--(\ref{it:finrpd}) are maximal in an abstract sense, while those of \reth{finitebn} are maximal with respect to inclusion. This is partly related to the fact that up to powers and conjugacy, $B_{n}(\St)$ has just three conjugacy classes of finite order elements, while $B_{n}(\rp)$ has many more. It could happen that a subgroup of $B_{n}(\rp)$ that is abstractly isomorphic to a proper subgroup of one of the groups given in \reth{finitebnrp2} be maximal with respect to inclusion. This is the subject of work in progress~\cite{GG14}. 

The proof of \reth{finitebnrp2} is obtained by combining \reco{embed} with \reth{finitebn}. In this way, we establish a list of possible finite subgroups of $B_n(\rp)$. Some of these possibilities are not realised (notably $\tonestar$ is not realised if $n\equiv 2 \pmod 3$, despite apparently being compatible with the embedding). The final step is to prove that the subgroups given in the statement of \reth{finitebnrp2} are indeed realised for the given values of $n$. This is achieved in a similar manner to that of the finite subgroups of $B_{n}(\St)$. As for $\St$, it is also possible to give explicit algebraic realisations of the dicyclic subgroups of $B_n(\rp)$. For example, we obtain $\ang{a,\garside}\cong \dic{8n}$ and $\ang{b,\garside a^{-1}}\cong \dic{8(n-1)}$ using \req{propgarside}~\cite[Proposition~15]{GG11}. Explicit realisations of $\tonestar$ and $\oonestar$ have been found in $B_{3}(\rp)$~\cite{GG14}, and applying \reco{embed} to them yields isomorphic copies in $B_{6}(\St)$.

As an application of \reth{finitebnrp2} and the short exact sequence~\reqref{mcg} for $\rp$, one may also obtain an alternative proof of the classification of the finite subgroups of $\operatorname{\mathcal{MCG}}(\rp,n)$ due to Bujalance, Cirre and Gamboa~\cite{BCG}. 

\begin{thm}[\cite{BCG}]\label{th:mcg}
Let $n\geq 2$. The finite subgroups of $\operatorname{\mathcal{MCG}}(\rp,n)$ are abstractly isomorphic to the subgroups of the following groups:
\begin{enumerate}
\item the dihedral group $\dih{4n}$ of order $4n$.
\item the dihedral group $\dih{4(n-1)}$ if $n\geq 3$.
\item $\sn[4]$ if $n\equiv 0,1\pmod 3$.
\item $\an[5]$ if $n\equiv 0,1,6,10\pmod{15}$.
\end{enumerate}
\end{thm}

One useful fact that is used to classify the virtually cyclic subgroups of $B_{n}(\St)$ is the knowledge of the centraliser and normaliser of its maximal finite cyclic and dicyclic subgroups. Note that if $i\in\brak{0,1}$, the centraliser of $\alpha_{i}$, considered as an element of $B_{n}$, is equal to $\ang{\alpha_{i}}$~\cite{BDM,GW}. A similar equality holds in $B_{n}(\St)$ and is obtained using \req{mcg} and the corresponding result for $\mcg$, which is due to Hodgkin~\cite{Ho}.
\begin{prop}[\cite{GG12}]\label{prop:genhodgkin1}
Let $i\in\brak{0,1,2}$, and let $n\geq 3$. 
\begin{enumerate}[(a)]
\item\label{it:centalphai} The centraliser of $\ang{\alpha_{i}}$ in $B_{n}(\St)$ is equal to $\ang{\alpha_{i}}$, unless $i=2$ and $n=3$, in which case it is equal to $B_{3}(\St)$.
\item\label{it:normcyclic} The normaliser of $\ang{\alpha_{i}}$ in $B_{n}(\St)$ is equal to:
\begin{equation*}
\begin{cases}
\ang{\alpha_{0},\garside} \cong \dic{4n} & \text{if $i=0$}\\
\ang{\alpha_{2},\alpha_{0}^{-1}\garside \alpha_{0}} \cong \dic{4(n-2)} & \text{if $i=2$}\\
\ang{\alpha_{1}}\cong \Z_{2(n-1)} & \text{if $i=1$,}
\end{cases}
\end{equation*}
unless $i=2$ and $n=3$, in which case it is equal to $B_{3}(\St)$.
\item If $i\in\brak{0,2}$, the normaliser of the standard copy of $\dic{4(n-i)}$ in $B_{n}(\St)$ is itself, except when $i=2$ and $n=4$, in which case the normaliser is equal to $\alpha_{0}^{-1} \sigma_{1}^{-1} \ang{\alpha_{0},\garside[4]} \sigma_{1}\alpha_{0}$, and is isomorphic to $\quat[16]$. 
\end{enumerate}
\end{prop}

A related problem is that of knowing which powers of $\alpha_{i}$ are conjugate in $B_{n}(\St)$, for each $i\in\brak{0,1,2}$. The answer is that such powers are either equal or inverse:
\begin{prop}[\cite{GG12}]\label{prop:genhodgkin2}
Let $n\geq 3$ and $i\in \brak{0,1,2}$, and suppose that there exist $r,m\in \Z$ such that $\alpha_{i}^m$ and $\alpha_{i}^r$ are conjugate in $B_{n}(\St)$. 
\begin{enumerate}
\item\label{it:conjpowera} If $i=1$ then $\alpha_1^m=\alpha_1^r$.
\item\label{it:conjpowerb} If $i\in \brak{0,2}$ then $\alpha_{i}^m=\alpha_{i}^{\pm r}$.
\end{enumerate}
\end{prop}
Once more, this generalises a corresponding result in $\mcggen{\St}$~\cite{Ho}. Using \reth{murasugi}, \repr{genhodgkin2} implies that if $F$ is a finite cyclic subgroup of $B_{n}(\St)$ then that the only possible actions of $\Z$ on $F$ are the trivial action and multiplication by $-1$. This also has consequences for the possible actions of $\Z$ on dicyclic subgroups of $B_{n}(\St)$. 

\subsection{Virtually cyclic subgroups of the braid groups of $\St$ and $\rp$}\label{sec:vcyclics}

In view of the Farrell-Jones Fibred Isomorphism Conjecture (see \resec{gens}), in order to calculate the lower algebraic $K$-theory of the group rings of the braid groups of $\St$ and $\rp$, it is necessary to know their virtually cyclic subgroups. Recall that a group is said to be \emph{virtually cyclic} if it contains a cyclic subgroup of finite index. It is clear from the definition that any finite subgroup is virtually cyclic, hence it suffices to concentrate on the \emph{infinite} virtually cyclic subgroups of these braid groups, which are in some sense their `simplest' infinite subgroups. The classification of the virtually cyclic subgroups of these braid groups is an interesting problem in its own right, and helps us to understand better the structure of these two groups. For the whole of this section, we refer the reader to~\cite{GG12} for more details.

Recall that by results of Epstein and Wall~\cite{Ep,W}, any infinite virtually cyclic group $G$ is isomorphic to $F\rtimes \Z$ or to $G_{1}\bigast_{F} G_{2}$, where $F$ is finite and $[G_{i}:F]=2$ for $i\in\brak{1,2}$. We shall say that $G$ is of \emph{Type~I} or \emph{Type~II} respectively. This enables us to establish a list of the possible infinite virtually cyclic subgroups of a given infinite group $\Gamma$, providing one knows its finite subgroups (which by Theorems~\ref{th:finitebn} and~\ref{th:finitebnrp2} is the case for our braid groups). The real difficulty lies in deciding whether the groups belonging to this list are indeed realised as subgroups of $\Gamma$. 

Let $n\geq 4$. In the case of $P_{n}(\St)$, as we saw in \resec{finitesubgp}, $\ang{\ft}$ is the only non-trivial finite subgroup, and since it is equal to the centre of $P_{n}(\St)$ by \repr{nocentres}(\ref{it:propsbns1}), it is then easy to see that the infinite virtually cyclic subgroups of $P_{n}(\St)$ are isomorphic to $\Z$ or to $\Z_{2}\times \Z$. The classification of the virtually cyclic subgroups of $P_{n}(\rp)$ was obtained in~\cite{GG8}, using \repr{finitepn}. Although the structure of the finite subgroups of $P_{n}(\rp)$ differs for $n=3$ and $n\geq 4$, up to isomorphism, the infinite virtually cyclic subgroups of $P_{n}(\rp)$ are the same for all $n\geq 3$:
\begin{thm}[\cite{GG8}]\label{th:vcp}
Let $n\geq 3$. The isomorphism classes of the infinite virtually cyclic subgroups of $P_{n}(\rp)$ are $\Z$,
$\Z_{2}\times \Z$ and $\Z_{4} \bigast_{\Z_{2}} \Z_{4}$.
\end{thm}

One obtains the classification of the virtually cyclic subgroups of $P_{n}(\rp)$ as a immediate corollary of \repr{finitepn} and \reth{vcp}~\cite{GG8}. One of the key results needed in the proof of \reth{vcp} is that $P_{n}(\rp)$ has no subgroup isomorphic to $\Z_{4} \times \Z$, which follows in a straightforward manner from \repr{order4}(\ref{it:order4b}). This fact allows us to eliminate several potential Type~I  and Type~II subgroups.

We now turn to the case of $B_{n}(\St)$. As we observed previously in \resec{sphere}, if $n\leq 3$ then $B_{n}(\St)$ is a known finite group, and so we shall suppose in what follows that $n\geq 4$. If $G$ is a group, let $\aut{G}$ (resp.\ $\out{G}$) denote the group of its automorphisms (resp.\ outer automorphisms). We define the following two families of virtually cyclic groups. 
\begin{defn}\label{def:v1v2}
Let $n\geq 4$.
\begin{enumerate}[(1)]
\item\label{it:mainIdef} Let $\mathbb{V}_{1}(n)$ be the family comprised of the following Type~I virtually cyclic groups:
\begin{enumerate}[(a)]
\item\label{it:mainzq} $\Z_{q}\times \Z$, where $q$ is a strict divisor of $2(n-i)$, $i\in \brak{0,1,2}$, and $q\neq n-i$ if $n-i$ is odd.
\item\label{it:mainzqt} $\Z_{q}\rtimes_{\rho} \Z$, where $q\geq 3$ is a strict divisor of $2(n-i)$, $i\in \brak{0,2}$, $q\neq n-i$ if $n$ is odd, and $\rho(1)\in \aut{\Z_{q}}$ is multiplication by $-1$.
\item\label{it:maindic} $\dic{4m}\times \Z$, where $m\geq 3$ is a strict divisor of $n-i$ and $i\in \brak{0,2}$.
\item\label{it:maindict} $\dic{4m}\rtimes_{\nu} \Z$, where $m\geq 3$ divides $n-i$, $i\in \brak{0,2}$, $(n-i)/m$ is even, and where $\nu(1)\in \aut{\dic{4m}}$ is defined by:
\begin{equation}\label{eq:actdic4m}
\left\{
\begin{aligned}
\nu(1)(x)&=x\\
\nu(1)(y)&=xy
\end{aligned}\right.
\end{equation}
for the presentation~\reqref{presdic} of $\dic{4m}$. 

\item\label{it:mainq8} $\quat\rtimes_{\theta} \Z$, for $n$ even and $\theta\in \operatorname{Hom}(\Z,\aut{\quat})$, for the following actions:
\begin{enumerate}[(i)]
\item $\theta(1)=\id$.
\item\label{it:mainIcii} $\theta=\alpha$, where $\alpha(1)\in \aut{\quat}$ is given by $\alpha(1)(i)=j$ and $\alpha(1)(j)=k$, where $\quat=\brak{\pm 1, \pm i, \pm j, \pm k}$.
\item\label{it:mainIciii} $\theta=\beta$, where $\beta(1)\in \aut{\quat}$ is given by $\beta(1)(i)=k$ and $\beta(1)(j)=j^{-1}$.
\end{enumerate}

\item\label{it:maint} $\tonestar \times \Z$ for $n$ even.
\item\label{it:maing} $\tonestar \rtimes_{\omega} \Z$ for $n\equiv 0,2 \bmod 6$, where $\omega(1)\in \aut{\tonestar}$ is the automorphism defined in terms of the presentation \reqref{preststar} by:
\begin{equation}\label{eq:nontrivacttstar}
\left\{
\begin{aligned}
P &\mapsto QP\\
Q &\mapsto Q^{-1}\\
X &\mapsto X^{-1}.
\end{aligned}\right.
\end{equation}

\item\label{it:maino} $\oonestar \times \Z$ for $n\equiv 0,2 \bmod 6$.
\item\label{it:maini} $\istar \times \Z$ for $n\equiv 0,2,12,20 \bmod{30}$.
\end{enumerate}
\item\label{it:mainIIdef} Let $\mathbb{V}_{2}(n)$ be the family comprised of the following Type~II virtually cyclic groups:
\begin{enumerate}[(a)]
\item\label{it:mainIIa} $\Z_{4q}\bigast_{\Z_{2q}} \Z_{4q}$, where $q$ divides $(n-i)/2$ for some $i\in\brak{0,1,2}$.

\item\label{it:mainIIb} $\Z_{4q}\bigast_{\Z_{2q}} \dic{4q}$, where $q\geq 2$ divides $(n-i)/2$ for some $i\in\brak{0,2}$.

\item\label{it:mainIIc} $\dic{4q}\bigast_{\Z_{2q}} \dic{4q}$, where $q\geq 2$ divides $n-i$ strictly for some $i\in\brak{0,2}$.

\item\label{it:mainIId} $\dic{4q}\bigast_{\dic{2q}} \dic{4q}$, where $q\geq 4$ is even and divides $n-i$ for some $i\in\brak{0,2}$. 

\item\label{it:mainIIe} $\oonestar \bigast_{\tonestar} \oonestar$, where $n\equiv 0,2 \bmod{6}$. 
\end{enumerate}
\end{enumerate}
Finally, let $\mathbb{V}(n)$ be the family comprised of the elements of $\mathbb{V}_{1}(n)$ and $\mathbb{V}_{2}(n)$. In what follows, $\rho, \nu,\alpha,\beta$ and $\omega$ will denote the actions defined in parts~(\ref{it:mainIdef})(\ref{it:mainzqt}), (\ref{it:mainIdef})(\ref{it:maindict}), (\ref{it:mainIdef})(\ref{it:mainq8})(\ref{it:mainIcii}), (\ref{it:mainIdef})(\ref{it:mainq8})(\ref{it:mainIciii}) and (\ref{it:mainIdef})(\ref{it:maing}) respectively.
\end{defn}

Up to a finite number of exceptions, we may then classify the infinite virtually cyclic subgroups of $B_{n}(\St)$.
\begin{thm}[\cite{GG12}]\label{th:main}
Suppose that $n\geq 4$. 
\begin{enumerate}[(1)]
\item\label{it:mainI} 
If $G$ is an infinite virtually cyclic subgroup of $B_{n}(\St)$ then $G$ is isomorphic to an element of $\mathbb{V}(n)$.

\item\label{it:mainII} Conversely, let $G$ be an element of $\mathbb{V}(n)$. Assume that the following conditions hold:
\begin{enumerate}[(a)]
\item\label{it:excepa} if $G\cong \quat \rtimes_{\alpha} \Z$ then $n\notin \brak{6,10,14}$.

\item if $G\cong \tonestar \times \Z$ then $n\notin \brak{4,6,8,10,14}$.

\item if $G\cong \oonestar \times \Z$ or $G\cong \tonestar \rtimes_{\omega} \Z$ then $n\notin \brak{6,8,12,14,18,20,26}$.

\item\label{it:excepd} if $G\cong \istar \times \Z$ then $n\notin \brak{12,20,30,32,42,50,62}$.

\item\label{it:excepe} if $G\cong \oonestar \bigast_{\tonestar} \oonestar$ then $n\notin \brak{6,8,12,14,18,20,24,26,30,32,38}$.
\end{enumerate}
Then there exists a subgroup of $B_{n}(\St)$ isomorphic to $G$.

\item\label{it:mainIII} Let $G$ be equal to $\tonestar\times \Z$ (resp.\ $\oonestar\times \Z$) if $n=4$ (resp.\ $n=6$). Then $B_{n}(\St)$ has no subgroup isomorphic to $G$.
\end{enumerate}
\end{thm}

\begin{rem}\label{rem:exceptions}
Together with \reth{finitebn}, \reth{main} yields a complete classification of the virtually cyclic subgroups of $B_{n}(\St)$ with the exception of a the thirty-eight cases for which the problem of their existence is open, given by the excluded values of $n$ in the above conditions~(\ref{it:mainII})(\ref{it:excepa})--(\ref{it:excepe}) but removing the two cases of part~(\ref{it:mainIII}) which we know not to be realised.
\end{rem}

The proof of \reth{main} is divided into two stages.  In conjunction with \reth{finitebn}, Epstein and Wall's results give rise to a family $\mathcal{VC}$ of virtually cyclic groups with the property that any infinite virtually cyclic subgroup of $B_{n}(\St)$ belongs to $\mathcal{VC}$. The first stage is to show that any such subgroup belongs in fact to the subfamily $\mathbb{V}(n)$ of $\mathcal{VC}$. This is achieved in several ways: the analysis of the centralisers and normalisers of the finite order elements of $B_{n}(\St)$ given in Propositions~\ref{prop:genhodgkin1} and~\ref{prop:genhodgkin2}; the study of the (outer) automorphism groups of the finite subgroups of \reth{finitebn}; and the periodicity of $B_{n}(\St)$ given by  \repr{per24}. Putting together these reductions allows us to prove \reth{main}(\ref{it:mainI}). The structure of the finite subgroups of $B_{n}(\St)$ imposes strong constraints on the possible Type~II subgroups, and the proof in this case is more straightforward than that for the Type~I subgroups. The second stage of the proof consists in proving the realisation of the elements of $\mathbb{V}(n)$ as subgroups of $B_{n}(\St)$ and to proving parts~(\ref{it:mainII}) and (\ref{it:mainIII}) of \reth{main}. The construction of the elements of $\mathbb{V}(n)$ involving finite cyclic and dicyclic groups as subgroups of $B_{n}(\St)$ is largely algebraic, and relies heavily on equations~\reqref{fundaa} and~\reqref{fundab} that describe the action by conjugation of the $\alpha_{i}$ on the generators of $B_{n}(\St)$. In contrast, the realisation of the elements of $\mathbb{V}(n)$ involving the binary polyhedral groups is geometric in nature, and occurs on the level of mapping class groups via the relation~\reqref{mcg} and the constructions of the finite subgroups of $B_{n}(\St)$ of \reth{finitebn}.

Since the open cases of \rerem{exceptions} only occur for even values of $n$, the complete classification of the infinite virtually cyclic subgroups of $B_{n}(\St)$ for all $n\geq 5$ odd follows directly from \reth{main}.
\begin{thm}[\cite{GG12}]\label{th:mainodd}
Let $n\geq 5$ be odd. Then up to isomorphism, the following groups are the infinite virtually cyclic subgroups of $B_{n}(\St)$.
\begin{enumerate}[(I)]
\item \begin{enumerate}[(a)]
\item $\Z_{m} \rtimes_{\theta} \Z$, where $\theta(1)\in\brak{\id, -\id}$, $m$ is a strict divisor of $2(n-i)$, for $i\in \brak{0,2}$, and $m\neq n-i$.
\item $\Z_{m} \times \Z$, where $m$ is a strict divisor of $2(n-1)$.
\item $\dic{4m} \times \Z$, where $m\geq 3$ is a strict divisor of $n-i$ for $i\in \brak{0,2}$.
\end{enumerate}
\item 
\begin{enumerate}[(a)]
\item $\Z_{4q}\bigast_{\Z_{2q}}\Z_{4q}$, where $q$ divides $(n-1)/2$.
\item $\dic{4q}\bigast_{\Z_{2q}}\dic{4q}$, where $q\geq 2$ is a strict divisor of $n-i$, and $i\in \brak{0,2}$
\end{enumerate}
\end{enumerate}
\end{thm}

Since in \reth{main} we are considering the realisation of the various subgroups up to isomorphism, one may ask whether each of the given elements of $\mathbb{V}_{2}(n)$ is unique up to isomorphism. It turns out that with with the exception of $\quat[16]\bigast_{\quat} \quat[16]$, abstractly there is only one way (up to isomorphism) to embed the amalgamating subgroup in each of the two factors, in other words for all of the other elements of $\mathbb{V}_{2}(n)$, the group is unique up to isomorphism~\cite{GG12}. Note that this result refers to abstract isomorphism classes of the given Type~II groups, and does not depend on the fact that the amalgamated products occurring as elements of $\mathbb{V}_{2}(n)$ are realised as subgroups of $B_{n}(\St)$. In the exceptional case of $\quat[16]\bigast_{\quat} \quat[16]$, abstractly there are two isomorphism classes defined respectively by:
\begin{equation*}
K_{1}=\setangr{x,y,a,b}{x^{4}=y^{2},\; a^{4}=b^{2},\; yxy^{-1}= x^{-1},\; bab^{-1}=a^{-1},\; x^{2}=a^{2},\; y=b}.
\end{equation*}
and
\begin{equation*}
K_{2}=\setangr{x,y,a,b}{x^{4}=y^{2},\; a^{4}=b^{2},\; yxy^{-1}= x^{-1},\; bab^{-1}=a^{-1},\; x^{2}=b,\; y=a^{2}b}.
\end{equation*}
If $n\geq 4$ is even, both $K_{1}$ and $K_{2}$ are realised as subgroups of $B_{n}(\St)$, with the possible exception of $K_{2}$ if $n\in\brak{6,14,18,26,30,38}$~\cite{GG12}. 

Using \req{mcg}, another consequence of \reth{main} is the classification of the virtually cyclic subgroups of $\mcg$, with a finite number of exceptions (see~\cite[Theorem~14]{GG12} for more details).

A similar analysis of the isomorphism classes of the infinite virtually cyclic subgroups of $B_{n}(\rp)$ is the subject of work in progress~\cite{GG14,GG16}.



\section{$K$-theory of surface braid groups}\label{sec:ktheory}

In this section, we indicate how the results of the previous sections may be used to compute the lower algebraic $K$-theory of the group rings of surface braid groups. In \resec{gens}, we start by recalling two conjectures of Farrell and Jones, whose validity for a given group provides a recipe to calculate its lower $K$-groups. In \resec{fjconjasp}, we outline the proof of the fact that surface braid groups of aspherical surfaces satisfy the Farrell-Jones conjecture, and in \resec{ktheorys2rp2}, we shall see how to extend this result to the braid groups of $\St$ and $\rp$. In order to calculate the lower algebraic $K$-theory of a group using this approach, one needs to be able to determine the lower $K$-groups of its virtually cyclic subgroups, as well as certain \emph{Nil groups} that are related to these subgroups. In \resec{genrems}, we recall some general methods that one may use to determine these lower $K$- and Nil groups. Finally, in \resec{ktheoryresults}, we state and outline the proofs of the known results, namely the lower $K$-groups of braid groups of aspherical surfaces, and of $P_{n}(\St)$, $P_{n}(\rp)$ and $B_{4}(\St)$.

\subsection{Generalities}\label{sec:gens}

Let $G$ be a discrete group and let $\Z[G]$ denote its integral group ring. The approach to the algebraic $K$-theoretical calculations of $\Z[G]$, which we outline in this section, consists in using the Farrell-Jones (Fibred) Isomorphism Conjecture that proposes to compute the $K$-groups of $\Z[G]$ from two sources: first, the algebraic $K$-theory of the class of virtually cyclic subgroups of $G$, and secondly, homological data.

\begin{defn}
A collection $\mathcal{F}$ of subgroups of $G$ is called a \emph{family} if:
\begin{enumerate}[(a)]
\item if $H\in\mathcal{F}$ and $A\leq H$ then $A\in\mathcal{F}$, and
\item if $H\in\mathcal{F}$ and $g\in G$ then $gHg^{-1}\in\mathcal{F}$.
\end{enumerate}
\end{defn}

The collection of finite subgroups of $G$, denoted $\mathit{\mathcal{F}in}$, and that of the virtually cyclic subgroups of $G$, denoted $\mathcal{VC}$, are examples of families of $G$. Given a family $\mathcal{F}$ of subgroups of $G$, a \emph{universal space} for $G$ with isotropy in $\mathcal{F}$ is a $G$-space $E\mathcal{F}$ that satisfies the following properties:
\begin{enumerate}
  \item the fixed set $E\mathcal{F}^H$ is non empty and contractible for all $H\in \mathcal{F}$, and
  \item the fixed set $E\mathcal{F}^H$ is empty for all $H\not\in \mathcal{F}$.
\end{enumerate} 
Universal spaces exist and are unique up to $G$-homotopy~\cite{tD}. If $\mathcal{F}$ consists of the trivial subgroup of $G$, the corresponding universal space is the universal space for principal $G$-bundles, and if $\mathcal{F}=\mathit{\mathcal{F}in}$, the corresponding universal space is the universal space for proper actions. If $\mathcal{F}=\mathcal{VC}$, we denote the corresponding universal space by $\underline{\underline{E}}G$. Although universal spaces exist for any family of subgroups of $G$, models for $E\mathcal{VC}$ that are suitable for making computations are still sparse, but there are some constructions for hyperbolic groups~\cite{JL} and $\operatorname{\emph{CAT}}(0)$ groups~\cite{F,L}.

Let $R$ be a ring with unit, and let $\operatorname{\emph{Or}}_{\mathcal{F}}$ is the orbit category of the group $G$ restricted to the family $\mathcal{F}$. J.~Davis and W.~L\"uck constructed a functor $\map{\mathbb{K}}{\operatorname{\emph{Or}}_{\mathcal{F}}(G)}[\mathit{Spectra}]$~\cite{DL}, whose value at the orbit $G/H$ is the non-connective algebraic $K$-theory spectrum of Pedersen-Weibel~\cite{PW}, and which satisfies the fundamental property that $\pi_i(\mathbb{K}(G/H))=K_i(\Z[H])$. The $K$-theoretical formulation of the Farrell-Jones isomorphism conjecture is as follows (one may consult~\cite{DL,JS} for more details).

\begin{fjconjecture}
Let $G$ be a discrete group. Then the assembly map
$$
\map{A_\mathcal{VC}}{H^G_n(\evc{G};\dbK)}[H^G_n(pt;\dbK)\cong K_n(\Z[G])],
$$
induced by the projection $\evc{G}\to \mathit{pt}$ is an isomorphism, where $H^G_n(-;\dbK)$ is a generalised equivariant homology theory with local coefficients in the functor $\dbK$, and $\evc{G}$ is a model for the universal space for the family $\mathcal{VC}$.
\end{fjconjecture}

A version of IC that is suitable for more general situations is the \emph{Fibred Farrell-Jones Conjecture} (FIC), which we now describe.  Given a group homomorphism $\map{\phi}{K}[G]$ and a family $\mathcal{F}$ of subgroups of a group $G$ that is also closed under finite intersections, the induced family on $K$ by $\phi$ is defined by: 
$$
\phi^*\mathcal{F}=\setr{H\leq K}{\phi(H)\in\mathcal{F}}.
$$

\begin{ficconj}[\cite{BLR}]
Let $G$ be a discrete group and let $\mathcal{F}$ be a family of subgroups of $G$. The pair $(G,\mathcal{F})$ is said to satisfy the \emph{Fibred Isomorphism Conjecture} if for all group homomorphisms $\map{\phi}{K}[G]$,  the assembly map
$$
\map{A_{\phi^*\mathcal{F}}}{H^K_n(E{\phi^*\mathcal{F}};\dbK)}[H^K_n(pt;\dbK)]
$$ 
is an isomorphism for all $n\in\Z$.
\end{ficconj}

Note that the validity of FIC implies that of IC by taking $K=G$ and $\phi=\id$. Two of the fundamental properties of FIC are as follows. 

\begin{thm}[\cite{BLR2}]\label{th:ficsubgroup}
If $G$ is a group that satisfies FIC and $H$ is a subgroup of $G$ then $H$ also satisfies FIC.
\end{thm}

\begin{thm}[\cite{BLR2}]\label{th:ficfiber}
Let $\map{f}{G}[Q]$ be a surjective group homomorphism. Assume that $(Q,\mathcal{VC}(Q))$ satisfies FIC and that IC is satisfied for all $H\in f^*\mathcal{VC}(Q)$. Then $(G,\mathcal{VC}(G))$ satisfies FIC.
\end{thm}

The Fibred Isomorphism Conjecture has been verified for word hyperbolic groups by A.~Bartels, W.~L\"uck and H.~Reich~\cite{BLR}, for $\operatorname{\emph{CAT}}(0)$ groups by C.~Wegner~\cite{We}, and for $\operatorname{SL}_n(\Z), n\geq 3$, by A.~Bartels, W.~L\"uck, H.~Reich and H.~Rueping~\cite{BLRR}. We record two of these results for future reference.

\begin{thm}[\cite{BLR}]\label{th:fichyp}
If $G$ is a hyperbolic group in the sense of Gromov then $G$ satisfies FIC.
\end{thm}

\begin{thm}[\cite{We}]\label{th:ficcat0}
If $G$ is a $\operatorname{\emph{CAT}}(0)$ group then $G$ satisfies FIC.
\end{thm}

The validity of the Fibred Isomorphism Conjecture has recently been shown for braid groups by D.~Juan-Pineda and L.~S\'anchez~\cite{JS} (see Theorems~\ref{th:Braidfic},~\ref{th:fics2} and~\ref{th:ficrp}). We will sketch the proofs in Sections~\ref{sec:fjconjasp} and~\ref{sec:ktheorys2rp2}. The original isomomorphism conjecture by T.~Farrell and L.~Jones was stated in~\cite{FJ0}. They proved several cases of the conjecture for the \emph{pseudoisotopy} functor. Here we shall only treat the case of the conjecture for the \emph{algebraic} $K$-theory functor.

\subsection{The $K$-theoretic Farrell-Jones Conjecture for braid groups of aspherical surfaces}\label{sec:fjconjasp}

In this section, we outline the ingredients needed to prove that braid groups of the plane or a compact surface other than the sphere or the projective plane satisfy FIC. The main tools that we shall require are the concepts of poly-free and strongly poly-free groups, which we now recall.

\begin{defn}\label{pf}
 A group $G$ is said to be \emph{poly-free} if there exists a filtration $1=G_0\subset G_1\subset\cdots\subset G_n=G$ of normal subgroups such that each quotient $G_{i+1}/G_i$ is a finitely-generated free group.
\end{defn}

The following result is due to D.~Juan-Pineda and L.~S\'anchez~\cite{JS}.

\begin{thm}[\cite{JS}]\label{th:ficpf}
 If $G$ is a poly-free group then $G$ satisfies FIC.
\end{thm}

The proof uses induction on the length of the filtration and the fact that the initial induction step is applied to a hyperbolic group.

Suppose first that $M$ is either the complex plane or a compact surface with non-empty boundary. Taking $r=1$ in \req{fnles} yields the following Fadell-Neuwirth fibration:
$$
F_{m+1,n-1}(\Int{M})\to F_{m,n}(\Int{M})\to F_{m,1}(\Int{M}),
$$
so by \reth{fnses}, we obtain the short exact sequence~\reqref{fnses}:
$$
1\to P_{n-1}(M\setminus \mathcal{Q}_{m+1})\to P_{n}(M\setminus \mathcal{Q}_m)\to \pi_{1}(M\setminus \mathcal{Q}_m)\to 1.
$$
It thus follows that for all $i\in {1,\ldots,n}$, $P_{i-1}(M\setminus \mathcal{Q}_{n-i+m+1})$ is normal in $P_{i}(M\setminus \mathcal{Q}_{n+m-i})$, and the corresponding quotient is isomorphic to the free group $\pi_{1}(M\setminus \mathcal{Q}_m)$ that is of finite rank. Setting $G_i=P_{i}(M\setminus \mathcal{Q}_{n-i})$ for all $i\in {0,1,\ldots,n}$ gives rise to a filtration that yields a poly-free structure for $P_{n}(M)$, and applying \reth{ficpf}, we obtain the following:

\begin{thm}[\cite{JS}]
Assume that $M=\mathbb{C}$ or that $M$ is a compact surface with non-empty boundary. Then the pure braid group $P_n(M)$ is poly-free, and thus satisfies FIC.
\end{thm}

Now suppose that $M$ is a compact aspherical surface with empty boundary. Taking $m=0$ and $r=1$ in \req{fnles} gives rise to the Fadell-Neuwirth fibration $F_{0,n-1}(M\setminus\mathcal{Q}_{1})\to F_{0,n}(M)\stackrel{p}{\to} F_{0,1}(M)=M$, and by \reth{fnses} induces the following short exact sequence:
$$
1\to P_{n-1}(M\setminus\mathcal{Q}_{1})\to P_n(M)\stackrel{p_{\#}}{\to} \pi_1(M)\to 1.
$$
Since $M$ is aspherical, the group $\pi_1(M)$ is finitely-generated Abelian or hyperbolic, and so satisfies FIC by Theorems~\ref{th:fichyp} and~\ref{th:ficcat0}. Now $\ker{p_{\#}}\cong P_{n-1}(M\setminus\mathcal{Q}_{1})$ is poly-free and $p_{\#}^{-1}(C)\cong P_{n-1}(M\setminus\mathcal{Q}_{1})\rtimes C$ where $C$ is any cyclic subgroup of $\pi_1(M)$, which is also poly-free, hence in both cases they satisfy FIC. \reth{ficfiber} then implies that $P_n(M)$ satisfies FIC. Putting together the two cases gives:
\begin{thm}[\cite{JS}]
Assume that $M=\mathbb{C}$ or that $M$ is a compact surface other than the sphere or the projective plane. Then the pure braid group $P_n(M)$ satisfies FIC for all $n\geq 1$.
\end{thm}

The next step is to go from $P_n(M)$ to $B_n(M)$. The idea is to embed the given group in a larger group (a wreath product in fact) that satisfies FIC and then apply \reth{ficsubgroup}. We start by adding one more property to  the definition of poly-free group.

\begin{defn}[\cite{AFR}]\label{spf}
A group $G$ is called \emph{strongly poly-free} (SPF) if it is poly-free and the following condition holds: for each $g\in G$ there exists a compact surface $M$ and a diffeomorphism $\map{f}{M}$ such that the action $C_{g}$ by conjugation of $g$ on $G_{i+1}/G_i$ may be realised geometrically, \emph{i.e.} the following diagram commutes:
$$
\xymatrix{\pi_1(M) \ar[r]^{f_\#}\ar[d]_{\varphi} & \pi_1(M) \\ G_{i+1}/G_i \ar[r]^{C_g} & G_{i+1}/G_i\ar[u]_{\varphi^{-1}}}
$$
where $\phi$ is a suitable isomorphism.
\end{defn}

The following result was proved in~\cite{AFR}.

\begin{thm}[\cite{AFR}]\label{th:Braidspf}
Assume that $M=\mathbb{C}$ or that $M$ is a compact surface with non-empty boundary. Then $P_n(M)$ is an SPF group for all $n\geq 1$.
\end{thm}

One of the main theorems in~\cite{JS} is the following:
\begin{thm}[\cite{JS}]\label{th:ficspfwrfinite}
Let $G$ be an SPF group, and let $H$ be a finite group. Then the wreath product $G\wr H$ satisfies FIC.
\end{thm}

We also recall the following result due to A.~Bartels, W.~L\"uck and H.~Reich~\cite{BLR2}.
\begin{lem}[\cite{BLR2}]\label{lem:kernelvc}
 Let $1\to K\to G\to Q\to 1$ be a short exact sequence of groups. Assume that $K$ is virtually cyclic and that $Q$ satisfies FIC. Then $G$ satisfies FIC.
\end{lem}

Moreover, given a finite extension of a group of the form 
$$
1\to G\to\Gamma\to H\to 1,
$$
where $H$ is a finite group, it follows that there is an injective homomorphism $\Gamma\lhra G\wr H$~\cite[Algebraic Lemma]{FR}. Since $P_{n}(M)$ is of finite index in $B_{n}(M)$ by \req{permseq}, it follows from Theorems~\ref{th:Braidspf} and~\ref{th:ficspfwrfinite} and the above observation that:

\begin{thm}[\cite{JS}]\label{th:Braidfic}
Assume that $M=\mathbb{C}$ or that $M$ is a compact surface other than the sphere or the projective plane. Then the full braid group $B_n(M)$ satisfies FIC for all $n\geq 1$.
\end{thm}

\subsection{The Farrell-Jones Conjecture for the braid groups of $\St$ and $\rp$}\label{sec:ktheorys2rp2}

The results of \resec{fjconjasp} treat the case of the braid groups of all surfaces with the exception of $\St$ and $\rp$. In this section, we outline the proof of the fact that the braid groups of these two surfaces also satisfy FIC.

Let $n\in \N$. Recall from \resec{sphere} that $P_n(\St)$ is trivial for $n=1,2$, and that $P_3(\St)\cong \Z_{2}$, hence these groups satisfy trivially FIC. So suppose that $n>3$. Taking $m=0$, $r=3$ and $M=\St$ in \req{fnles}, we obtain the following fibre bundle:
$$
F_{2,n-3}(\mathbb{C})\approx F_{3,n-3}(\St)\to F_{0,n}(\St)\to F_{0,3}(\St),
$$
and by \reth{fnses}, its long exact sequence in homotopy yields the Fadell-Neuwirth short exact sequence:
$$
1\to P_{n-3}(\mathbb{C}\setminus \mathcal{Q}_2) \to P_{n}(\St)\to P_{3}(\St)\to 1.
$$
Observe that $G=P_{n-3}(\mathbb{C}\setminus \mathcal{Q}_2)$ is an SPF group as it is part of the filtration of $P_{n-3}(\mathbb{C})$, hence Theorems~\ref{th:ficsubgroup} and~\ref{th:ficspfwrfinite} imply that $\pi_1(F_{0,n}(\St))=P_n(\St)$ satisfies FIC. In~\cite{Mi}, S.~Mill\'an-Vossler proved that $B_n(\St)$ fits in an extension of the form:
$$
1\to G\to B_n(\St)/\langle\ft\rangle\to \sn\to 1
$$
(this is a consequence of \req{permseq} and Propositions~\ref{prop:nocentres} and~\ref{prop:splitpns}), so $B_n(\St)/\langle\ft\rangle$ satisfies FIC by Theorems~\ref{th:ficsubgroup} and~\ref{th:ficspfwrfinite}. Taking $M=\St$ in \req{mcg} and applying \relem{kernelvc}, we see that $B_n(\St)$ satisfies FIC. Summing up these considerations, we obtain:
\begin{thm}[\cite{JS}]\label{th:fics2}
Both $P_n(\St)$ and $B_n(\St)$ satisfy FIC for all $n\geq 1$.
\end{thm}

The situation for $\rp$ is similar. Consider first the case of $P_{n}(\rp)$. By \resec{sphere}, $P_1(\rp)\cong \Z_{2}$, $P_2(\rp)\cong \quat$ and $P_3(\rp)\cong \F[2]\rtimes\quat$. It follows that $P_1(\rp)$ and $P_2(\rp)$ satisfy FIC as they are finite, and that $P_3(\rp)$ also satisfies FIC by \reth{fichyp} since it is (virtually) hyperbolic. Now let $n>3$. Taking the short exact sequence~\reqref{sesgennr} with $M=\rp$ and $r=2$ gives rise to the following short exact sequence:
$$
1\to G\to P_n(\rp)\to \quat\to 1,
$$
where $G=P_{n-2}(\rp\setminus \mathcal{Q}_{2})$ is an SPF group. It follows once more from Theorems~\ref{th:ficspfwrfinite} and~\ref{th:ficsubgroup} that $P_n(\rp)$ satisfies FIC for all $n>3$. Passing to the case of $B_{n}(\rp)$, note that $B_1(\rp)\cong \Z_{2}$ and $B_2(\rp)\cong\quat[16]$ by \resec{sphere}. Now $G$ is not normal in $B_{n}(\rp)$, but the intersection $H$ of its conjugates in $B_{n}(\rp)$ is a finite-index normal subgroup of both $G$ and $B_{n}(\rp)$, and for all $n\geq 3$, $B_n(\rp)$ fits in a short exact sequence: 
$$
1\to H\to B_n(\rp)\to B_n(\rp)/H\to 1,
$$
where $B_n(\rp)/H$ is finite. Since $G$ is SPF, it follows from~\cite{Mi} that $H$ is also SPF, and we conclude from \reth{ficspfwrfinite} and~\cite[Algebraic Lemma]{FR} that $B_n(\rp)$ satisfies FIC. We record these results as follows.
\begin{thm}[\cite{JS}]\label{th:ficrp}
Both $P_n(\rp)$ and $B_n(\rp)$ satisfy FIC for all $n\geq 1$.
\end{thm}

\subsection{General remarks for computations}\label{sec:genrems}

As we mentioned before, the validity of FIC should, in principle, furnish the necessary tools needed to compute the algebraic $K$-groups of the group rings for surface braid groups. We will concentrate in this section on \emph{lower} $K$-groups, that is $K_i(-)$ for $i\leq 1$.  Recall that the domain of the assembly map in the statement of IC is
\begin{align}\label{eq:eqhomology} 
H^G_n(\evc{G};\dbK).
\end{align}
This is an extraordinary equivariant homology theory whose coefficients are the functor $\dbK$. The input of $\dbK$ consists of the orbits of the type $G/V$, where $V$ varies over the virtually cyclic subgroups of $G$, and its values at these orbits are the spectra $\dbK(G/V)$ whose homotopy groups are given by $\pi_i(\dbK(G/V))\cong K_i(\Z[V])$. On the other hand, there is an Atiyah-Hirzebruch-type spectral sequence that computes the equivariant homology groups of \req{eqhomology} whose $E_2$-term is given by:
$$
 E^{p,q}_2\cong H_p(\bvc{G};\brak{K_q}),
 $$
where this is now an ordinary homology theory whose local coefficients are the algebraic $K$-groups of the virtually cyclic subgroups of $G$, and which appear as isotropy at different subcomplexes of $\bvc{G}=\evc{G}/G$. In summary, in order to compute $H^G_n(\evc{G};\dbK)$, we need to understand the following:
\begin{enumerate}[(a)]
\item the algebraic $K$-groups $K_i(\Z[V])$ for all $i\leq n$ and all virtually cyclic subgroups $V$ of $G$.
\item the spaces $\evc{G}$ and $\bvc{G}$.
\item how these groups and spaces are assembled together. This is encoded in the spectral sequence.
\end{enumerate}
 
Let $V$ be a virtually cyclic group. As indicated in \resec{vcyclics}, $V$ is either finite, of Type~I (so is isomorphic to a semidirect product of the form $F\rtimes \Z$, where $F$ is finite), or of Type~II (so is isomorphic to an amalgam of the form $G_{1}\bigast_{F} G_{2}$ where $F$ is of index $2$ in both $G_{1}$ and $G_{2}$.). In the Type~II case, $V$ fits in a short exact sequence of the form: 
$$
1\to F\to V\to \dih{\infty}\to 1,
$$
where $F$ is a finite group and $\dih{\infty}$ is the infinite dihedral group. The computation of the algebraic $K$-theory groups for each of these cases is currently an active area of study. In general, finite groups may be treated with induction-restriction methods, see~\cite{O}. We shall comment on the case of the finite subgroups of $B_{n}(\St)$ later on. In order to study the algebraic $K$-groups of Type~I and Type~II groups, we need some background. 
 
Let $R$ be an associative ring with unit, and let $R[t]$ denote its polynomial ring. Let $\map{\epsilon}{R[t]}[R]$ be the augmentation map induced by $t\mapsto 1$, and let $\map{\epsilon_{\ast}}{K_i(R[t])}[K_i(R)]$ be the homomorphism induced on $K$-groups.
 
 \begin{defn}
Let $R$ be an associative ring with unit. The Bass \emph{Nil} groups of $R$ are defined by:
$$
NK_i(R)=\ker{\epsilon_{\ast}}.
$$
\end{defn}
 
The Bass Nil groups appear in the study of $K$-groups of virtually cyclic groups via the Bass-Heller-Swan fundamental theorem:
\begin{thm}[Bass, Heller and Swan~\cite{B}]\label{th:BassHS}
Let $R$ be an associative ring with unit, and let $R[t,t^{-1}]$ be its Laurent polynomial ring. Then for all $i\in \Z$,
\begin{equation}\label{eq:bhs}
K_i(R[t,t^{-1}])\cong K_i(R)\oplus K_{i-1}(R)\oplus NK_i(R)\oplus  NK_i(R).
\end{equation}
\end{thm}
  
Observe that if a group $G$ is of the form $F\times\Z$ for some group $F$, its group ring may be described as follows:
$$
  \Z[G]=\Z[F\times\Z]\cong \Z[F][t,t^{-1}].
  $$
From \reth{BassHS}, we thus obtain:
\begin{cor}\label{cor:bhs}
The algebraic $K$-groups of a group $V=F\times \Z$ are of the form:
\begin{equation}\label{eq:bhsgen}
K_i(\Z[V])\cong K_i(\Z[F])\oplus K_{i-1}(\Z[F])\oplus NK_i(\Z[F])\oplus  NK_i(\Z[F]).
\end{equation}
  \end{cor}
  
If $V$ is as above and virtually cyclic, so $F$ is finite, \req{bhsgen} tells us that we need to compute the $K$-groups of the group ring $\Z[F]$ as well as the Bass Nil groups. If on the other hand, $V$ is a non-trivial semi-direct product of the form $V=F\rtimes_\alpha \Z$, where $\alpha$ denotes the action of $\Z$ on $F$, the corresponding group ring is the \emph{twisted} Laurent polynomial ring $\Z[F]_\alpha[t,t^{-1}]$. This case has been studied by T.~Farrell and W.~C.~Hsiang in~\cite{FHs}. They found a formula similar to that of \req{bhs} of Bass-Heller-Swan, but the terms $NK_i(\Z[F])\oplus NK_i(\Z[F])$ should be replaced by:
  $$
  NK_i(\Z[F],\alpha)\oplus NK_i(\Z[F],\alpha^{-1}),
  $$
which are similar groups that take into account the action of $\Z$ on $F$. These are now known as Farrell-Hsiang \emph{twisted Nil} groups. Together with the Bass Nil groups, these  Nil groups are the subject of investigation, full computations are few and far between, and they are in general very large groups due to the following fact:
\begin{thm}[\cite{Farr,Ra}]
Let $R$ be a ring. Then both the Bass Nil and Farrell-Hsiang Nil groups are either trivial or are not finitely generated.
\end{thm}

 The case of virtually cyclic groups of the form $V=A \bigast_F B$ is handled by the foundational work of F.~Waldhausen~\cite{Wa}. There is a long exact sequence of the form:
 \begin{multline*}
\cdots\to K_{n}(\Z [F])\to K_{n}(\Z [A])\oplus K_{n}(\Z [B])\to K_{n}(\Z [V])/\operatorname{\textbf{Nil}}_n^W\to K_{n-1}(\Z [F])\\
\to K_{n-1}(\Z [A])\oplus K_{n-1}(\Z [B])\to K_{n-1}(\Z [V])/\operatorname{\textbf{Nil}}_{n-1}^W\to \cdots,
\end{multline*}
where the term $\operatorname{\textbf{Nil}}_{n}^W$ denotes the \emph{Waldhausen Nil groups} defined by:
 \begin{equation*}
\operatorname{\textbf{Nil}}_n^W=\operatorname{Nil}^W_n(\Z[F];\Z[A\setminus F],\Z[B\setminus F]).
\end{equation*}
A somewhat better description of the Waldhausen Nil groups $\operatorname{\textbf{Nil}}_{n}^W$ is given in the work of J.~Davis, K.~Khan and A.~Ranicki~\cite{DKR} who identify these groups with  Farrell-Hsiang Nil groups of a group of the form $F\rtimes \Z$ for a suitable subgroup isomorphic to $\Z$ of the infinite dihedral group $\dih{\infty}=V/F$.
 
Some general results for algebraic $K$-groups for group rings of finite groups are known. We record some of them in the following proposition.

 \begin{prop}
Let $F$ be a finite group. Then:
\begin{enumerate}[(a)]
\item The groups $K_i(\Z[F])$ are finitely-generated Abelian groups for all $i\geq -1$.
\item The groups $K_i(\Z[F])$ vanish for $i<-1$.
\item The groups $NK_i(\Z[F])$ vanish for $i<0$.
\end{enumerate}
\end{prop}

The first part is proved in~\cite{Ku} if $i\geq 0$ and in~\cite{C} if $i=-1$, the second part is proved in~\cite{C}, and the third part in~\cite{C,FJ}.
 
On the other hand, the $NK_i(\Z[F])$ are non trivial for $i=0,1$ even for simple finite virtually cyclic groups, such as $F=\Z_2\times \Z_2$ or $\Z_{4}$~\cite{Weib}. It is therefore a challenge to decide whether the algebraic $K$-groups of infinite virtually cyclic groups are finitely-generated groups. The only known case that is always finitely generated is in degree $-1$: 
 \begin{prop}[\cite{FJ}]\label{negkvc}
   Let $V$ be a virtually cyclic group. Then: 
\begin{enumerate}
\item $ K_{-1}(\Z[V])$ is a finitely-generated group that is generated by the images of the homomorphisms $K_{-1}(\Z[G])\to K_{-1}(\Z[V])$ induced by the inclusions $G\lhra V$, where $G$ runs over the conjugacy classes of the finite subgroups of $V$.
\item the groups $K_i(\Z[V])$ are trivial for $i<-1$.
\end{enumerate}
\end{prop}
 
 We finish this section by recalling the lower $K$-groups of the integers $\Z$, which is fundamental for many of the calculations that follow.
 \begin{prop}\label{prop:kintegers}
For the ring $\Z$, the following results hold:
   \begin{enumerate}[(a)]
      \item\label{it:kgroupa} $K_i(\Z)$ is a finitely-generated Abelian group for all $i\in\Z$.
      \item\label{it:kgroupb} $K_1(\Z)=\Z_2$ and $K_0(\Z)=\Z$.
      \item\label{it:kgroupc} $K_i(\Z)=0$ for all $i<0$.
      \item\label{it:kgroupd} $NK_i(\Z)=0$ for all $i\in\Z$.
      \item\label{it:kgroupe} $K_i(\Z[\Z])\cong K_i(\Z)$ for all $i\in\Z$.
   \end{enumerate}
 \end{prop}
 
The proof of~(\ref{it:kgroupa}) may be found in~\cite{Q0}, and that of~(\ref{it:kgroupd}) follows from the regularity of $\Z$ and the work of D.~Quillen who showed that the Nil groups of regular rings vanish~\cite{Q1}. Part~(\ref{it:kgroupb}) is a consequence of the fact that $K_1(\Z)$ is just the units of $\Z$, and that every finitely-generated projective module over $\Z$ is free, and part~(\ref{it:kgroupc}) follows from the equality $\dim(\Z)=0$. Finally, part~(\ref{it:kgroupe}) is implied by the previous results and the Bass-Heller-Swan theorem (\reth{BassHS}).
 
We are interested in the non-trivial lower $K$-groups. Given a group $G$, we define $\widetilde{K}_i(\Z[G])$ to be the Whitehead group $\operatorname{Wh}(G)$ if $i=1$, the reduced $K_0$-group $\widetilde{K}_0(\Z[G])$ if $i=0$, and the usual $K_i$-groups if $i<0$. The results stated are valid for these reduced groups and for $i\leq 1$, and some of the computational results will be given for these reduced groups.  In this context, we may reinterpret \repr{kintegers} by saying that $\widetilde{K}_i(\Z)=0$ and $\widetilde{K}_i(\Z[\Z])=0$ for all $i\leq 1$.
 

\subsection{Computational results}\label{sec:ktheoryresults}

We now gather together the information obtained in the preceding sections. We start with the case of torsion-free braid groups, which by \reco{torsion} are precisely the braid groups of the complex plane or compact surfaces other than $\St$ or $\rp$. In this case, the only virtually cyclic subgroups of $G$ are trivial or infinite cyclic. By \repr{kintegers}, the reduced lower $K$-groups of $\Z$ and of $\Z[\Z]$ vanish, and the coefficients of the spectral sequence needed to compute the equivariant homology groups of \req{eqhomology}, whose coefficients are the reduced $K$-groups, are all trivial, so this spectral sequence collapses, thus yielding the trivial group. Hence:

\begin{thm}[\cite{AFR,JS}]\label{th:kthasp}
Let $G$ be the braid group (pure or full) of the complex plane or of a compact surface without boundary different from $\St$ and $\rp$. Then $\widetilde{K}_i(\Z[G])=0 \text{ for all } i\leq 1$.
\end{thm}
 
We now turn to the case of the pure braid groups of $\St$ and $\rp$. From the discussion just before the statement of \reth{vcp}, if $n\geq 4$, the infinite virtually cyclic subgroups $V$ of $P_n(\St)$ are isomorphic to $\Z$ or $\Z\times \Z_{2}$ and it is well known that $\widetilde{K}_i(\Z[V])=0$ for these two groups, using~\repr{kintegers} and \reco{bhs} for example. Since $P_1(\St)$ and $P_2(\St)$ are trivial and $P_3(\St)=\Z_2$ and the reduced lower $K$-groups of these groups also vanish, we have the following:
\begin{thm}[\cite{JM}]
For all  $i\leq 1$ and $n\geq 1$, $\widetilde{K}_i(\Z[P_n(\St)])=0$.
 \end{thm} 
 
The case of $P_{n}(\rp)$ is somewhat more involved. The reason is that by \repr{finitepn}, $\quat$ is realised as a subgroup of $P_{n}(\rp)$ if $n\in \brak{2,3}$, and its reduced $K$-group is non trivial in degree $0$.  More precisely, if $i\leq 1$,
\begin{equation*}
\widetilde{K}_i(\Z[\quat])=
\begin{cases} 
\Z_2 & \text{if $i=0$}\\
0 & \text{otherwise.}
 \end{cases}   
\end{equation*}
Since $P_{1}(\rp)\cong \Z_{2}$ and $P_{2}(\rp)\cong \quat$, we thus obtain the lower $K$-groups of these two groups. So assume that $n\geq 3$. With the exception of $\quat$, the reduced lower $K$-groups of the other finite subgroups of $P_{n}(\rp)$, as well as those of the infinite virtually cyclic subgroups given by \reth{vcp}, are trivial. From this, one may show that the reduced lower algebraic $K$-groups of $P_{n}(\rp)$ are as follows.
\begin{thm}[\cite{JM2}]
Suppose that $n\geq 3$ and $i\leq 1$. Then:
\begin{equation*}
\widetilde{K}_i(\Z[P_n(\rp)])=
\begin{cases}
\Z_2& \text{if $n=3$ and $i=0$}.\\
0 & \text{otherwise.}
\end{cases}
\end{equation*}
\end{thm}
 
The situation for the braid groups of both $\St$ and $\rp$ is currently the subject of investigation. By \reth{main}, the virtually cyclic subgroups of $B_n(\St)$ are known for all $n>3$, with the exception of a small number of cases. Many of the reduced lower $K$-groups of the finite subgroups of $B_n(\St)$ have been carried out. The $\widetilde{K}_{0}$-groups of the binary polyhedral groups and of the dicyclic groups $\dic{4m}$, $m\leq 13$, were computed in~\cite{Sw}. The Whitehead group of all finite subgroups of $B_{n}(\St)$ and the $K_{-1}$-groups of the binary polyhedral groups and of many dicyclic groups were determined in~\cite{GJM}. We remark that these $K_{-1}$-groups exhibit new structural phenomena that had not appeared previously in the study of the lower algebraic $K$-theory of other groups, such as the existence of torsion. These calculations are somewhat involved and require techniques from different areas. 

Passing to the case of the computation of the lower algebraic $K$-theory of $B_{n}(\St)$, $n\geq 4$, the only complete result so far is that for $n=4$~\cite{GJM}. We outline the steps in this case. A first important observation is that $B_4(\St)$ is isomorphic to an amalgamated product of the form $\quat[16]\bigast_{\quat} \tonestar$~\cite{GJM}. By \reth{finitebn} and~\cite[Proposition~1.5]{GGagt2}, the maximal finite subgroups of $B_{4}(\St)$ are isomorphic to $\tstar$ or $\quat[16]$, and there is a single conjugacy class of each. Moreover, we obtain the infinite virtually cyclic subgroups of $B_4(\St)$ from \reth{main}, and from this, one may deduce the maximal virtually cyclic subgroups of $B_{4}(\St)$:
\begin{thm}[\cite{GJM}]\label{th:maxvcb4}\mbox{}
\begin{enumerate}[(a)]
\item Every infinite maximal virtually cyclic subgroup of $B_4(\St)$ is isomorphic to $\quat[16] \bigast_{\quat} \quat[16]$ or to $\quat \rtimes \Z$ for one of the three possible actions (see part~(\ref{it:mainq8}) of the definition of the family $\mathbb{V}_{1}(n)$ in \resec{vcyclics}).
\item If $V$ is a finite maximal cyclic subgroup of $B_{4}(\St)$ then $V\cong \tstar$.
\item Let $G$ be a group that is isomorphic to $\quat \rtimes \Z$ for one of the three possible actions, or to $\quat[16] \bigast_{\quat} \quat[16]$. Then $B_{4}(\St)$ possesses both maximal and non-maximal virtually cyclic subgroups that are abstractly isomorphic to $G$.
\end{enumerate}
 \end{thm}
 
Calculations of the reduced lower algebraic $K$-groups of the groups given in \reth{maxvcb4} may be found in~\cite{GJM}. The next step is to find a model for $\evc{B_4(\St)}$. Since $B_{4}(\St)$ is an amalgam of finite groups, it follows that it is Gromov hyperbolic. If $G$ is a hyperbolic group, D.~Juan-Pineda and I.~Leary found a model for $\evc{G}$~\cite{JL}. In our case, this can be described as:
$$
\evc{B_4(\St)}= \mathbf{T}\bigast D,
$$
which is the join of a suitable tree $\mathbf{T}$ and a countable discrete set $D$. From this description, it also follows that the equivariant homology groups of \req{eqhomology} are isomorphic to:
 $$
 H^{B_4(\St)}_n(\mathbf{T};\{\dbK\})\oplus \left(\bigoplus_{V\in \operatorname{Max}(\mathcal{VC}(B_4(\St)))}NIL_n(V)\right),
 $$
 where $NIL_n$ denotes one of the Nil groups described above according to the type of infinite virtually cyclic group involved, and $\operatorname{Max}(\mathcal{VC}(B_4(\St)))$ is a set of representatives of the conjugacy classes of maximal infinite virtually cyclic subgroups of $B_4(\St)$. We summarise the final result for $B_4(\St)$ as follows.
\begin{thm}[\cite{GJM}]\label{th:b4s2kth}
The reduced lower algebraic $K$-groups for $B_4(\St)$ are given by
$$
\widetilde{K}_i(\Z[B_4(\St)])=
\begin{cases}
\Z\oplus \operatorname{Nil}_1, & \text{if $i=1$}\\
\Z_2\oplus \operatorname{Nil}_0,& \text{if $i=0$}\\
\Z_2\oplus \Z & \text{if $i=-1$}\\
0& \text{if $i<-1$,}
\end{cases}
$$
where for $i=0,1$, 
$$
\operatorname{Nil}_i\cong \bigoplus_{\infty}[ 2(\Z_2)^{\infty} \oplus W],
$$
$2(\Z_2)^{\infty}$ denotes two infinite countable direct sums of copies of $\Z_2$, and $W$ is an infinitely-generated Abelian group of exponent $2$ or $4$. 
 \end{thm}
 
Since the groups $\quat \rtimes \Z$ and $\quat[16]\bigast_{\quat}\quat[16]$ that appear in the statement of \reth{maxvcb4} appear as maximal subgroups of $B_4(\St)$, they contribute in a non-trivial manner via the Bass, Farrell-Hsiang and Waldhausen Nil groups to the reduced lower $K$-groups of  $\Z[B_4(\St)]$.

\subsection{Remarks}
\begin{enumerate}[(a)]
\item We have concentrated on the \emph{lower} algebraic $K_i$-groups, that is, in degrees $i\leq 1$. This is due to our lack of knowledge about $K_i(\Z[V])$ if $V$ is a virtually cyclic group if $i>1$. Little is known about the $K_{i}$-groups for $i>1$, even for finite groups. One example for $i=2$ may be found in~\cite{JLMP}.
 
\item In~\cite{GJM}, J.~Guaschi, D.~Juan-Pineda and S. Mill\'an-L\'opez developed techniques to compute reduced lower algebraic $K$-groups of many of the finite subgroups of $B_{n}(\St)$, in particular for small values of $n$. Some other results concerning these computations will appear in~\cite{GJM2}. How these subgroups are assembled to build up all of the reduced lower $K$-groups of a specific braid group $B_n(\St)$ for $n>4$ is the subject of work in progress. The main missing ingredient is the construction of a suitable model for $\evc{B_n(\St)}$. Note that the amalgamated product structure of $B_{4}(\St)$ is specific to this case, and we cannot hope for it to be carried over to braid groups with more strings.
     
\item The case of $B_{n}(\rp)$ is also still open if $n\geq 3$. However, many features are currently being studied: the classification of the virtually cyclic subgroups of $B_{n}(\rp)$~\cite{GG14,GG16}, as well as their $K$-groups and models for the corresponding universal spaces. 
     
\item In work in progress, it has been proved by D.~Juan-Pineda and L.~S\'anchez that if $G$ is a hyperbolic group, then $\operatorname{rank}(K_i(\Z[G]))< \infty$ for all $i\in \Z$. From this we have that $\operatorname{rank}(K_i(B_4(\Z[\St])))<\infty$ for all $i\in\Z$.
 \end{enumerate}
 




\begin{thebibliography}{999}
\renewcommand{\bibname}{References} 
\addcontentsline{toc}{section}{\bibname} 
{\small


\bibitem{AM} A.~Adem and R.~J.~Milgram, Cohomology of finite groups, Springer-Verlag, New York-Heidelberg-Berlin (1994).

\bibitem{AS} A.~Adem and J.~H.~Smith, Periodic complexes and group actions, \emph{Ann.\ Math.} \textbf{154}  (2001), 407--435.

\bibitem{ATD} Algebraic topology discussion list, January 2004, \url{http://www.lehigh.edu/~dmd1/pz119.txt}.


\bibitem{AFR} C.~S.~Aravinda, F.~T.~Farrell and S.~K.~Roushon, Algebraic $K$-theory of pure braid groups, \emph{Asian J.~Math.} \textbf{4} (2000), 337--343.

\bibitem{A1} E.~Artin, Theorie der Z\"opfe, \emph{Abh.\ Math.\ Sem.\ Univ.\ Hamburg} \textbf{4} (1925), 47--72.

\bibitem{A2} E.~Artin, Theory of braids, \emph{Ann.\ Math.} \textbf{48} (1947), 101--126.

\bibitem{A3} E.~Artin, Braids and permutations, \emph{Ann.\ Math.} \textbf{48} (1947), 643--649.

\bibitem{Bar} V.~Bardakov, Linear representations of the group of conjugating automorphisms and the braid groups of some manifolds, \emph{Siberian Math.~J.} \textbf{46} (2005), 13--23.     

\bibitem{Bar2} V.~Bardakov, Linear representations of the braid groups of some manifolds, \emph{Acta Appl.\ Math.} \textbf{85} (2005), 41--48.

\bibitem{bmvw} V.~Bardakov, R.~Mikhailov, V.~V.~Vershinin and J.~Wu, Brunnian braids on surfaces, \emph{Algebr.\ Geom.\ Topol.} \textbf{12} (2012) 1607--1648.

\bibitem{BLR} A.~Bartels, W.~L\"uck and H.~Reich, The $K$-theoretic Farrell-Jones conjecture for hyperbolic groups, \emph{Invent.\ Math.} \textbf{172} (2008), 29--70.

\bibitem{BLR2} A.~Bartels, W.~L\"uck and H.~Reich, On the Farrell-Jones conjecture and its applications, \emph{J.~Topol.} \textbf{1} (2008), 57--86.

\bibitem{BLRR} A.~Bartels, W.~ L\"uck, H.~Reich and  H.~Rueping, $K$- and $L$-theory of group rings over $GL_n(\Z)$, preprint, \url{arXiv:1204.2418}.

\bibitem{B} H.~Bass, Algebraic $K$-theory, W.~A.~Benjamin Inc., New York-Amsterdam, 1968.

\bibitem{Baue} H.~J.~Baues, Obstruction theory on homotopy classification of maps, Lecture Notes in Mathematics \textbf{628}, Springer-Verlag, Berlin, 1977.

\bibitem{Be} P.~Bellingeri, On presentation of surface braid groups, \emph{J.~Algebra} \textbf{274} (2004), 543-563.

\bibitem{BGG} P.~Bellingeri, S.~Gervais and J.~Guaschi, Lower central series of Artin-Tits and surface braid groups, \emph{J.~Algebra} \textbf{319} (2008), 1409--1427.

\bibitem{BG} P.~Bellingeri and E.~Godelle, Positive presentations of surface braid groups, \emph{J.~Knot Theory Ramif.} \textbf{16} (2007), 1219--1233.

\bibitem{BGG2} P.~Bellingeri, E.~Godelle and J.~Guaschi, Exact sequences, lower central series and representations of surface braid groups, preprint, \url{arXiv math:1106.4982}.

\bibitem{BCHWW} A.~J.~Berrick, F.~R.~Cohen, E.~Hanbury, Y.-L.~Wong and J.~Wu, Braids: Introductory Lectures on Braids, Configurations and Their Applications, Lecture Notes Series, Institute for Mathematical Sciences, National University of Singapore, Vol.~19, World Scientific, 2010.

\bibitem{BCWW}  A.~J.~Berrick, F.~R.~Cohen, Y.-L.~Wong and J.~Wu, Configurations, braids, and homotopy groups, \emph{J.~Amer.\ Math.\ Soc.} \textbf{19}  (2006), 265--326.

\bibitem{BDM} D.~Bessis, F.~Digne and J.~Michel, Springer theory in braid groups and the Birman-Ko-Lee monoid, \emph{Pacific J.~Math.} \textbf{205} (2002), 287--309.

\bibitem{Big2} S.~Bigelow, Braid groups are linear, \emph{J.\ Amer.\ Math.\ Soc.} \textbf{14} (2001), 471--486.

\bibitem{Bi1} J.~S.~Birman, On braid groups, \emph{Comm.\ Pure Appl.\ Math.} \textbf{22} (1969), 41--72.

\bibitem{Bi1a} J.~S.~Birman, Mapping class groups and their relationship to braid groups, \emph{Comm.\ Pure Appl.\ Math.} \textbf{22} (1969), 213--238.

\bibitem{Bi2} J.~S.~Birman, Braids, links and mapping class groups, \emph{Ann.\ Math.\ Stud.} \textbf{82}, Princeton University Press, 1974.

\bibitem{Bi3} J.~S.~Birman, Mapping class groups of surfaces, in Braids (Santa Cruz, CA, 1986), 13--43, \emph{Contemp.\ Math.} \textbf{78}, Amer.\ Math.\ Soc., Providence, RI, 1988. 

\bibitem{BB} J.~S.~Birman and T.~E.~Brendle, Braids: a survey, in Handbook of knot theory, 19--103, edited by W.~Menasco and M.~Thistlethwaite, Elsevier B.~V., Amsterdam, 2005.

\bibitem{BCP} C.-F.~B\"odigheimer, F.~R.~Cohen and M.~D.~Peim, Mapping class groups and function spaces, Homotopy methods in algebraic topology (Boulder, CO, 1999), \emph{Contemp.\ Math.} \textbf{271}, 17--39, Amer.\ Math.\ Soc., Providence, RI, 2001.

\bibitem{BRW} S.~Boyer, D.~Rolfsen and B.~Wiest, Orderable $3$-manifold groups, \emph{Ann.\ Inst.\ Fourier} \textbf{55} (2005), 243--288.

\bibitem{Bri} E.~Brieskorn, Sur les groupes de tresses (d'apr\`es V.~I.~Arnol'd), S\'eminaire Bourbaki, 24\`eme ann\'ee 
(1971/1972), Exp.\ No.\ 401, Lecture Notes in Mathematics \textbf{317}, Springer, Berlin, 1973, 21--44. 

\bibitem{BrS} E.~Brieskorn and K.~Saito, Artin-Gruppen und Coxeter-Gruppen, \emph{Invent.\ Math.} \textbf{17} (1972), 245--271.

\bibitem{Br} K.~S.~Brown, Cohomology of groups, Graduate Texts in Mathematics \textbf{87}, Springer-Verlag, New York-Berlin (1982).

\bibitem{BCG} E.~Bujalance, F.~J.~Cirre and J.~M.~Gamboa, Automorphism groups of the real projective plane with holes and their conjugacy classes within its mapping class group, \emph{Math.\ Ann.} \textbf{332} (2005), 253--275.

\bibitem{BZ} G.~Burde and H.~Zieschang, Knots, Second edition, de Gruyter Studies in Mathematics, \textbf{5}, Walter de Gruyter \& Co., Berlin, 2003. 

\bibitem{C} D.~W.~Carter, Lower $K$-theory of finite groups, \emph{Comm.\ Algebra} \textbf{8} (1980), 1927--1937.

\bibitem{Ch} W.-L.~Chow, On the algebraical braid group, \emph{Ann.\ Math.} \textbf{49} (1948) 654--658.

\bibitem{Coh} F.~R.~Cohen, Introduction to configuration spaces and their applications, in~\cite{BCHWW}, 183--261.

\bibitem{CG} F.~R.~Cohen and S.~Gitler, On loop spaces of configuration spaces, \emph{Trans.\ Amer.\ Math.\ Soc.} \textbf{354} (2002), 1705--1748.

\bibitem{Cox} H.~S.~M.~Coxeter, Regular complex polytopes, Second edition, Cambridge University Press, Cambridge, 1991.

\bibitem{CM} H.~S.~M.~Coxeter and W.~O.~J.~Moser, Generators and relations for discrete groups, Ergebnisse der Mathematik und ihrer Grenzgebiete, Vol.~14, Fourth edition, Springer-Verlag, Berlin, 1980.

\bibitem{Cr} J.~Crisp, Injective maps between Artin groups, in Geometric group theory down under, 1996, 119--137, Eds.\ J.~Cossey, C.~F.~Miller~III, W.~D.~Neumann, M.~Shapiro, de Gruyter, 1999.

\bibitem{DKR} J.~Davis, K.~Khan and A.~Ranicki, Algebraic $K$-theory over the infinite dihedral group: an algebraic approach, \emph{Algebr.\ Geom.\ Topol.} \textbf{11} (2011), 2391--2436.

\bibitem{DL} J.~Davis and W.~L\"uck, Spaces over a category and assembly maps in isomorphism conjectures in $K$- and $L$-theory, \emph{$K$-Theory} \textbf{15} (1998), 201--252.

\bibitem{Deh} P.~Dehornoy, Braid groups and left distributive operations, \emph{Trans.\ Amer.\ Math.\ Soc.} \textbf{345} (1994), 115--150.


\bibitem{DDGKM} P.~Dehornoy, F.~Digne, E.~Godelle, D.~Krammer and J.~Michel, Garside Theory, in preparation.

\bibitem{DDRW} P.~Dehornoy, I.~Dynnikov, D.~Rolfsen and B.~Wiest, Ordering braids, Mathematical Surveys and Monographs \textbf{148}, American Mathematical Society, Providence, RI, 2008.

\bibitem{Del} P.~Deligne, Les immeubles des groupes de tresses g\'en\'eralis\'es, \emph{Invent.\ Math.} \textbf{17} (1972), 273--302.

\bibitem{Dy} J.~L.~Dyer, The algebraic braid groups are torsion-free: an algebraic proof, \emph{Math.\ Z.} \textbf{172} (1980), 157--160.

\bibitem{Ep} D.~B.~A.~Epstein, Ends, in Topology of $3$-manifolds and related topics (Proc.\ Univ.\ Georgia Institute, 1961), 110--117, Prentice-Hall, Englewood Cliffs, N.J., 1962.

\bibitem{Fa} E.~Fadell, Homotopy groups of configuration spaces and the string problem of Dirac, \emph{Duke Math.~J.} \textbf{29} (1962), 231--242.

\bibitem{FH1} E.~Fadell and S.~Y.~Husseini,  Geometry and topology of configuration spaces, Springer Monographs in Mathematics, Springer-Verlag, Berlin, 2001.

\bibitem{FaN} E.~Fadell and L.~Neuwirth, Configuration spaces, \emph{Math.\ Scand.} \textbf{10} (1962), 111--118.

\bibitem{FvB} E.~Fadell and  J.~Van~Buskirk, The braid groups of $\mathbb{E}^2$ and $\St$, \emph{Duke Math.~J.} \textbf{29} (1962), 243--257.

\bibitem{FaRa} M.~Falk and R.~Randell,  The lower central series of a fiber-type arrangement, \emph{Invent.\ Math.} \textbf{82} (1985), 77--88.

\bibitem{FaRa2} M.~Falk and R.~Randell, Pure braid groups and products of free groups, in Braids (Santa Cruz, CA, 1986), 217--228, \emph{Contemp.\ Math.} \textbf{78}, Amer.\ Math.\ Soc., Providence, RI, 1988.

\bibitem{FM} B.~Farb and D.~Margalit, A primer on mapping class groups, Princeton Mathematical Series \textbf{49}, Princeton University Press, Princeton, NJ, 2012.

\bibitem{F} D.~Farley, Constructions of $EVC$ and $EFBC$ for groups acting on $CAT(0)$ spaces, \emph{Algebr.\ Geom.\ Topol.} \textbf{10} (2010),  2229--2250.

\bibitem{Farr} F.~T.~Farrell, The nonfiniteness of Nil, \emph{Proc.\ Amer.\ Math.\ Soc.} \textbf{65} (1977), 215--216. 

\bibitem{FHs} F.~T.~Farrell and W.~C.~Hsiang, The Whitehead group of poly-(finite or cyclic) groups, \emph{J.~London Math.\ Soc.} \textbf{24} (1981), 308--324.

\bibitem{FJ0}  F.~T.~Farrell and L.~E.~Jones,  Isomorphism conjectures in algebraic $K$-theory, \emph{J.~Amer.\ Math.\ Soc.} \textbf{6} (1993), 249--297.

\bibitem{FJ}  F.~T.~Farrell and L.~E.~Jones, The lower algebraic $K$-theory of virtually infinite cyclic groups, \emph{$K$-Theory} \textbf{9} (1995), 13--30.

\bibitem{FR} F.~T.~Farrell and S.~K.~Roushon, The Whitehead groups of braid groups vanish, \emph{Internat.~Math.~Res.~Notices} \textbf{10} (2000), 515--526.

\bibitem{FZ} E.~M.~Feichtner and G.~M.~Ziegler, The integral cohomology algebras of ordered configuration spaces of spheres, \emph{Doc.\ Math.} \textbf{5} (2000), 115--139.

\bibitem{FGRW} R.~A.~Fenn, M.~T.~Greene, D.~Rolfsen, C.~Rourke and B.~Wiest, Ordering the braid groups, \emph{Pac.\ J.\ Math.} \textbf{191} (1999), 49--74.

\bibitem{FoN} R.~H.~Fox and L.~Neuwirth, The braid groups, \emph{Math.\ Scand.} \textbf{10} (1962), 119--126.

\bibitem{Ga} F.~A.~Garside, The braid group and other groups, \emph{Quart.\ J.~Math.\ Oxford} \textbf{20} (1969), 235--254.

\bibitem{GVB} R.~Gillette and J.~Van Buskirk, The word problem and consequences for the braid groups and mapping class groups of the $2$-sphere, \emph{Trans.\ Amer.\ Math.\ Soc.} \textbf{131} (1968), 277--296. 

\bibitem{Go} C.~H.~Goldberg, An exact sequence of braid groups, \emph{Math.\ Scand.} \textbf{33} (1973), 69--82.

\bibitem{GG1} D.~L.~Gon\c{c}alves and J.~Guaschi, On the structure of surface pure braid groups, \emph{J.~Pure Appl.\ Algebra} \textbf{182} (2003), 33--64 (due to a printer's error, this article was republished in its entirety with the reference \textbf{186} (2004) 187--218).

\bibitem{GG2} D.~L.~Gon\c{c}alves and J.~Guaschi, The roots of the full twist for surface braid groups, \emph{Math.\ Proc.\ Camb.\ Phil.\ Soc.} \textbf{137} (2004), 307--320.

\bibitem{GG3} D.~L.~Gon\c{c}alves and J.~Guaschi, The braid groups of the projective plane, \emph{Algebr.\ Geom.\ Topol.} \textbf{4} (2004), 757--780.

\bibitem{GG4} D.~L.~Gon\c{c}alves and J.~Guaschi, The braid group $B_{n,m}(\St)$ and the generalised Fadell-Neuwirth short exact sequence, \emph{J.~Knot Theory Ramif.} \textbf{14} (2005), 375--403.

\bibitem{GG5} D.~L.~Gon\c{c}alves and J.~Guaschi, The quaternion group as a subgroup of the sphere braid groups, \emph{Bull.\ London Math.\ Soc.} \textbf{39} (2007), 232--234.

\bibitem{GGgeom} D.~L.~Gon\c{c}alves and J.~Guaschi, The braid groups of the projective plane and the Fadell-Neuwirth short exact sequence, \emph{Geom.\ Dedicata} \textbf{130} (2007), 93--107.

\bibitem{GGagt2} D.~L.~Gon\c{c}alves and J.~Guaschi,  The classification and the conjugacy classes of the finite subgroups of the sphere braid groups, \emph{Algebr.\ Geom.\ Topol.} \textbf{8} (2008), 757--785.

\bibitem{GG6} D.~L.~Gon\c{c}alves and J.~Guaschi, The lower central and derived series of the braid groups of the sphere, \emph{Trans.\ Amer.\ Math.\ Soc.} \textbf{361} (2009), 3375--3399.

\bibitem{GGjktr2} D.~L.~Gon\c{c}alves and J.~Guaschi, The lower central and derived series of the braid groups of the finitely-punctured sphere, \emph{J.~Knot Theory Ramif.} \textbf{18} (2009), 651--704.

\bibitem{GG7} D.~L.~Gon\c{c}alves and J.~Guaschi, Braid groups of non-orientable surfaces and the Fadell-Neuwirth short exact sequence, \emph{J.~Pure Appl.\ Algebra} \textbf{214} (2010), 667--677.

\bibitem{GG8} D.~L.~Gon\c{c}alves and J.~Guaschi, Classification of the virtually cyclic subgroups of the pure braid groups of the projective plane, \emph{J. Group Theory} \textbf{13} (2010), 277--294.

\bibitem{GG9} D.~L.~Gon\c{c}alves and J.~Guaschi, The Borsuk-Ulam theorem for maps into a surface, \emph{Topology Appl.} \textbf{157} (2010), 1742--1759.

\bibitem{GG10} D.~L.~Gon\c{c}alves and J.~Guaschi, The lower central and derived series of the braid groups of the projective plane, \emph{J.~Algebra} \textbf{331} (2011), 96--129.

\bibitem{GG11} D.~L.~Gon\c{c}alves and J.~Guaschi, Surface braid groups and coverings, \emph{J.~London Math.\ Soc.} \textbf{85} (2012), 855--868.

\bibitem{GG13} D.~L.~Gon\c{c}alves and J.~Guaschi, Minimal generating and normally generating sets for the braid and mapping class groups of the disc, the sphere and the projective plane, \emph{Math.~Z.}, to appear.

\bibitem{GG12} D.~L.~Gon\c{c}alves and J.~Guaschi, The classification of the virtually cyclic subgroups of the sphere braid groups, monograph to appear in the series SpringerBriefs in Mathematics, \url{arXiv math:1110.6628}.

\bibitem{GG15} D.~L.~Gon\c{c}alves and J.~Guaschi, Some homotopy properties of the inclusion $F_n(S) \hooklongrightarrow S^n$ for S either $\St$ or $\rp$ and the virtual cohomological dimension of $B_n(S)$ and $P_n(S)$, work in progress.

\bibitem{GG14} D.~L.~Gon\c{c}alves and J.~Guaschi, Conjugacy classes of finite subgroups of the braid groups of the projective plane, work in progress.

\bibitem{GG16} D.~L.~Gon\c{c}alves and J.~Guaschi, The classification of the virtually cyclic subgroups of the braid groups of the projective plane, work in progress.

\bibitem{GM1} J.~Gonz\'alez-Meneses, New presentations of surface braid groups, \emph{J.~Knot Theory Ramif.} \textbf{10} (2001), 431--451. 

\bibitem{GM2} J.~Gonz\'alez-Meneses, Ordering pure braid groups on closed surfaces,  \emph{Pac.\ J.\ Math.} \textbf{203} (2002), 369--378.

\bibitem{GM4} J.~Gonz\'alez-Meneses, Basic results on braid groups, \emph{Ann.\ Math.\ Blaise Pascal} \textbf{18} (2011), 15--59.

\bibitem{GMP} J.~Gonz\'alez-Meneses and L.~Paris, Vassiliev invariants for braids on surfaces, \emph{Trans.\ Amer.\ Math.\ Soc.}  \textbf{356} (2004), 219--243.

\bibitem{GW} J.~González-Meneses and B.~Wiest, On the structure of the centralizer of a braid, \emph{Ann.\ Sci.\ \'Ecole Norm.\ Sup.} \textbf{37} (2004), 729--757.

\bibitem{GL} E.~A.~Gorin and V.~J.~Lin, Algebraic equations with continuous coefficients and some problems of the algebraic theory of braids, \emph{Math.\ USSR Sbornik} \textbf{7} (1969), 569--596.

\bibitem{GJM} J.~Guaschi, D.~Juan-Pineda and S.~Mill\'an-L\'opez, The lower algebraic $K$-theory of the braid groups of the sphere, preprint, \url{arXiv:1209.4791}.

\bibitem{GJM2} J.~Guaschi, D.~Juan-Pineda and S.~Mill\'an-L\'opez, The lower algebraic $K$-theory of the finite subgroups of the braid groups of the sphere, work in progress.

\bibitem{Ham0} M.-E.~Hamstrom, Homotopy properties of the space of homeomorphisms on $P^{2}$ and the Klein bottle, \emph{Trans.\ Amer.\ Math.\ Soc.} \textbf{120} (1965), 37--45.

\bibitem{Ham} M.-E.~Hamstrom, Homotopy groups of the space of homeomorphisms on a $2$-manifold, \emph{Illinois J.~Math.} \textbf{10} (1966), 563--573.

\bibitem{Ha} V.~L.~Hansen, Braids and Coverings: selected topics, \emph{London Math.\ Society Student Text} \textbf{18}, Cambridge University Press, 1989.

\bibitem{Hat} A.~Hatcher, Algebraic topology, Cambridge University Press, Cambridge, 2002.

\bibitem{Ho} L.~Hodgkin, $K$-theory of mapping class groups: general $p$-adic $K$-theory for punctured spheres, \emph{Math.~Z.} \textbf{218} (1995), 611--634.

\bibitem{Hu} S.~T.~Hu, Homotopy theory, Pure and Applied Mathematics, Vol.~VIII, Academic Press, New York, 1959.

\bibitem{I} N.~V.~Ivanov, Mapping class groups, in Handbook of geometric topology, 523--633, North-Holland, Amsterdam, 2002.

\bibitem{J} D.~L.~Johnson, Presentation of groups, LMS Lecture Notes \textbf{22} (1976), Cambridge University Press.

\bibitem{Jon1} V.~F.~R.~Jones, Braid groups, Hecke algebras and type $\text{II}_1$ factors, in Geometric methods in operator algebras (Kyoto, 1983), Pitman Res.\ Notes Math.\ Ser., \textbf{123}, 242--273, Longman Sci.\ Tech., Harlow, 1986.

\bibitem{Jon2} V.~F.~R.~Jones, Hecke algebra representation of braid groups and link polynomials, \emph{Ann.\ Math.} \textbf{126} (1987), 335--388. 

\bibitem{JLMP} D.~Juan-Pineda,  J.-F.~Lafont, S.~Mill\'an-Vossler and S.~Pallekonda, Algebraic $K$-theory of virtually free groups, \emph{Proc.\ Roy.\ Soc.\ Edinburgh Sect.~A} \textbf{141} (2011), 1295--1316.

\bibitem{JL} D.~Juan-Pineda and I.~Leary, On classifying spaces for the family of virtually cyclic subgroups, in Recent developments in algebraic topology, \emph{Contemp.\ Math.} \textbf{407} (2001), 135--145.

\bibitem{JM} D.~Juan-Pineda and S. Mill\'an-L\'opez, Invariants associated to the pure braid groups of the sphere, \emph{Bol.\ Soc.\ Mat.\ Mexicana} \textbf{3} (2006), 27--32.

\bibitem{JM2} D.~Juan-Pineda and S.~Mill\'an-L\'opez, The Whitehead group and the lower algebraic $K$-theory of braid groups of $\St$ and $\rp$, \emph{Algebr.\ Geom.\ Topol.} \textbf{10} (2010), 1887--1903.

\bibitem{JS} D.~Juan-Pineda and J.~S\'anchez, The $K$-theoretic Farrell-Jones Isomorphism conjecture for braid groups, submitted for publication.

\bibitem{Ka} C.~Kassel, L'ordre de Dehornoy sur les tresses, in S\'eminaire Bourbaki, Vol.\ 1999/2000, \emph{Ast\'erisque} \textbf{276} (2002), 7--28.

\bibitem{KT} C.~Kassel and V.~Turaev, Braid groups, Graduate Texts in Mathematics \textbf{247}, Springer, New York, 2008.

\bibitem{KP} R.~P.~Kent~IV and D.~Peifer, A geometric and algebraic description of annular braid groups, \emph{Int.~J.~Algebra and Computation} \textbf{12} (2002) 85--97.

\bibitem{KR} D.~M.~Kim and D.~Rolfsen, An ordering for groups of pure braids and fibre-type hyperplane arrangements, \emph{Canad.\ J.\ Math.} \textbf{55} (2003), 822--838.

\bibitem{Kor} M.~Korkmaz, On the linearity of certain mapping class groups, \emph{Turkish J.\ Math.} \textbf{24} (2000), 367--371.

\bibitem{Kr1} D.~Krammer, The braid group $B_4$ is linear, \emph{Invent.\ Math.} \textbf{142} (2000), 451--486.

\bibitem{Kr2} D.~Krammer, Braid groups are linear, \emph{Ann.\ Math.} \textbf{155} (2002), 131--156.

\bibitem{Ku} A.~Kuku, Higher algebraic $K$-theory, Handbook of Algebra~\textbf{4}, 3--74, Elsevier/North-Holland, Amsterdam, 2006.

\bibitem{Lad} Y.~Ladegaillerie, Groupes de tresses et probl\`eme des mots dans les groupes de tresses, \emph{Bull.\ Sci.\ Math. (2)} \textbf{100} (1976), 255--267.
 
\bibitem{Lam} S.~Lambropoulou,  Braid structures in knot complements, handlebodies and $3$-manifolds, in Knots in Hellas~'98 (Delphi), 274--289, Ser.\ Knots Everything, \textbf{24}, World Sci.\ Publishing, River Edge, NJ, 2000.

\bibitem{L} W.~L\"uck, On the classifying space of the family of virtually cyclic subgroups for $\operatorname{CAT}(0)$-groups. \emph{M\"unster J.~Math.} \textbf{2} (2009), 201--214.

\bibitem{Mag0} W.~Magnus, \"Uber Automorphismen von Fundamentalgruppen berandeter Fl\"achen, \emph{Math.\ Ann.} \textbf{109} (1934), 617--646.

\bibitem{Mag} W.~Magnus, Braid groups: a survey, in Proceedings of the Second International Conference on the Theory of Groups (Australian Nat.\ Univ., Canberra, 1973), 463--487, Lecture Notes in Maths., Vol.~372, Springer, Berlin, 1974.

\bibitem{MKS} W.~Magnus, A.~Karrass and D.~Solitar, Combinatorial group theory, Second revised edition, Dover Publications Inc., New York, 1976.

\bibitem{Man} S.~Manfredini, Some subgroups of Artin's braid group, \emph{Topology Appl.} \textbf{78} (1997), 123--142.

\bibitem{McC} G.~S.~McCarty Jr., Homeotopy groups, \emph{Trans.\ Amer.\ Math.\ Soc.} \textbf{106} (1963), 293--304.

\bibitem{Mi}  S.~Mill\'an-Vossler, The lower algebraic $K$-theory of braid groups on $\St$ and $\rp$, VDM Verlag Dr.~Muller, Germany, 2008.

\bibitem{Mo} S.~Moran, The mathematical theory of knots and braids, an introduction, North-Holland Mathematics Studies \textbf{82}, North-Holland Publishing Co., Amsterdam, 1983.

\bibitem{M} K.~Murasugi, Seifert fibre spaces and braid groups, \emph{Proc.\ London Math.\ Soc.} \textbf{44} (1982), 71--84.

\bibitem{MK} K.~Murasugi and B.~I.~Kurpita,  A study of braids, Mathematics and its Applications \textbf{484}, Kluwer Academic Publishers, Dordrecht, 1999.

\bibitem{Ne} M.~H.~A.~Newman, On a string problem of Dirac, \emph{J.\ London Math.\ Soc.} \textbf{17} (1942), 173--177. 

\bibitem{Oc} O.~Ocampo Uribe, Grupos de tran\c{c}as brunnianas, Ph.D thesis, Universidade de S\~ao Paulo, Brazil, work in progress.

\bibitem{O} R.~Oliver, Whitehead groups of finite groups, London Mathematical Society Lecture Note Series \textbf{132}, Cambridge University Press, Cambridge, 1988.

\bibitem{P} L.~Paris, Braid groups and Artin groups, in Handbook of Teichm\"uller theory Vol.~II, 389--451, \emph{IRMA Lect.\ Math.\ Theor.\ Phys.} \textbf{13}, Eur.\ Math.\ Soc., Z\"urich, 2009.

\bibitem{PR} L.~Paris and D.~Rolfsen, Geometric subgroups of surface braid groups \emph{Ann.\ Inst.\ Fourier} \textbf{49} (1999), 417--472.

\bibitem{PW} E.~Pedersen and C.~Weibel, A non-connective delooping of algebraic $K$-theory, Springer Lecture Notes in Math.\ \textbf{1126}, Springer-Verlag, Berlin-Heidelberg-New York, 166--181, 1985.

\bibitem{Q1} D.~Quillen, Higher algebraic $K$-theory: I, Cohomology of groups and algebraic K-theory, 413--478, \emph{Adv.\ Lect.\ Math.} \textbf{12}, Int.\ Press, Somerville, MA, 2010.

\bibitem {Q0} D.~Quillen, Finite generation of the groups $K_i$ of rings of algebraic integers, Cohomology of groups and algebraic $K$-theory, 479--488, \emph{Adv.\ Lect.\ Math.}, \textbf{12}, Int.\ Press, Somerville, MA, 2010.

\bibitem{Ra} R.~Ramos, Non-finiteness of twisted nils, \emph{Bol.\ Soc.\ Mat.\ Mexicana} \textbf{13} (2007), 55--64. 

\bibitem{R} D.~Rolfsen, Tutorial on the braid groups, in~\cite{BCHWW}, 1--30.

\bibitem{RolW} D.~Rolfsen and B.~Wiest, Free group automorphisms, invariant orderings and topological applications, \emph{Algebr.\ Geom.\ Topol.} \textbf{1} (2001), 311--320.

\bibitem{RZ} D.~Rolfsen and J.~Zhu, Braids, orderings and zero divisors, \emph{J.~Knot Theory Ramif.} \textbf{7} (1998), 837--841.

\bibitem{RouW} C.~Rourke and B.~Wiest, Order automatic mapping class groups, \emph{Pacific J.\ Math.} \textbf{194} (2000), 209--227.

\bibitem{S} G.~P.~Scott, Braid groups and the group of homeomorphisms of a surface, \emph{Proc.\ Camb.\ Phil.\ Soc.} \textbf{68} (1970), 605--617.

\bibitem{SW} H.~Short and B.~Wiest, Orderings of mapping class groups after Thurston, \emph{Enseign.\ Math.} \textbf{46} (2000), 279--312.

\bibitem{St} M.~Stukow, Conjugacy classes of finite subgroups of certain mapping class groups, Seifert fibre spaces and braid groups, \emph{Turkish J.\ Math.} \textbf{2} (2004), 101--110.

\bibitem{Sw} R.~G.~Swan, Projective modules over binary polyhedral groups, \emph{J.~Reine Angew.\ Math.} \textbf{342} (1983), 66--172.

\bibitem{Th} J.~G.~Thompson, Note on $H(4)$, \emph{Comm.\ Algebra} \textbf{22} (1994), 5683--5687.

\bibitem{tD} T.~tom~Dieck, Transformation groups and representation theory, Springer Lecture Notes in Mathematics \textbf{766}, Springer, Berlin, 1979.

\bibitem{vB} J.~Van~Buskirk, Braid groups of compact $2$-manifolds with elements of finite order, \emph{Trans.\ Amer.\ Math.\ Soc.} \textbf{122} (1966), 81--97.

\bibitem{Ve} V.~V.~Vershinin, Braids, their properties and generalizations, Handbook of Algebra~\textbf{4}, Elsevier/North-Holland, Amsterdam, 2006, 427--465.

\bibitem{Wa} F.~Waldhausen, Algebraic $K$-theory of generalized free products, Part I, \emph{Ann. Math.} \textbf{108} (1978), 135-–204.

\bibitem{W} C.~T.~C.~Wall, Poincar\'e complexes~I, \emph{Ann.\ Math.} \textbf{86} (1967), 213--245.

\bibitem{We} C.~Wegner, The $K$-theoretic Farrell-Jones conjecture for $CAT(0)$-groups, \emph{Proc.\ Amer.\ Math.\ Soc.} \textbf{140} (2012), 779--793.

\bibitem{Weib} C.~Weibel, $NK_0$ and $NK_1$ of the groups $C_4$ and $D_4$, \emph{Comment.\ Math.\ Helv.} \textbf{8} (2009), 339--349 (addendum).

\bibitem{Wh1} G.~W.~Whitehead, Elements of homotopy theory, Graduate Texts in Mathematics \textbf{61}, Springer-Verlag, New York, 1978.

\bibitem{Wo} J.~A.~Wolf, Spaces of constant curvature, sixth edition, AMS Chelsea Publishing, vol.~372, 2011.

\bibitem{Z1} O.~Zariski,  On the Poincar\'e group of rational plane curves, \emph{Amer.\ J.\ Math.} \textbf{58} (1936), 607--619.

\bibitem{Z2} O.~Zariski, The topological discriminant group of a Riemann surface of genus $p$, \emph{Amer.\ J.\ Math.} \textbf{59} (1937), 335--358.

}

\end{thebibliography}
\end{document}